\documentclass[11pt]{article}

\usepackage{lineno,hyperref}
\usepackage{color}
\usepackage{here}
\usepackage{listings}
\usepackage{ifthen}
\usepackage{multirow}
\usepackage{array}
\usepackage{alterqcm}
\usepackage{pifont}
\usepackage{exscale, relsize}
\usepackage{soul}
\usepackage{makeidx}
\usepackage{graphicx}
\usepackage{graphics}
\usepackage{newlfont}
\usepackage{t1enc}
\usepackage[french2]{minitoc}
\usepackage{exscale,relsize}
\usepackage{subfig}
\usepackage{amsmath}
\usepackage{soul}
\usepackage[affil-it]{authblk}

\newtheorem{lem}{Lemma}[section]
\newtheorem{thm}{Theorem}[section]

\newtheorem{rmq}{Remark}[section]

\newtheorem{note}{Note}[section]
\begin{document}

\title{Tree-Grass interactions dynamics and Pulse Fires: mathematical and numerical studies}

\author{A. Tchuint\'{e} Tamen$^{\star,+}$, Y. Dumont $^{\dagger}$\footnote{Corresponding author: yves.dumont@cirad.fr},  S. Bowong$^{\ddagger,\star,+}$   , J. J. Tewa$^{\star, +}$,  P. Couteron$^{\ddagger\ddagger}$}
  \affil{
  $^{\star}$ LIRIMA, GRIMCAPE, Faculty of Science, University of Yaounde 1, Cameroon \\
  $^{\dagger}$ CIRAD, Umr AMAP, Montpellier, France \\
  $^{\ddagger}$ University of Douala, Cameroon \\
  $^{+}$ IRD, UMI 209, UMMISCO, IRD France Nord, Bondy, France \\
  $^{\ddagger\ddagger}$ IRD, Umr AMAP, Montpellier, France}

\maketitle
\begin{abstract}

Savannas are dynamical systems where grasses and trees can either dominate or coexist. Fires are known to be central in the functioning of the savanna biome though their characteristics are expected to vary along the rainfall gradients as observed in Sub-Saharan Africa. In this paper, we model the tree-grass dynamics using impulsive differential equations that consider fires as discrete events. This framework allows us to carry out a comprehensive qualitative mathematical analysis that revealed more diverse possible outcomes than the analogous continuous model. We investigated local and global properties of the equilibria and show that various states exist for the physiognomy of vegetation. Though several abrupt shifts between vegetation states appeared determined by fire periodicity, we showed that direct shading of grasses by trees is also an influential process embodied in the model by a competition parameter leading to bifurcations. Relying on a suitable nonstandard finite difference scheme, we carried out numerical simulations in reference to three main climatic zones as observable in Central Africa.
\end{abstract}






\
\section{Introduction}

In savannas, trees and grasses typically coexist 
\cite{Menaut1983}. Fire is  recognized as playing a major part in the dynamics of this biome. The nature of grass-tree interactions
and fire regimes strongly vary along environmental gradients in
tropical savannas. Fire is more intense in wet than in arid
savannas, where lower water availability leads to lower grass, i.e.
fuel load, production. Thus fire is expected to control Tree-Grass
dynamics in wet savannas \cite{Frost1986}. But two
hypotheses for Tree-Grass coexistence have been introduced during
these last decades. First, Walter (1971) \cite{Walter1971ecology}
proposed the idea that trees and grass exploit two different rooting
niches. Grasses are rooted in superficial soil layers and first use the
incoming water, whereas tree roots are situated in subsoil, so
that trees could grow only where enough water reached deeper soil
horizons. This idea was developed analytically by Walker and
Noy-Meir (1982) \cite{Walker1982aspects} using a Lotka-Volterra
theory of co-existence between competitors. The second hypothesis says that grass-tree coexistence
	is driven by limited opportunities for seedling to escape both
	droughts and flame zone into the adult stage (Hochberg et al. 1994 \cite{Hochberg1994influences}; Higgins et al. (2000)
	\cite{Higgins2000fire}). In areas where tree seedlings succeed to establish in spite of competition with grasses, they are burnt by frequent grass fires (Higgins et al. 2000
\cite{Higgins2000fire}).\par

Savannas fires are frequent, up to occurring every 1-5 years in
wet savannas (Frost \& Robertson 1985) though the fire return time is usually a decreasing function of mean annual precipitation. Fuel load made of dead aerial grass parts typically ranges
between $2$ and $10$ t.ha$^{-1}$ of dry matter (DM) (Lacey et al. 1982
\cite{Lacey1982fire}; Stronach \& Mac-Naughton 1989
\cite{Stronach1989grassland}; Menaut et al. 1991
\cite{Menaut1991biomass}; Mordelet 1993
\cite{Mordelet1993influence}) and flame height is usually $2$-$3$
metres high (Frost \& Robertson 1987 \cite{Frost1987fire}). Although
the fire burns most or all the aboveground grass biomass, the
large underground root systems of perennial grass species enable most of the tufts to survive even the most intense fires and to rapidly establish new shoots before the onset of the rainy season. In contrast to grasses, trees which are less than $2$m height may either succumb to fire or have to resprout from roots and have their growth delayed (Bond and Midgley, 2001
\cite{Bond2001ecology}). Mature trees ($>8$m) and shrubs beyond $2$m
are more fire resistant and only experience partial die-back (Menaut \&
César 1979 \cite{Menaut1979structure}; Gillon 1983
\cite{Gillon1983fire}). Early fires (in the beginning of the dry
season) are less violent than late fires and have a lower impact on
tree regeneration (Abbadie et al. 2006 \cite{Abbadie2006lamto}).\par

Africa is a land of extreme contrasts in rainfall distribution and
the time of year during which rainfall occurs (Janowiak 1987
\cite{Janowiak1988}). When soil resource supply is
temporally variable, trees and grasses will experience two distinct
phases of resource availability: pulse periods when resources are
high and most growth and biomass accumulation (fuel load) occurs,
and inter-pulse periods when resources are too low for most tree and
grass to take up and most mortality due to resource deficits takes
place (Goldberg \& Novoplansky, 1997 \cite{Goldberg1997}; Noy-Meir,
1973 \cite{Noy1973desert}). Hence essential resource availability
(e.g. water) is discontinuously available and the availability of
these resources impact the ecosystem as discrete pulse events
interspersed among long periods of limited resource availability
(Schwinning et al. 2004 \cite{Schwinning2004}).\par

Fires are sudden event
that consume trees and grass biomass (Scheiter 2008
\cite{Scheiter2009Grasstree}).
The broad objective of this 
study is to examine the influence of pulse events with regard to
 fires impact on the Tree-Grass dynamics along the rainfall
gradient in Africa.  Tree-grass
savanna models can not be studied without the important role of fires
(Tilman 1994 \cite{Tilman1994}; Higgins et al. 2000
\cite{Higgins2008physically}; Sankaran et al. 2004
\cite{Sankaran2004tree}, 2005 \cite{Sankaran2005determinants};
D'Odorico et al. 2006 \cite{DOdorico2006probabilistic}; Accatino et
al. 2010 \cite{Accatino2010tree}; Beckage et al. 2011
\cite{Beckage2011grass}; Staver et al. 2011 \cite{Staver2011tree}; Yatat et al. 2014 \cite{yatat2014} and Tchuinte
et al. 2014 \cite{Tchuinte2014}).
This paper extends
our earlier work (Tchuinte et al. 2014 \cite{Tchuinte2014}) where we consider a continuous tree-grass interaction model that featured a fairly generic family of non-linear functions of grass biomass to model fire intensity and its impact on tree. We have shown
that the continuous model is able to predict a variety of dynamical outcomes.
Notably, the number of equilibria featuring Tree-Grass coexistence
depends on the characteristics of fairly generic Monod functions used to
model the fire impact on tree dynamics. Moreover, we have shown that
various bistability situations occur among forest, grassland and
Tree-Grass (i.e. savanna) equilibria (for more detail see Tchuinte et al. 2014
\cite{Tchuinte2014}). Of course, in practice, fires are not
continuous. Recent studies of the
interactions between fire and vegetation are based on stochastic
approaches because of the random and unpredictable nature of fire
occurrences (D'Odorico et al., 2006
\cite{DOdorico2006probabilistic}; Beckage et al. 2011
\cite{Beckage2011grass}).\par

In section 2 we will present the model with pulse fires. The theoretical analysis is developed in  section
3. We show that the system admits
four equilibria among which two trivial equilibria (the bare soil
and the forest equilibria), and two periodic equilibria (the
periodic grassland and the periodic savanna equilibria). We show
that there are various bistabilities: between forest and grassland;
between forest and savanna. Local and global stabilities are distinguished
using classical tools such as Floquet multipliers and comparison
theorem. We highlight thresholds that summarize the dynamics of the model and explain the theoretical meaning of these thresholds.
Prior to illustrate our theoretical results numerically, in section
4, based on the scheme developed in \cite{yatat2014}, we develop a reliable nonstandard finite difference method (NSFD)
that preserves the qualitative properties of the system
(Anguelov et al. 2012 \cite{Anguelov2012}, 2013  \cite{Anguelov2013}, 2014 \cite{Anguelov2014}).  Section 6
concludes the paper.
Some mathematical details are included in appendices.


\section{The mathematical model}


In savanna environment, fire intensity is tightly linked to dried
grass biomass that remains during the dry season (Higgins et al.
2008 \cite{Higgins2008physically}). During the last decades the
effects of fire on vegetation dynamics have been studied (Scholes
and Walker 1993 \cite{Scholes1993african}; Higgins et al. 2000
\cite{Higgins2000fire}). Most of the models associated to or derived from these studies are ordinary
differential equations (ODE) which assume that fires occur
continuously with a fixed frequency. However fires are sudden event
that consume  grass biomasses and kill or harm tree seedlings (Scheiter 2008
\cite{Scheiter2009Grasstree}). The season of burning and the time
between recurring fires determine trees and grass physioniomies in
most ecosystems and especially in the savanna biome (Thonicke et al., (2001) \cite{Thonicke2001}). In
this paper we present a new Tree-Grass model that aim to contribute to
our understanding about  how pulse fire shapes vegetation dynamics in
fire-prone savanna-like ecosystems. We consider fire as discrete
events and derive the following impulsive differential system

 \begin{equation}
 \left\{
 \begin{array}{l}
 \left.
 \begin{array}{l}
 \displaystyle \frac{dG}{dt}=\gamma_{G}G\left(1-\displaystyle\frac{G}{K_{G}}\right)-\delta_{G0}G-\gamma_{TG}TG,\\
 \\
 \displaystyle \frac{dT}{dt}=\gamma_{T}T\left(1-\displaystyle\frac{T}{K_{T}}\right)-\delta_{T}T,\\
 \end{array}
 \right\}, t\neq t_{n}, n=1, 2, ..., N_{\tau} \\
 \left.
 \begin{array}{l}\\
 \Delta G(t_{n})=G(t_{n}^{+})-G(t_{n})=-\lambda_{fG}G(t_{n}),\\
 \\
 \Delta T(t_{n})=T(t_{n}^{+})-T(t_{n})=-\lambda_{fT}\omega(\lambda_{fG}G(t_{n}))T(t_{n}),\\
 \end{array}
 \right\}, t= t_{n}, n=1, 2, ..., N_{\tau} \\
 \\
 G(t_{0}^{+})= G_{0},\\
 \\
 T(t_{0}^{+})= T_{0},
 \end{array}
 \right.
 \label{Impuleq1}
 \end{equation}

where,\par

\begin{itemize}
\item $T$ and $G$ are tree and grass biomasses respectively,
\item
$\tau=\frac{1}{f}$ is the period of time between two consecutive
fires, and $f$ is the frequency of fire,
\item $N_{\tau}$ is a countable number of fire occurrence,
\item $t_{n}=n\tau,$ $n=1, 2, ..., N_{\tau}$
are called moments of impulsive effects of fire, and satisfy $0\leq
t_{1}<t_{2}<...<t_{k}<t_{N_{\tau}}$ ,
\item $\omega(G)$
is a generic non-linear functional  which expresses fire intensity
as an increasing function of grass biomass. Other than smoothness, it
satisfies the following three conditions: $(i)$ fires spread if
and only if fuel is available ($\omega(0)=0$), $(ii)$ fire-impact
increases with fuel available ($\omega(G)\geq 0$,
$\omega^{'}(G)>0$), and $(iii)$ there is boundary effects
$\lim\limits_{G\rightarrow\infty}\omega(G)<1$.
\end{itemize}
Other parameters used are listed in the following table.

\begin{table}[H]
{\small
\begin{center}
 \caption{ {\small Parameter symbols and names used to initialize the model}}
\renewcommand{\arraystretch}{1.5}
\begin{tabular}{cccc}
\hline\hline
Symbol & Parameter name & Units \\
 \hline\hline
 $\gamma_{G}$  & Grass biomass production per unit of grass biomass per year  &  yr$^{-1}$\\
 $\delta_{G0}$ &   Grass biomass loss by herbivory (grazing) or human action & yr$^{-1}$\\
 $K_{G}$ & Carrying capacity of grass biomass & t.ha$^{-1}$\\
 $\mu_{G}=\displaystyle\frac{\gamma_{G}}{K_{G}}$ & Additional death due to grass-grass competition  & ha.t$^{-1}$.yr$^{-1}$\\
 $\lambda_{fG}$ &  loss of grass biomass due to fire & - \\
$\gamma_{T}$& Tree biomass production per unit of tree biomass per year &  yr$^{-1}$\\
$\delta_{T}$ & Tree biomass loss by herbivory (browsing) or human action & yr$^{-1}$\\
$K_{T}$ & Carrying capacity of tree biomass  & t.ha$^{-1}$ \\
$\mu_{T}=\displaystyle\frac{\gamma_{T}}{K_{T}}$ & Additional death due to tree-tree competition  & ha.t$^{-1}$.yr$^{-1}$\\
$\lambda_{fT}$ &   loss of tree biomass due to fire & -\\
$\gamma_{TG}$ & grass mortality due to tree/grass competition  & ha.t$^{-1}$.yr$^{-1}$\\
\hline\hline
\end{tabular}
\label{Params_symbols}
\end{center}}
\end{table}

We suppose that solutions of $(\ref{Impuleq1})$ is right continuous
at $t_{n}$, $n=1, 2, ..., N_{\tau}$, that is
$G(t_{n}^{+})=\lim\limits_{h\rightarrow 0^{+}}G(t_{n}+h)=G(t_{n})$
and $T(t_{n}^{+})=\lim\limits_{h\rightarrow
0^{+}}T(t_{n}+h)=T(t_{n})$, where
 $G(t_{n}^{+})$
and $T(t_{n}^{+})$ are the biomass values for grasses and trees
instantly after impulsive fire. Immediately following each fire
pulse, system (\ref{Impuleq1}) evolves from its new initial state
without being further affected by the fire scheme until the next
pulse is applied. In agreement with empirical experience (Abbadie et al., 2006 \cite{Abbadie2006lamto}), we assume that the level of destruction of the
tree biomass depends on the available grass biomass through 
$\omega(G)$.\par

\section{Theoretical analysis}

System (\ref{Impuleq1}) belongs to basic theory of IDE (Bainov 1993
\cite{Bainov1993}) and their applications in Ecology. IDEs generally
describe phenomena which are subject to abrupt or instantaneous
changes. Model (\ref{Impuleq1}) derives from the family of impulsive
Kolmogorov-type population dynamics in the theory of mathematical
biology which the general form is given by

\begin{equation}
 \left\{
 \begin{array}{l}
\displaystyle \frac{dx}{dt}=x_{i}(t)f_{i}(t,x(t)), \hspace{0.1cm}
t\neq t_{k},\\
x_{i}(t_{k}^{+})=I_{ik}(t_{k}, x(t_{k})), \hspace{0.1cm} t=t_{k}
 \end{array}
 \right.
\label{Impuleq2}
 \end{equation}
 where $x_{i}(t)$
 represents the density or size of species $x_{i}$ at the time $t$, $x(t)=(x_{1}(t), x_{2}(t), ..., x_{n}(t))$, $f(t,x)$ is an
 n-dimensional real
 functional defined by
 \begin{displaymath}
 f(t, x)=(f_{1}(t, x), f_{2}(t, x), ..., f_{n}(t, x))
\end{displaymath}
 and $I_{k}(t_{k},x)$ is also an n-dimensional
 real functional defined by

 \begin{displaymath}
I_{k}(t_{k}, x)=(I_{1k}(t_{k}, x), I_{2k}(t_{k}, x), ...,
I_{nk}(t_{k}, x)).
\end{displaymath}

This family of models has recently attracted the attention of
several authors (\cite{Baek2010}, \cite{Bunimovich2008},
\cite{Ahmad2007}, \cite{Lakmeche2000}, \cite{He2009},
\cite{Wang2009}, \cite{Zhang2008}, and the references cited
therein). The main study subjects are the permanence, persistence
and extinction of species, the local and global asymptotic stability
of systems, the existence and uniqueness of positive periodic
solution and almost periodic solution, and the bifurcation and
dynamical complexity, etc. However, in all models investigated in
the literature,  authors generally do not consider the non-linear
impulses i.e. non-linear form of the function $I_{k}(t_{k},x)$. They
mostly focused on the quasi-linear impulses. Since for the
continuous model (see Tchuinte et al. 2014 \cite{Tchuinte2014}) the
non-linear shape brings a wealth of possibilities, notably for the
existence of various positive Tree-Grass equilibria. We retain this option although it could render
the model difficult to study. In our model, we consider a generic
non-linear functional response $\omega(G)$, which expresses the
causality between grass biomass and fire intensity as to model the
impact of fire on the woody biomass.  Depending of the fuel accumulation (grass
biomass  i.e. G), $\omega(G)$ could take  the general sigmoidal form
$\displaystyle\frac{G^{\theta}}{\alpha^{\theta}+G^{\theta}}$,
$\theta>0$, or some equivalent form. However, in our study
$\omega(G)$  is principally treated as a non-linear increasing function. In this regard, since
our Tree-Grass model does not contain more than two populations in
competition, it may appear to be simple mathematically at first
sight, but it is, in fact very challenging and complicated due to the
nonlinear impulsive functions.\par

Prior to analyzing the above model $(\ref{Impuleq1})$, it is
important to show positivity and boundedness for solutions as they
represent biomasses. Positivity implies that the populations survive
and boundedness may be interpreted as a natural restriction to
growth as a consequence of limited resources. Then, model
$(\ref{Impuleq1})$ requires that trajectories remain positive and
that trajectories do not tend to infinity with increasing time.
\par
Set $X_{G}=\displaystyle\frac{r_{G}}{\mu_{G}}$ and
$Y_{T}=\displaystyle\frac{r_{T}}{\mu_{T}}$, where
$r_{G}=\gamma_{G}-\delta_{G0}$ and $r_{T}=\gamma_{T}-\delta_{T}$
 are  net primary production of
grass and tree biomasses respectively. The following lemma holds.

\begin{lem}
\label{lem_compact}
 When $\gamma_{G}>\delta_{G0}$ and
$\gamma_{T}>\delta_{T}$, the compact
\begin{displaymath}
\mathcal{B}=\left\{(G,T)^{T}\in\mathbf{R}_{+}^{2}/G\leq X_{G}=K_{G}\left(1-\displaystyle\frac{\delta_{G0}}{\gamma_{G}}\right), T\leq
Y_{T}=K_{T}\left(1-\displaystyle\frac{\delta_{T}}{\gamma_{T}}\right) \right\}
\end{displaymath}
is positively invariant and attracting for the system
(\ref{Impuleq1}). Note that solutions of system (\ref{Impuleq1}) are
bounded below by 0, and above by the carrying capacities of grass and
tree biomasses.
 \end{lem}
$\mbox{\bf{Proof:}}$ See appendix A, page \pageref{AppendixA}.\par

When $t\neq n\tau$, the right-hand side of system (\ref{Impuleq1})
is locally Lipschitz continuous on $\mathcal{B}$.
Thus, system (\ref{Impuleq1}) has a unique solution.

\subsection{Equilibria}
System (\ref{Impuleq1}) has constant and periodic equilibria.\\
First of all, it is obvious that $E_{00}=(0; 0)$ and
$E_{01}=(0,Y_{T})$ are "trivial" equilibria of system
(\ref{Impuleq1}). While $E_{00}$ represents the bare soil, $E_{01}$
is the constant forest equilibrium, like in Tchuinte et al. 2014 \cite{Tchuinte2014}.\par
 Now, let us show the
existence of periodic solutions of the impulsive system
$(\ref{Impuleq1})$. The existence of the periodic grassland
equilibrium depends on the following threshold
\begin{displaymath}
  \mathcal{R}_{0,pulse}^{\tilde{G}_{e}}=\displaystyle\frac{r_{G}}{\displaystyle\frac{1}{\tau}\ln\left(\displaystyle\frac{1}{1-\lambda_{fG}}\right)}.
\end{displaymath}

 The following theorem holds.
 \begin{thm}
 \label{thm_grassland}(Semi-trivial periodic equilibrium)\par
 when $\mathcal{R}_{0,pulse}^{\tilde{G}_{e}}>1$, System $(\ref{Impuleq1})$ has a periodic grassland equilibrium $\tilde{E}_{10}=(\tilde{G}_{e}(t); 0)$, 
  where
\begin{displaymath}
  \tilde{G}_{e}(t)=X_{G}\displaystyle\frac{\{(1-\lambda_{fG})e^{r_{G}\tau}-1\}e^{r_{G}(t-n\tau)}}{\{(1-\lambda_{fG})e^{r_{G}\tau}-1\}e^{r_{G}(t-n\tau)}+\lambda_{fG}e^{r_{G}\tau}},\hspace{0.5cm}t\in
[n\tau, (n+1)\tau[, n=0,1,2,...
  \end{displaymath}
\end{thm}
$\mbox{\bf{Proof:}}$ See appendix B, page \pageref{AppendixB}.\par

Let us now show that there exists an unique positive periodic Tree-Grass equilibrium $(\tilde{G}^{*}(t); \tilde{T}^{*}(t))$.
  We set
  
  \begin{equation}
  	\label{Impuleq55a}
  	\begin{array}{lcl}
  		m(t,n\tau,T(n\tau))&=&r_{G}(t-n\tau)+\displaystyle\frac{\gamma_{TG}}{\mu_{T}}\ln\left[\displaystyle\frac{1}{1+\displaystyle\frac{T(n\tau)}{Y_{T}}\left[e^{r_{T}(t-n\tau)}-1\right]}\right],
  	\end{array}
  \end{equation}
  
  \begin{equation}
   \chi(t,n\tau,T(n\tau))=e^{m(t,n\tau,T(n\tau))},
  \label{Impuleq44a}
  \end{equation}
  \begin{equation}
   \label{Impuleq44b}
   G^{*}=\displaystyle\frac{[(1-\lambda_{fG})\chi((n+1)\tau,n\tau,T^{*})-1]}{\mu_{G}\int_{n\tau}^{(n+1)\tau}\chi(u,n\tau,T^{*})du},\hspace{0.5cm}n=0,1,2,...,
   \end{equation}
   and,
   \begin{equation}
   \label{Impuleq44c}
  T^{*}=
  \displaystyle\frac{Y_{T}\{(1-\lambda_{fT}\omega(\lambda_{fG}G^{*}))e^{r_{T}\tau}-1\}}{(e^{r_{T}\tau}-1)}
   \end{equation}
    \begin{equation}
\mathcal{R}_{0,pulse}^{*}=\displaystyle\frac{r_{T}}{\displaystyle\frac{1}{\tau}\ln\left(\displaystyle\frac{1}{1-\lambda_{fT}\omega(\lambda_{fG}\tilde{G}_{e}(\tau))}\right)},
  \end{equation}
 where
 \begin{displaymath}
 \tilde{G}_{e}(\tau)=X_{G}\displaystyle\frac{(1-\lambda_{fG})e^{r_{G}\tau}-1}{e^{r_{G}\tau}-1}.
 \end{displaymath}

  \begin{thm}(Uniqueness of the non-trivial periodic equilibrium)\par
\label{thm11} When  $\mathcal{R}_{0,pulse}^{\tilde{G}_{e}}>1$ and
$\mathcal{R}_{0,pulse}^{*}>1$, then
 system $(\ref{Impuleq1})$ has a unique positive Tree-Grass periodic equilibrium $\tilde{E}_{11}^{*}=(\tilde{G}^{*}(t); \tilde{T}^{*}(t))$, where
\begin{displaymath}
\tilde{G}^{*}(t)=\displaystyle\frac{\chi(t,n\tau,T^{*})G^{*}}{1+\mu_{G}G^{*}\int_{n\tau}^{t}\chi(u,n\tau,T^{*})du}\hspace{0.25cm}\mbox{and}\hspace{0.25cm}\tilde{T}^{*}(t)=
\displaystyle\frac{e^{r_{T}(t-n\tau)}T^{*}}{1+\displaystyle\frac{T^{*}}{Y_{T}}(e^{r_{T}(t-n\tau)}-1)},
\end{displaymath}
 $n\tau\leq t<(n+1)\tau,\hspace{0.5cm}n=0,1,2,...,$ such that 
 $\tilde{G}^{*}(\tau)=G^{*}$ and  $\tilde{T}^{*}(\tau)=T^{*}$.
\end{thm}
$\mbox{\bf{Proof:}}$ See  appendix C, page \pageref{AppendixC}.\par

\subsection{Local stability of the "trivial" equilibria}

 Set
 \begin{displaymath}
 \mathcal{R}_{01}=\displaystyle\frac{r_{G}}{r_{T}}\times\displaystyle\frac{\mu_{T}}{\gamma_{TG}},\hspace{0.5cm}\mbox{and}\hspace{0.5cm}
 \tilde{\mathcal{R}}_{0,\mathcal{R}_{01}}=\mathcal{R}_{0,pulse}^{\tilde{G}_{e}}\left(1-\displaystyle\frac{1}{\mathcal{R}_{01}}\right).
 \end{displaymath}

 \begin{thm}
\label{thm12}(Local stability of constant equilibria)
\begin{enumerate}
\item $E_{00}=(0; 0)$ is always unstable.
\item If $\mathcal{R}_{01}\leq 1$, then the forest equilibrium $E_{01}$ is locally asymptotically stable (LAS) (similarly as in the continuous model (see
Tchuinte et al. (2014) \cite{Tchuinte2014}).
\item If $\mathcal{R}_{01}>1$ and  $\tilde{\mathcal{R}}_{0,\mathcal{R}_{01}}<1$, then $E_{01}$ is LAS.
This situation is specific for the impulse model. The continuous
model does not imply the stability of the forest when
$\mathcal{R}_{01}>1$.
\item If  $\mathcal{R}_{01}>1$ and $\tilde{\mathcal{R}}_{0,\mathcal{R}_{01}}>1$, then $E_{01}$ is unstable.
\end{enumerate}
 \end{thm}
$\mbox{\bf{Proof:}}$ See appendix D, page \pageref{AppendixD}.\par

\subsection{Local stability of periodic equilibria}

We begin to investigate the local asymptotic stability of the
periodic grassland equilibrium of system $(\ref{Impuleq1})$. To
complete this subsection, we show the local stability of the
periodic savanna equilibrium.\par
 Recall that

\begin{displaymath}
\mathcal{R}_{0,pulse}^{\tilde{G}_{e}}=\displaystyle\frac{r_{G}}{\displaystyle\frac{1}{\tau}\ln\left(\displaystyle\frac{1}{1-\lambda_{fG}}\right)},\hspace{0.2cm}\mbox{and}\hspace{0.2cm}
\mathcal{R}^{*}_{0,pulse}=\displaystyle\frac{r_{T}}{\displaystyle\frac{1}{\tau}\ln\left(\displaystyle\frac{1}{1-\lambda_{fT}\omega(\lambda_{fG}\tilde{G}_{e}(\tau))}\right)}
\end{displaymath}
 where
 \begin{displaymath}
 \tilde{G}_{e}(t)=X_{G}\displaystyle\frac{\{(1-\lambda_{fG})e^{r_{G}\tau}-1\}e^{r_{G}(t-n\tau)}}{\{(1-\lambda_{fG})e^{r_{G}\tau}-1\}e^{r_{G}(t-n\tau)}+\lambda_{fG}e^{r_{G}\tau}},\hspace{0.5cm}t\in [n\tau, (n+1)\tau[, n=0,1,2,....
 \end{displaymath}

\begin{thm}
\label{thm13}  If $\mathcal{R}^{\tilde{G}_{e}}_{0,pulse}>1$ and
$\mathcal{R}^{*}_{0,pulse}<1$, then the periodic grassland
equilibrium $\tilde{E}_{10}=(\tilde{G}_{e}(t), 0)$ is locally
asymptotically stable.
\end{thm}
$\mbox{\bf{Proof:}}$ See appendix E, page \pageref{AppendixE}.\par

Now, we investigate local properties of the periodic savanna
equilibrium. Set
 \begin{displaymath}
 \mathcal{R}_{0,stable}^{*}=\displaystyle\frac{r_{T}}{\displaystyle\frac{1}{\tau}\ln\left(\displaystyle\frac{1}{1-\lambda_{fT}\omega(\lambda_{fG}G^{*})}\right)},\hspace{0.5cm}\mathcal{\tilde{R}}_{0,\mathcal{R}_{01}}^{G^{*}}=\displaystyle\frac{1}{\mathcal{R}_{01}}\left(1-\displaystyle\frac{1}{\mathcal{R}_{0,stable}^{*}}\right),
 \end{displaymath}
 and
\begin{displaymath}
\mathcal{\tilde{R}}_{0,stable}^{**}=\mathcal{\tilde{R}}_{0,\mathcal{R}_{01}}^{G^{*}}+\displaystyle\frac{1}{\mathcal{R}_{0,pulse}^{G_{e}}}+\displaystyle\frac{2}{X_{G}}\left(\displaystyle\frac{1}{\tau}\int_{0}^{\tau}\tilde{G}^{*}(u)du\right)
\end{displaymath}
where  $G^{*}$ and $T^{*}$ are defined in  (\ref{Impuleq44b}) and
(\ref{Impuleq44c}) respectively.  The threshold $
\mathcal{R}_{0,stable}^{*}$ represents the net production of tree
biomass relative to fire-induced tree biomass loss at the mixed
Tree-Grass equilibrium.\par

The following theorem holds.
 \begin{thm}
 \label{thm14} When  $\mathcal{R}_{0,stable}^{*}>1$ and $\mathcal{\tilde{R}}_{0,stable}^{**}>1$,
 the  savanna periodic equilibrium  $\tilde{E}_{11}^{*}=(\tilde{G}^{*}(t);\tilde{T}^{*}(t))$ is LAS.
 with
 \begin{displaymath}
 \tilde{G}^{*}(t)=\displaystyle\frac{\chi(t,n\tau,T^{*})G^{*}}{1+\mu_{G}G^{*}\int_{n\tau}^{t}\chi(u,n\tau,T^{*})du}\hspace{0.25cm}\mbox{and}\hspace{0.25cm}\tilde{T}^{*}(t)= \displaystyle\frac{e^{r_{T}(t-n\tau)}T^{*}}{1+\displaystyle\frac{T^{*}}{Y_{T}}(e^{r_{T}(t-n\tau)}-1)},
 \end{displaymath}
  $n\tau\leq t<(n+1)\tau,\hspace{0.5cm}n=0,1,2,...$, where  the expression of  $\chi$ is given in (\ref{Impuleq44a}).
 \end{thm}
$\mbox{\bf{Proof:}}$ See appendix F, page \pageref{AppendixF}.\par

The local stability is sufficient when there are multiple stable
states. However, for an unique equilibrium, the global stability is
necessary to ensure that all trajectories converge to the
equilibrium.

\begin{rmq}
We compute $G^{*}$ and $T^{*}$  using (\ref{Impuleq44b}) and
(\ref{Impuleq44c}) and a specific command ("fzero") in matlab which determines the fixed point. This allows us to obtain $\mathcal{R}_{0,stable}^{*}$. We use  numerical approximations of $\tilde{G}^{*}(t)$  to have $\mathcal{\tilde{R}}_{0,stable}^{**}$.
\end{rmq}

 \subsection{Global stability of equilibria}

In this section, we investigate the global stability of the forest
equilibrium and the periodic grassland equilibrium. The  following
theorem holds.
 \begin{thm}(Forest equilibrium GAS)\\
 \label{thm15}
The Forest  equilibrium $E_{01}=(0; Y_{T})$ is globally
asymptotically stable when
$\mathcal{R}_{0,pulse}^{\tilde{G}_{e}}<1$.
 \end{thm}
$\mbox{\bf{Proof:}}$ See appendix G, page \pageref{AppendixG}.\\

The global stability of the periodic grassland equilibrium is given
in the following theorem.
 \begin{thm}(Grassland  periodic equilibrium GAS)\\
 \label{thm16}
If $\mathcal{R}_{0,pulse}^{\tilde{G}_{e}}>1$ and
$\mathcal{R}_{0,pulse}^{*}<1$, then the grassland periodic
equilibrium $\tilde{E}_{10}=(\tilde{G}_{e}(t);0)$ is globally
asymptotically stable, where
 \begin{displaymath}
  \tilde{G}_{e}(t)=X_{G}\displaystyle\frac{\{(1-\lambda_{fG})e^{r_{G}\tau}-1\}e^{r_{G}(t-n\tau)}}{\{(1-\lambda_{fG})e^{r_{G}\tau}-1\}e^{r_{G}(t-n\tau)}+\lambda_{fG}e^{r_{G}\tau}},\hspace{0.5cm}t\in [n\tau, (n+1)\tau[,
  n=0,1,2,...,
  \end{displaymath}
  \begin{displaymath}
\mathcal{R}_{0,pulse}^{\tilde{G}_{e}}=\displaystyle\frac{r_{G}}{\displaystyle\frac{1}{\tau}\ln\left(\displaystyle\frac{1}{1-\lambda_{fG}}\right)},\hspace{0.2cm}\mbox{and}\hspace{0.2cm}
\mathcal{R}^{*}_{0,pulse}=\displaystyle\frac{r_{T}}{\displaystyle\frac{1}{\tau}\ln\left(\displaystyle\frac{1}{1-\lambda_{fT}\omega(\lambda_{fG}\tilde{G}_{e}(\tau))}\right)}.
\end{displaymath}
 \end{thm}
$\mbox{\bf{Proof:}}$
See appendix H.\par

\begin{rmq}
Consider $G^{*}(\tau)\leq G_{e}(\tau)$. We show that $\mathcal{R}_{0,pulse}^{*}\leq \mathcal{R}_{0,stable}^{*}$. Therefore
\begin{itemize}
	\item $\mathcal{R}_{0,stable}^{*}<1\Rightarrow \mathcal{R}_{0,pulse}^{*}<1$,
	\item $\mathcal{R}_{0,pulse}^{*}>1\Rightarrow\mathcal{R}_{0,stable}^{*}>1$.
\end{itemize} 
\end{rmq}

 Thresholds and their ecological meaning are recalled in the following table \ref{tb2}

\begin{table}[H]
	{\small
		\begin{center}
			\caption{The thresholds and their ecological meaning}
			\begin{tabular}{|l||c|}
				\hline
				Thresholds & Ecological meaning\\
				\hline
				&  the net primary production of grasses relative to the grass\\
				$
				\mathcal{R}_{01}=\displaystyle\frac{r_{G}}{r_{T}}\times\displaystyle\frac{\mu_{T}}{\gamma_{TG}}$
				&  production loss due to tree/grass competition
				\\
				&  throughout their life at the close forest equilibrium\\
				\hline\hline
				&\\
				$
				\mathcal{R}_{0,pulse}^{\tilde{G}_{e}}=\displaystyle\frac{r_{G}}{\displaystyle\frac{1}{\tau}\ln\left(\displaystyle\frac{1}{1-\lambda_{fG}}\right)}$
				& It is the  net primary production of grasses after fire\\
				&\\
				\hline\hline
				& \\
				$\tilde{\mathcal{R}}_{0,\mathcal{R}_{01}}=\mathcal{R}_{0,pulse}^{\tilde{G}_{e}}\left(1-\displaystyle\frac{1}{\mathcal{R}_{01}}\right)$ & It is the mixed threshold of the two previous thresholds\\
				&\\
				\hline\hline
				& represents the net production of tree biomass relative\\
				$\mathcal{R}_{0,pulse}^{*}=\displaystyle\frac{r_{T}}{\displaystyle\frac{1}{\tau}\ln\left(\displaystyle\frac{1}{1-\lambda_{fT}\omega(\lambda_{fG}\tilde{G}_{e}(\tau))}\right)}$ & to fire-induced biomass loss at the period of fire\\
				& at the grassland equilibrium\\
				\hline\hline
				& It is the net production of tree biomass relative\\
				$\mathcal{R}_{0,stable}^{*}=\displaystyle\frac{r_{T}}{\displaystyle\frac{1}{\tau}\ln\left(\displaystyle\frac{1}{1-\lambda_{fT}\omega(\lambda_{fG}G^{*})}\right)}$ & to fire-induced biomass loss at the period of fire\\
				& at the mixed Tree-Grass equilibrium\\
				\hline\hline
			\end{tabular}
			\label{tb2}
		\end{center}}
	\end{table}

We summarize all local and global properties of the impulsive model
(\ref{Impuleq1}) in table \ref{tabfinal}.

\begin{table}[H]
	{\footnotesize
		\begin{center}
			\caption{Long term behaviour of model (\ref{Impuleq1})}
			\renewcommand{\arraystretch}{1.5}
			\begin{tabular}{|l|c|c|c|c|c|c|c|c|c|}
				\cline{1-10}
				\multicolumn{6}{|c|}{\bf Thresholds} & \multirow{2}{1.5cm}{\bf Equilibria}  & \multirow{2}{1cm}{\bf Stable} & \multirow{2}{1.5cm}{\bf Unstable} & \multirow{2}{0.6cm}{\bf Case}\\
				\cline{1-6}
				$\mathcal{R}_{01}$ &  $\mathcal{R}_{0,pulse}^{\tilde{G}_{e}}$ & $\mathcal{\tilde{R}}_{0,\mathcal{R}_{01}}$ & $\mathcal{R}_{0,pulse}^{*}$ & $\mathcal{R}_{0,stable}^{*}$& $\mathcal{\tilde{R}}_{0,stable}^{**}$ & & & &  \\
				\hline
				\cline{1-2} \multirow{9}{1cm}{$<1$} & <1 &  
				&-&-&-& $E_{00}$, $E_{01}$ & $E_{01}$ $\mathbf{(GAS)}$ & $E_{00}$ & \bf{I}\\
				\cline{2-2}	\cline{4-10}
				& \multirow{6}{1cm}{>1}  &  & <1 & - & - & $E_{00}$, $E_{01}$ & $E_{01}$ $\mathbf{(LAS)}$ & $E_{00}$ & \bf{II}\\
				& & - & & & & $\tilde{E}_{10}$& $\tilde{E}_{10}$ $\mathbf{(LAS)}$& & \\
				\cline{4-10}
				& & & \multirow{4}{1cm}{>1}  & \multirow{3}{1cm}{>1}  & >1 & $E_{00}$, $E_{01}$& $E_{01}$ $\mathbf{(LAS)}$ & $E_{00}$ & \bf{III}\\
				& & & & & & $\tilde{E}_{10}$, $\tilde{E}_{11}^{*}$ & $\tilde{E}_{11}^{*}$ $\mathbf{(LAS)}$& $\tilde{E}_{10}$  & \\
				\cline{6-10}
				& & & & &  <1 & $E_{00}$, $E_{01}$&  & $E_{00}$, $\tilde{E}_{10}$ & \bf{IV}\\
				\cline{5-6}
				& & & & <1 &- & $\tilde{E}_{10}$, $\tilde{E}_{11}^{*}$ & $E_{01}$ $\mathbf{(GAS)}$ & $\tilde{E}_{11}^{*}$  & \\
				\cline{5-10}
				& & & &  & & $E_{00}$, $E_{01}$ &  & $E_{00}$  & \\
				& & & & - &- & $\tilde{E}_{10}$ & $E_{01}$ $\mathbf{(GAS)}$ & $E_{01}$  & \bf{V}\\
				\hline
				\cline{1-1}
				\multirow{18}{1cm}{$\geq 1$} & <1 & 
				&-&-&-& $E_{00}$, $E_{01}$ & $E_{01}$ $\mathbf{(GAS)}$ & $E_{00}$ & \bf{VI}\\
				\cline{2-2}\cline{4-10}
				& \multirow{6}{1cm}{>1}  &  & <1 & - & - & $E_{00}$, $E_{01}$ & $E_{01}$ $\mathbf{(LAS)}$ & $E_{00}$ & \bf{VII}\\
				& & $\leq 1$ & & & & $\tilde{E}_{10}$& $\tilde{E}_{10}$ $\mathbf{(LAS)}$& & \\
				\cline{4-10}
				& & & \multirow{4}{1cm}{>1}  & \multirow{3}{1cm}{>1}  & >1 & $E_{00}$, $E_{01}$& $E_{01}$ $\mathbf{(LAS)}$ & $E_{00}$ & \bf{VIII}\\
				& & & & & & $\tilde{E}_{10}$, $\tilde{E}_{11}^{*}$ & $\tilde{E}_{11}^{*}$ $\mathbf{(LAS)}$& $\tilde{E}_{10}$  & \\
				\cline{6-10}
				& & & & &  <1 & $E_{00}$, $E_{01}$&  & $E_{00}$, $\tilde{E}_{10}$ & \bf{IX}\\
				\cline{5-6}
				& & & & <1 &- & $\tilde{E}_{10}$, $\tilde{E}_{11}^{*}$ & $E_{01}$ $\mathbf{(GAS)}$ & $\tilde{E}_{11}^{*}$  & \\
				\cline{5-10}
				& & & & - &- & $E_{00}$, $\tilde{E}_{10}$ &  & $E_{00}$  & \bf{X}\\
				& & & &  & & $E_{01}$ & $E_{01}$ $\mathbf{(GAS)}$ &  $\tilde{E}_{10}$ & \\
				\cline{2-10}
				& <1 & 
				&-&-&-& $E_{00}$, $E_{01}$ & $E_{01}$ $\mathbf{(GAS)}$ & $E_{00}$ & \bf{XI}\\
				\cline{2-2}	\cline{4-10}
				& & \multirow{2}{1cm}{}  & <1 & - & - & $E_{00}$, $E_{01}$ &  & $E_{00}$, $E_{01}$ & \bf{XII}\\
				& & >1 & & & & $\tilde{E}_{10}$ & $\tilde{E}_{10}$ $\mathbf{(GAS)}$&  & \\
				\cline{4-10}
				& & & >1 &  \multirow{2}{1cm}{} & >1 & $E_{00}$, $E_{01}$ & & $E_{00}$, $E_{01}$ & \bf{XIII}\\
				& >1 &  & & >1 & & $\tilde{E}_{10}$, $\tilde{E}_{11}^{*}$ & $\tilde{E}_{11}^{*}$ $\mathbf{(GAS)}$ & $\tilde{E}_{10}$ & \\
				\cline{6-10}
				& & & & & <1 & $E_{00}$, $E_{01}$ & &$E_{00}$, $\tilde{E}_{10}$ & \bf{XIV}\\
				\cline{5-6}
				& & & & <1 &- & $\tilde{E}_{10}$, $\tilde{E}_{11}^{*}$ & $E_{01}$ $\mathbf{(GAS)}$ & $\tilde{E}_{11}^{*}$  & \\
				\cline{5-10}
				& & & & - &- & $E_{00}$, $\tilde{E}_{10}$ &  & $E_{00}$  &\\
				& & & &  & & $E_{01}$ & $\tilde{E}_{10}$ $\mathbf{(GAS)}$ & $E_{01}$  & \bf{XV}\\
				\hline
			\end{tabular}
			\label{tabfinal}
		\end{center}}
	\end{table}

Table \ref{tabfinal} gives all possible configurations for the impulsive system. According to the values taken by the thresholds, it may be possible to anticipate the long term behaviour of the system. Since many configurations are possible, it is essential to highlights the parameters that may have an important impact on the thresholds. In the next sections, we briefly present the numerical algorithm we have chosen to perform numerical simulations and present some results emphasizing the importance of the competition parameter $\gamma_{TG}$.

\section{The numerical  algorithm}
In the previous section, solutions were searched in form of
analytical expression. However, many impulsive differential
equations can not be solved in this way or their solving is more
complicated in the mathematical point of view.\par

A nonstandard numerical scheme for solving the
impulsive differential equation is built.  The nonstandard approach relies on the following
important rules: the standard denominator $\Delta t$ in each
discrete derivative is replaced by a time-step function
$0<\varphi(\Delta t)<1$ ; such that $\varphi(\Delta t)=\Delta
t+\mathcal{O}(\Delta t)$; the nonlinear terms are approximated in a
non local way; for instance the nonlinear term $T(t_{n})G(t_{n})$ in
the  problem can be approximated by $T^{n}G^{n+1}$. For an overview
and some applications in Biology of the nonstandard finite difference method see for instance (Anguelov et al., 2012 \cite{Anguelov2012}; Anguelov et al., 2014 \cite{Anguelov2014}).

The nonstandard approximations for system (\ref{Impuleq1}) are given by
  \begin{equation}
    \label{Impuleq112}
    \left\{
    \begin{array}{lcl}
     \displaystyle\frac{G^{n+1}-G^{n}}{\varphi_{1}(\Delta t)}&=& (\gamma_{G}-\delta_{G0})G^{n}-\mu_{G}G^{n}G^{n+1}-\gamma_{TG}T^{n}G^{n+1},\\
     \\
    \displaystyle\frac{T^{n+1}-T^{n}}{\varphi_{2}(\Delta t)}&=& (\gamma_{T}-\delta_{T})T^{n}-\mu_{T}T^{n}T^{n+1},\\
        \end{array}
    \right.
  \end{equation}
and
  \begin{equation}
    \label{Impuleq113}
    \left\{
    \begin{array}{lcl}
    G^{n+}&=& (1-\lambda_{fG})G^{n+1},\\
    \\
    T^{n+}&=& (1-\lambda_{fT}\omega(\lambda_{fG}G^{n}))T^{n+1},\\
        \end{array}
    \right.
  \end{equation}
where
\begin{equation}
 \label{Impuleq114}
 \varphi_{1}(\Delta t)=\displaystyle\frac{e^{(\gamma_{G}-\delta_{G0})\Delta t}-1}{\gamma_{G}-\delta_{G0}}
\end{equation}
 and
 \begin{equation}
  \label{Impuleq115}
\varphi_{2}(\Delta t)=\displaystyle\frac{e^{(\gamma_{T}-\delta_{T})\Delta t}-1}{\gamma_{T}-\delta_{T}}.
 \end{equation}
Scheme (\ref{Impuleq112}-b) with the time-step function (\ref{Impuleq115}) is an exact scheme, between each fire event. Similarly, when $\gamma_{TG}=0$, the scheme (\ref{Impuleq112}-a) with the time-step function (\ref{Impuleq114}) is also an exact scheme, between each fire event.
Altogether the numerical algorithm (\ref{Impuleq112})-(\ref{Impuleq113}) is positively stable and elementary stable i.e. it preserves equilibria and local properties of each equilibrium  of system  (\ref{Impuleq1}).
Thus at least locally, we are sure that schemes
(\ref{Impuleq112})-(\ref{Impuleq113})  replicate the dynamics of
system (\ref{Impuleq1}).

\section{Numerical simulations and discussion}
\subsection{Environmental setting}

To illustrate our analytical results and highlight important ecological parameters, we will perform some numerical simulations. In fact, from the ecological point of view, parameters can change drastically according to the environmental features. For instance, in Cameroon, three different zones can be particularly highlighted. The first zone is Region
1 (R1) where the biomass production is low (the Mean Annual Precipitations
(MAP) is less than $650$ mm  by year). The second zone is  Region 2
(R2) ($650$-$1100$ mm/yr). It has more biomass production by year.
The last zone 
 (R3) ($1100-1800$ mm by year) where savannas are observed in some cases in the immediate vicinity of forests (Favier et al. 2012 \cite{Favier2012abrupt}). Model (\ref{Impuleq1}) has 11 parameters.
Eight of them $(\gamma_{G}, K_{G}, \delta_{G0}, \alpha,
\gamma_{T}, K_{T}, \delta_{T}, \gamma_{TG})$ specify vegetation
growth while the others are related to the fire characteristics.
Parameter values used are based on literature sources.

\begin{figure}[H]
\begin{center}
\includegraphics[scale=1]{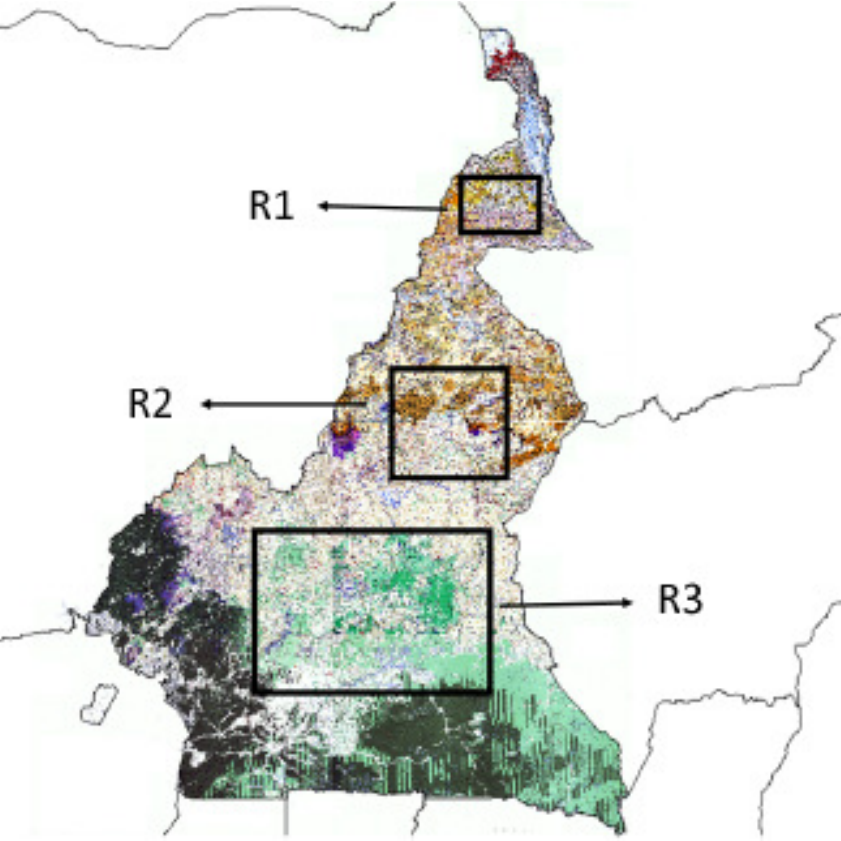}
\caption{\footnotesize {Overall vegetation map of Cameroon from Letouzey (1985). The three regions R1,
R2 and R3 are indicated by black rectangles.}} \label{CAMER}
\end{center}
\end{figure}

\begin{table}[H]
{\footnotesize
\begin{center}
 \caption{Vegetational Data}
\renewcommand{\arraystretch}{1.5}
\begin{tabular}{cccccc}
\hline \cline{1-6}
& \multicolumn{3}{c}{\bf Site name}& \\
\cline{2-4}
$\mathbf{Parameters}$  & $\mathbf{Semi-arid} $ & $\mathbf{Mesic}$ & $\mathbf{Humid}$ & $\mathbf{Units}$ & $\mathbf{Sources}$ \\
 & $\mathbf{R1 }$ & $\mathbf{R2}$ & $\mathbf{R3}$ &  & \\
   \hline\hline
   $K_{G}$ & $2-5$ & $8-10$ & $10-20$ & t.ha$^{-1}$ & \cite{penning1982productivity,Abbadie2006lamto}\\
 $\gamma_{G}$  & $0.4^{(*)}-1.5^{(**)}$ & $1.5-3$ & $3-4.6^{(***)}$ & yr$^{-1}$& $^{(*)}$, $^{(**)}$ \cite{Tucker1985}, $^{(***)}$ \cite{Menaut1979structure} \\
 $\delta_{G0}$ & $0$ & $0$ & $0-0.9$$^{(4*)}$ & yr$^{-1}$ & $^{(4*)}$  \cite{vanLangevelde2003}\\
 $\lambda_{fG}$ & $0.1-0.9$ &   $0.1-0.9$  &   $0.1-0.9$  & -& assumed, see also \cite{Abbadie2006lamto}\\
$\alpha$ & $2$ &$2$  &  $2$ & t.ha$^{-1}$& Assumed\\
$\theta$ & $2$ &$2$  &  $2$ & t.ha$^{-1}$& Assumed\\
\hline\hline
$K_{T}$& $10-25$ & $25-60$ & $60-115$ & t.ha$^{-1}$& \cite{mermoz2014biomass}\\
$\gamma_{T}$& $0.3^{(a)}-0.9$ & $0.9-1.2$  & $1.2-7.2^{(b)}$ &  yr$^{-1}$& $^{(a)}$ \cite{Tucker1985}, $^{(b)}$ \cite{Breman1995woody}\\
$\delta_{T}$ & $0$ & $0$ & $0-0.015^{(5*)}$ & yr$^{-1}$ & $^{(5*)}$ \cite{Hochberg1994influences}\\
$\lambda_{fT}$ & $0.1-0.5$ & $0.1-0.5$  & $0.1-0.4$  & - & Assumed\\
\hline\hline
$\tau$ & $\geq 10$ & $2-8$  & $0.5-2$ & yr & Overall expert-based knowledge, \\
& & & & &   in addition to \cite{Menaut1979structure,Favier2012abrupt}\\
$\gamma_{TG}$ & $(-0.01)-0.03$ & $0.01-0.08$ & $0.03-0.09$ & ha.t$^{-1}$.yr$^{-1}$&  \cite{Mordelet1993influence}, \cite{Abbadie2006lamto} \\
& & & & & (after reinterpretation)\\
\hline\hline
\end{tabular}
\label{3_regions}
\end{center}
} \footnotesize {\begin{note}\label{note_table} Range values for
parameters ($K_{G}$, $k_{T}$, $\lambda_{fG}$, $\lambda_{fT}$, $\tau$, $\gamma_{TG}$) used for the simulation runs. Early fire destroy only
$25\%$ (Abbadie et al., 2006 \cite{Abbadie2006lamto}) of grass
biomass, while late fires can destroy up to $90\%$ of
biomass.\end{note}}
\end{table}

Our numerical analysis focuses on Cameroon as part of Central Africa where we find a summary of African natural conditions, from humid equatorial climate near the Atlantic Ocean, up to the arid Sahelian tropical climate in the region of Lake Chad. \par

 The first site is R1.  It corresponds to
semi-arid zones.  Grasses may be
dominant while trees are generally of low stature.  Trees
are resource-limited, and the resource competition with grasses and between trees are
the key factor determining savanna existence (Baudena et al., 2014
\cite{Baudena2014forests}; Tchuinte et al., 2014
\cite{Tchuinte2014}). In Cameroon, R1 (small black rectangle in Figure
\ref{CAMER}) corresponds to the dry tropical climate of the extreme
North Country,  from Kaélé in Maroua and Mora, and from Yagoua to
Kousséri, Makary and Lake Chad. At the edges of Lake Chad, there is
only  3 months of rain with $500$ mm/yr.\par

The second region, R2, represents a mesic zone. Here the MAP varies
from $650$ mm/yr to $1100$ mm/yr. In Cameroon it is located at the
Northern part, from Adamawa to the Mandara Mountains (see the middle
black rectangle in fig \ref{CAMER}). There are two seasons covering
the whole of Adamawa plateau, from Banyo to Ngaroundere and
Meiganga. It has been argued, that mesic savannas are unstable compared to forest and disturbance-dependent with respect to fires (Sankaran et
al., 2005 \cite{Sankaran2005determinants}),
which prevent tree invasion, because they occur regularly
during the dry season. Grasses benefit from fire because
they recover faster than trees after fires, and profit of  open spaces to
growth. Thus grass-fire feedback is a characteristic feature that
leads to savanna or grassland persistence.\par

High rainfall occurring in wet African savannas directly reduce the role
of water as limiting factor. 
The last region (R3) that we are
interested in corresponds to humid areas where MAP varies between $1100$
mm/yr and $1800$ mm/yr. In this region, the high water resource available
enables high fuel (grass biomass) production and therefore fires are
more frequent and of greater impact on seedlings. The grass-fire feedback in R3 leads to a bistability of savanna and forest, as
shown using a simple continuous models (e.g. Tchuinte et al., 2014
\cite{Tchuinte2014}; Staver and Levin, 2012
\cite{Staver2012integrating}) and evidenced from remote sensing data by Favier et al. (2012) \cite{Favier2012abrupt}. In Cameroon, R3 (see the big black
rectangle in Figure \ref{CAMER}) encompasses two sub-zones: a sub-zone
of transition between
 equatorial and tropical climates, and a sub-zone which corresponds
 to the equatorial climate itself. Concerning the first one, the MAP is
 between $1100$ mm/yr and $1500$ mm/yr. It is observed from Bafia to
 Bertoua, Batouri and from Yoko to Betare Oya, Garoua Boulaï. The
 second sub-zone is a site  covering the entire South of the Country from
Yaounde (1564 mm/yr) to Yokadouma, from Ebolowa to Ambam, Mouloundou
and Ouesso (Congo). It extented near the Gabonese borber (1700mm).
In both the two sub-zones  there are four distinct seasons (two dry
seasons alternating with two wet seasons with unequal intensity). At the South Cameroon near the Gabonese
border (11 months of rainy season), there is a close canopy forest (see fig
\ref{CAMER}). Above R3 (at the highest end of the rainfall range which is $2000$ mm/yr), fires are totally suppressed and only
forests are observed, since grass growth is inhibited by tree shade.\par 

In the next section we present some numerical simulations. The fundamental tasks in studying disturbance are to discriminate between fluctuations that are extraordinary and those that are usual (McNaughton 1992 \cite{mcnaughton1992propagation}).
In our simulations, we assume that the ecological system is not impacted by Human and Animals (grazing and browsing), i.e. $\delta_{G0}=\delta_T=0$.

\subsection{Simulations to illustrate bifurcations due to $\gamma_{TG}$ in regions R1, R2, and R3}

\subsubsection{Simulations in region R1}

\begin{figure}[H]
\begin{center}
\includegraphics[scale=0.5]{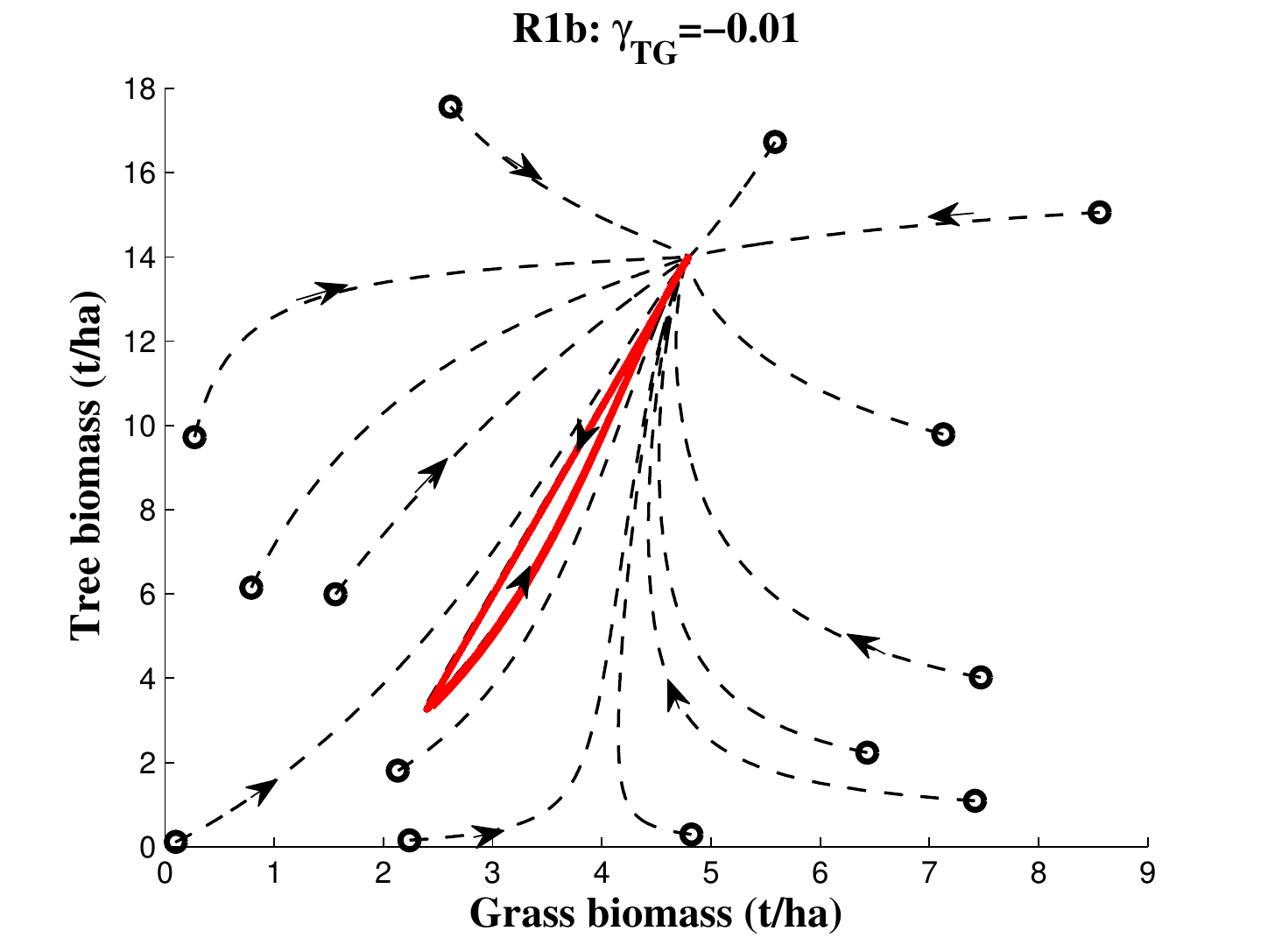}
\includegraphics[scale=0.5]{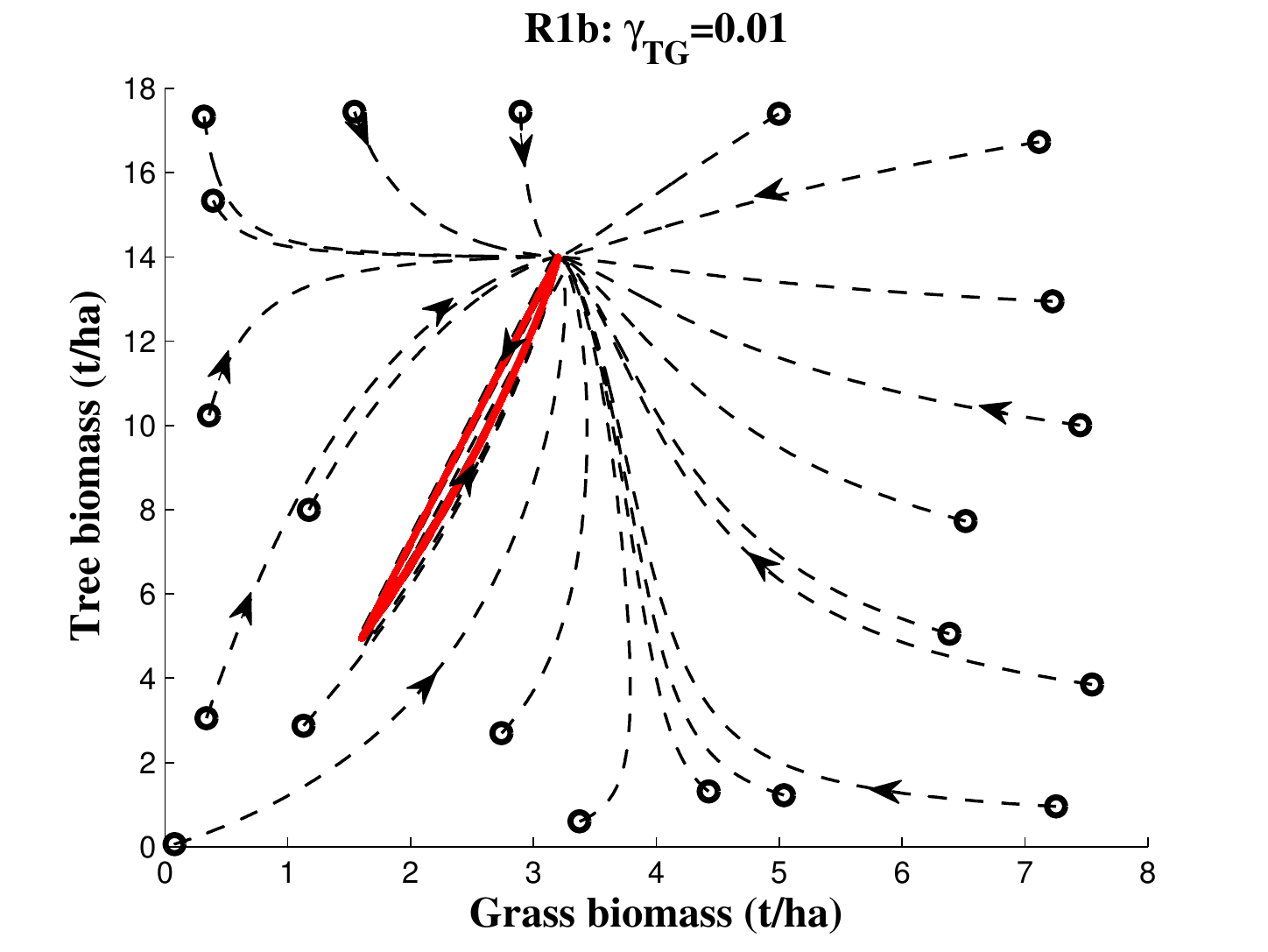}
\vspace{1cm}
\includegraphics[scale=0.5]{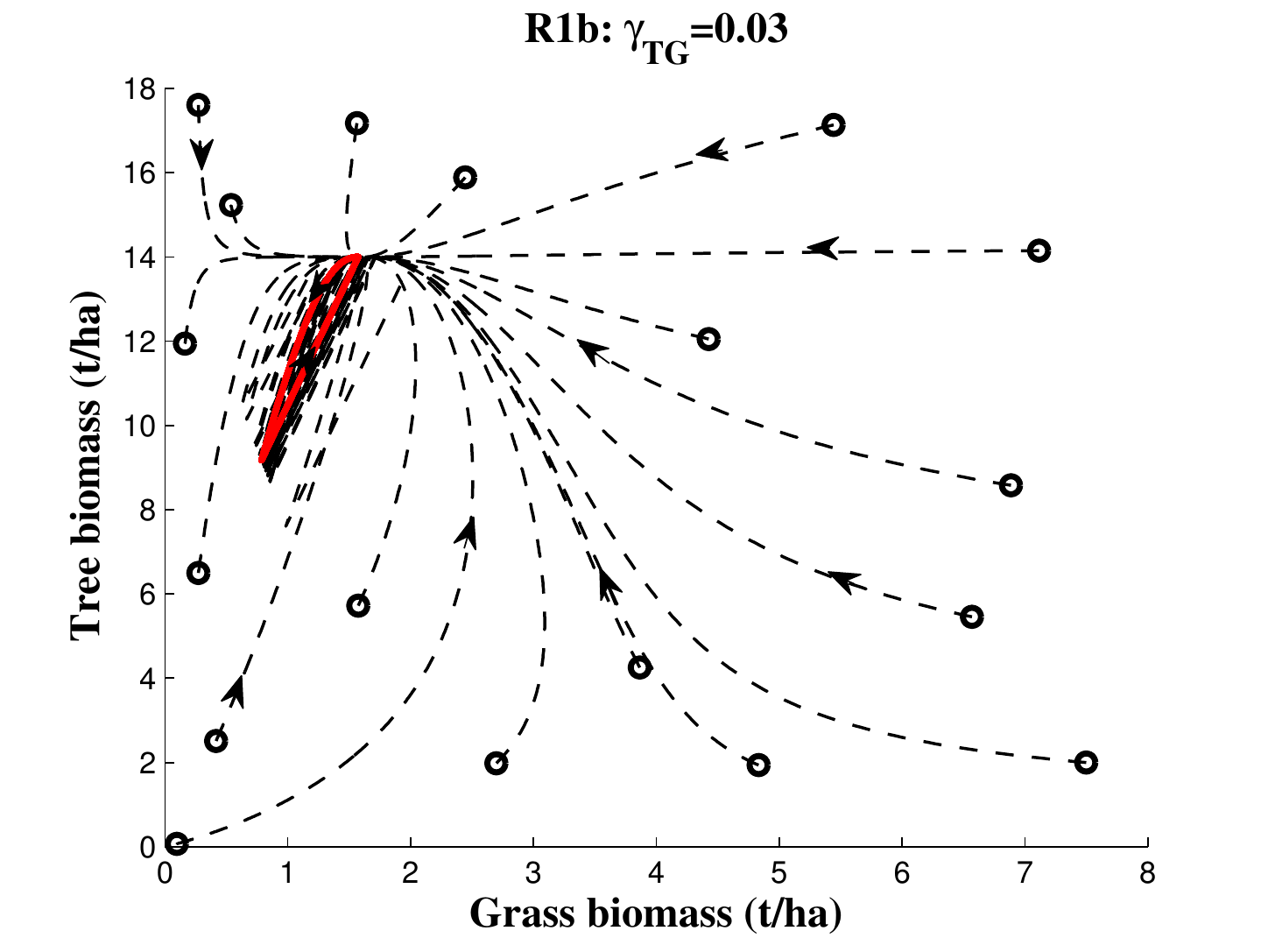}
\includegraphics[scale=0.5]{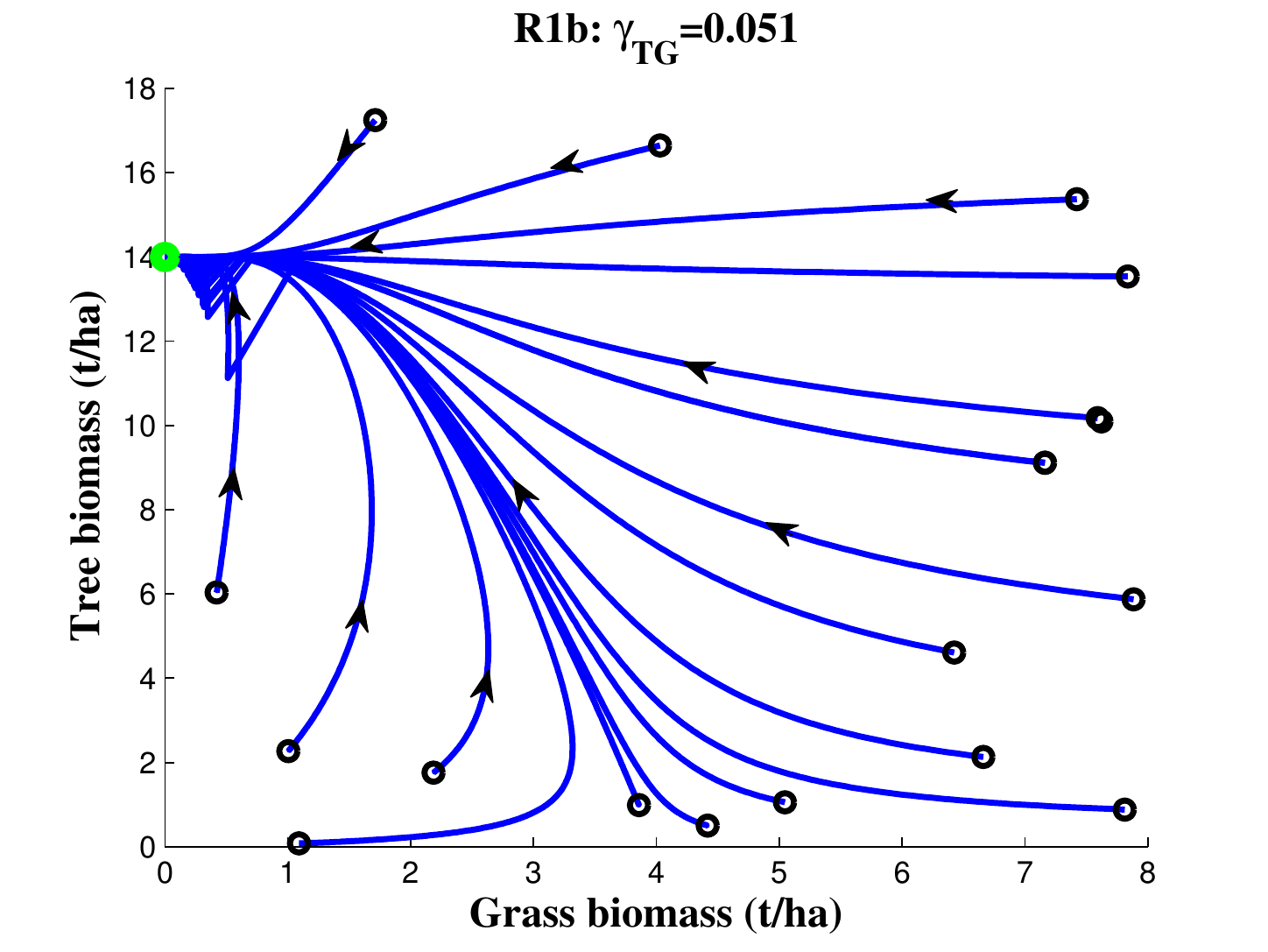}
\caption{\footnotesize { Phase diagrams  in R1.  Panels $\mathbf{a}$ and $\mathbf{b}$ ($\mathbf{c}$ and $\mathbf{d}$) show the impact of the tree/grass competition parameter (see panels panel $\mathbf{c}$ and
$\mathbf{d}$). When $\gamma_{TG}$ increases
	and exceeds a critical value, the attracting state of the system
	shifts from periodic savanna to a forest equilibrium (in fact, dense thickets).
}} \label{figR1_a}
\end{center}
\end{figure}

According to table  \ref{3_regions}, we choose the following values of parameters for region R1:
\begin{table}[H]
\caption{Parameters values related to figure \ref{figR1_a}}	
\begin{center}
\renewcommand{\arraystretch}{1.5}
\begin{tabular}{cccccccccc}
\hline
$K_{G}$ & $\gamma_{G}$ & $\delta_{G0}$ &  $K_{T}$&  $\gamma_{T}$ & $\delta_{T}$& $\alpha$ & $\tau$ & $\lambda_{fT}$ & $\lambda_{fG}$\\
\hline\hline
$4$ & $0.7$ & $0$ & $14$ & $0.75$ & $0$ & $2$ & $12$ & $0.9$ & $0.5$\\
\hline\hline
\end{tabular}
\label{params_R1}
\end{center}
\end{table}
Taking $\gamma_{TG}=[-0.01,0.01,0.03,0.051]$, and, using Table \ref{tabfinal}, we obtain Table \ref{R1}. Figure \ref{figR1_a} illustrates also the expected behaviours.

\begin{table}[H]
\caption{Thresholds Table related to Table \ref{params_R1} and Figure \ref{figR1_a}}	
\begin{center}
\renewcommand{\arraystretch}{1.5}
\begin{tabular}{cccccccccc}
	\hline
	\bf{Panel} & $\mathcal{R}_{01}$ & $\mathcal{R}_{0,pulse}^{\tilde{G}_{e}}$ & $\tilde{\mathcal{R}}_{0,\mathcal{R}_{01}}$  &  $\mathcal{R}_{0,pulse}^{*}$ & $\mathcal{R}_{0,stable}^{*}$ & $\mathcal{R}_{0,stable}^{**}$& \bf{Case}\\
	\hline\hline
	\bf{a} & $-$ & $ >1$ & $ >1$ & $ >1$ & $ >1$&$ >1$ & \bf{XIII} \\
	\hline
	\bf{b,c} & $ >1$ & $ >1$ & $ >1$ & $ >1$ & $ >1$ & $ >1$ & \bf{XIII} \\
	\hline
	\bf{d} & $<1$ &  $ >1$ & $ -$ & $ >1$ & -& - & \bf{V} \\
	\hline\hline
\end{tabular}
\label{R1}
\end{center}
\end{table}

In semi-arid areas,  Tree-Grass interactions are predominantly influenced by competition for soil water (e.g. Walker et al. 1981 \cite{Walker1981stability}). However, shading by tree foliage under arid climate can also increase grass production under the tree crown (Abbadie et al. 2006 \cite{Abbadie2006lamto}) and more generally can increase the water budget below the canopy (Barbier et al. 2008 \cite{Barbier2008}). Hence the influence of trees on grasses can range from facilitation to competition. Figure \ref{figR1_a} shows the influence of tree/grass interactions in R1. For higher values of $\gamma_{TG}$, trajectories converge to the forest equilibrium (see panel $\mathbf{d}$) corresponding, in fact to dense thickets. When $\gamma_{TG}$ is small, and even negative (positive effect on the grass biomass), the system converges to a periodic Tree-Grass coexistence equilibrium with fairly large amplitudes in the tree biomass (compare panels $\mathbf{a}$ and $\mathbf{b}$), or with small grass biomass (see panel $\mathbf{c}$). Here $\gamma_{TG}$ is an influential parameter since it permits a transition from savanna to forest. This result joins those of Sankaran et al. (2005) \cite{Sankaran2005determinants} which argued that in R1, savannas are stable in the sense that tree biomass and cover are primarily limited by resources (see
	panels $\mathbf{a}$ and $\mathbf{c}$ in Figure \ref{figR1_a}). Therefore, the competition parameter $\gamma_{TG}$ is an important
driver of Tree-Grass dynamics in R1 where fires return time are
typically higher than 10 years
and therefore are not necessary for grass-tree coexistence. Our
Tree-Grass impulsional model shows that in R1, there is only a
stable periodic grass-tree equilibrium or a stable forest
equilibrium (see panels $\mathbf{a}$ and $\mathbf{d}$ in Figure \ref{figR1_a}).\par

\subsubsection{Simulations in region R2}

\begin{figure}[H]
\begin{center}
\includegraphics[scale=0.5]{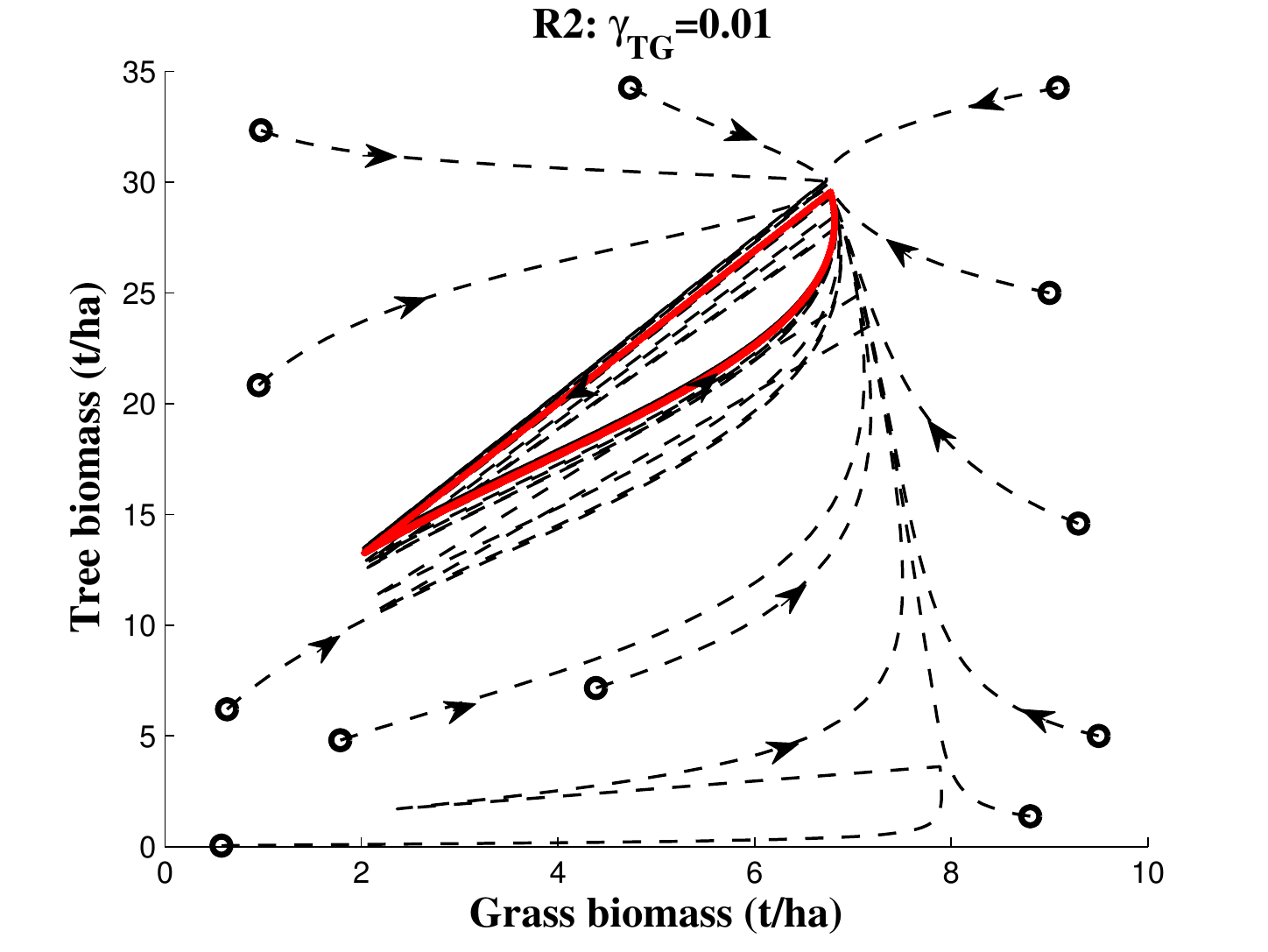}
\includegraphics[scale=0.5]{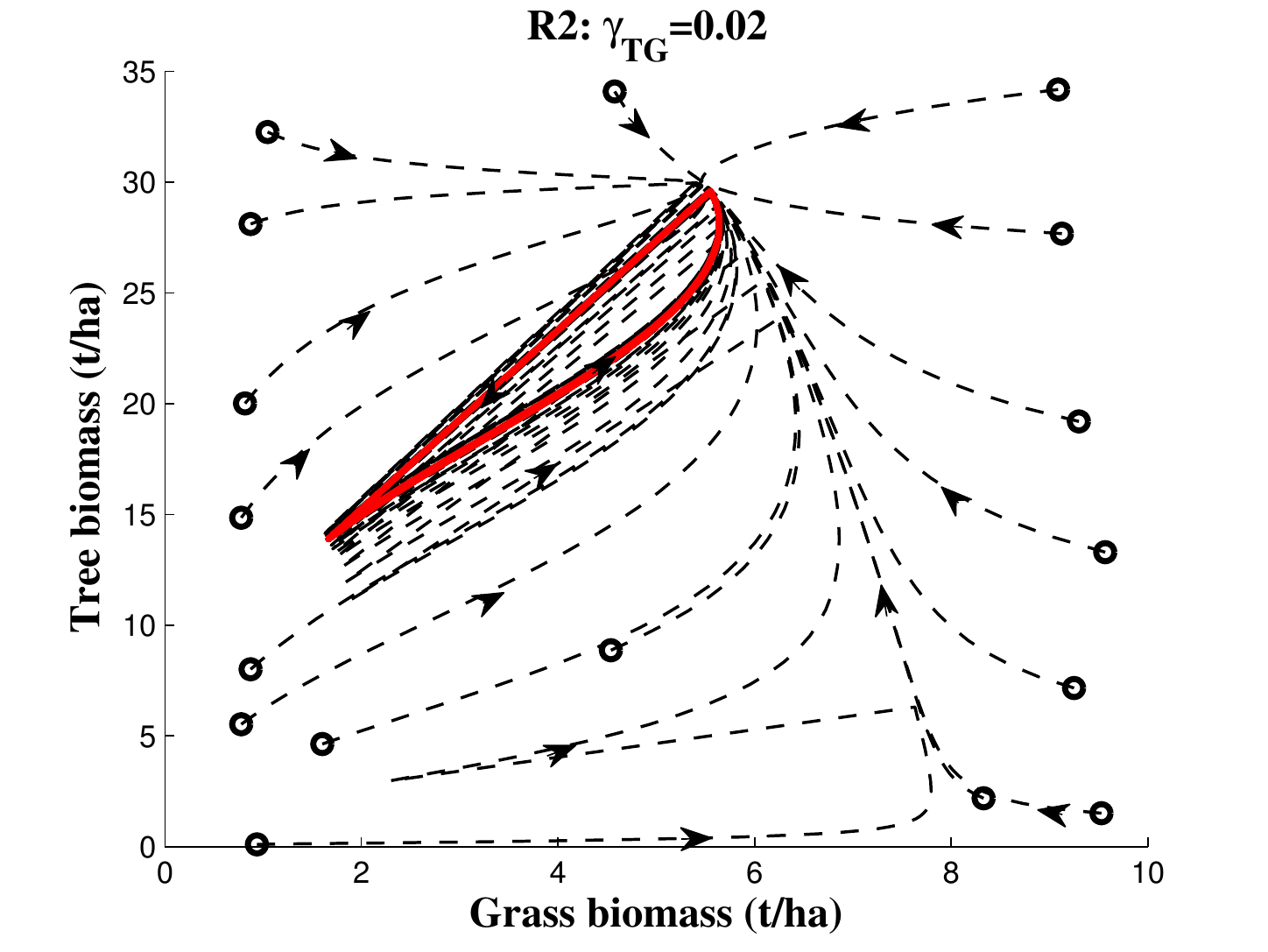}
\vspace{1cm}
\includegraphics[scale=0.5]{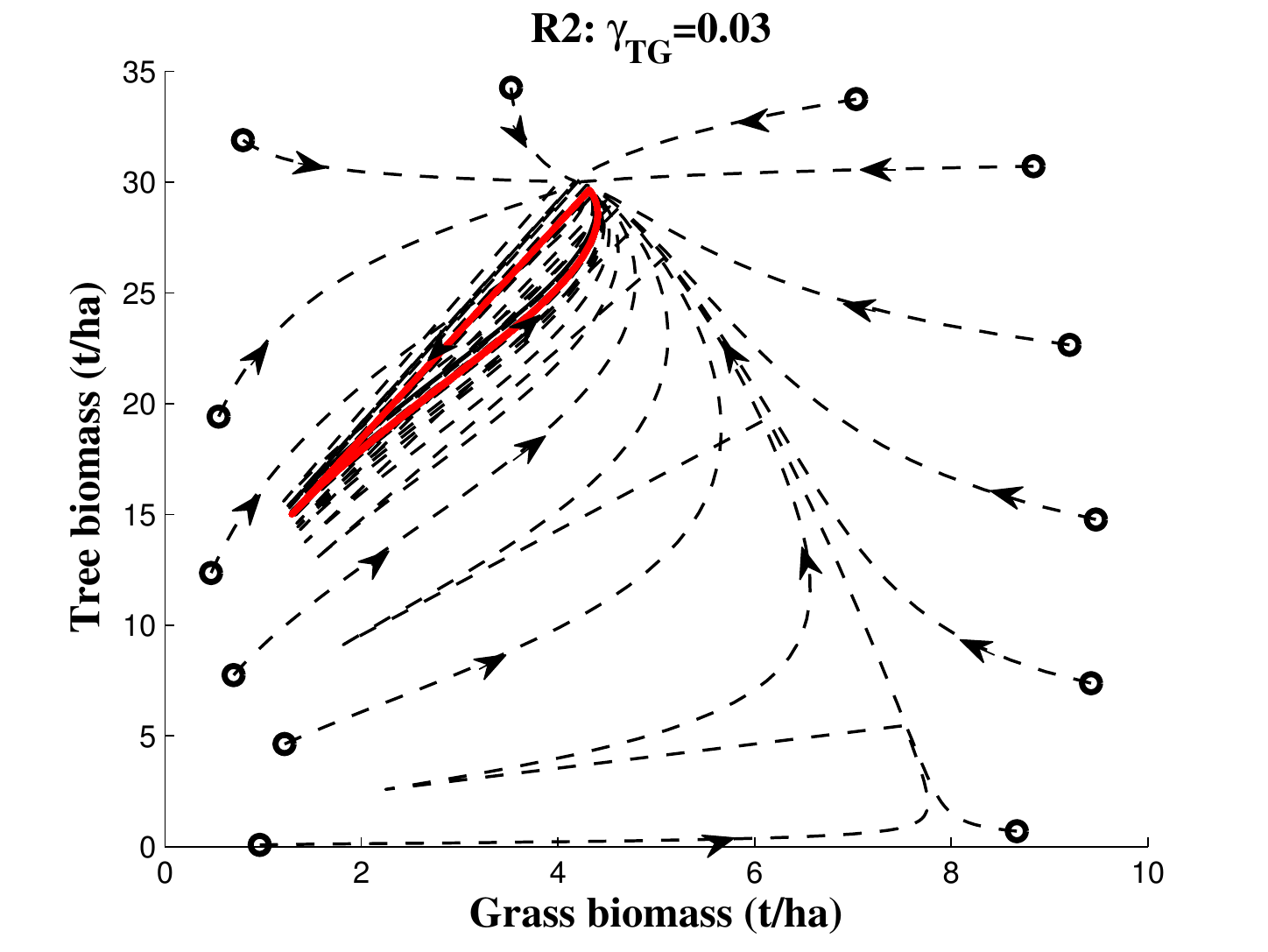}
\includegraphics[scale=0.5]{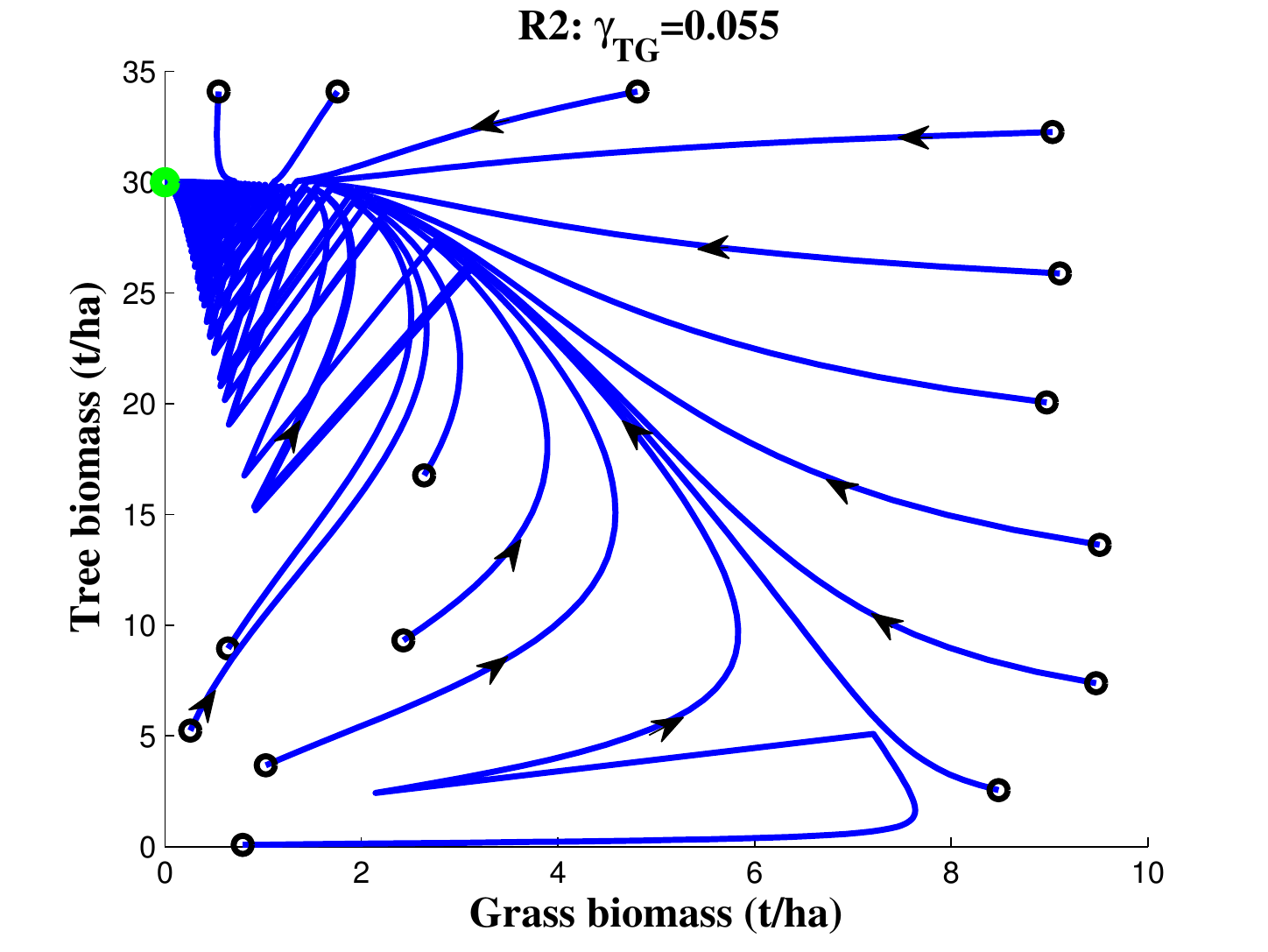}
\caption{Phase diagrams in R2 with a period of fire of $\tau=5$ years.} \label{figR2_a}
\end{center}
\end{figure}

Let us consider the following values of parameters according to table  \ref{3_regions}.
\begin{table}[H]
 	\caption{Parameters values related to figure \ref{figR2_a}}	
 \begin{center}
 \renewcommand{\arraystretch}{1.5}
 \begin{tabular}{cccccccccc}
 \hline
 $K_{G}$ & $\gamma_{G}$ & $\delta_{G0}$ &  $K_{T}$&  $\gamma_{T}$ & $\delta_{T}$& $\alpha$ & $\tau$ & $\lambda_{fT}$ & $\lambda_{fG}$\\
 \hline\hline
 $8$ & $1.9$ & $0$ & $30$ & $0.9$ & $0$ & $2$ & $5$ & $0.5$ & $0.6$\\
 \hline\hline
 \end{tabular}
 \end{center}
 \label{params_R2}
 \end{table}
According to Table \ref{tabfinal} and Table \ref{params_R2}, with $\gamma_{TG}=[0.01,0.02,0.03,0.055]$, we derive, in Table \ref{R2}, the behaviours of the Tree-Grass system. See also Figure \ref{figR2_a}, page \pageref{figR2_a}.

\begin{table}[H]
\caption{Thresholds Table related to Table \ref{params_R2} and Figure \ref{figR2_a}}	
\begin{center}
\renewcommand{\arraystretch}{1.5}
\begin{tabular}{ccccccccc}
	\hline
	\bf{Panel} & $\mathcal{R}_{01}$ & $\mathcal{R}_{0,pulse}^{\tilde{G}_{e}}$ & $\tilde{\mathcal{R}}_{0,\mathcal{R}_{01}}$ &  $\mathcal{R}_{0,pulse}^{*}$ & $\mathcal{R}_{0,stable}^{*}$ & $\mathcal{R}_{0,stable}^{**}$ &\bf{Case}\\
	\hline\hline
	\bf{a,b,c} & $>1$ & $ >1 $ & $>1$ & $>1$ & $>1$&$>1$ & \bf{XIII} \\
	\hline 
	\bf{d} & $>1$ & $>1$&  $ <1 $  & $>1$ & $-$ &$-$ & \bf{X} \\
	\hline\hline
\end{tabular}
\label{R2}
\end{center}
\end{table}
Figure \ref{figR2_a} below illustrates the bifurcation due to the competition parameter in Region (R2). The forest equilibrium is stable for $\gamma_{TG}=0.055$ a value in the upper range of plausible values (see panel $\mathbf{d}$). When $\gamma_{TG}$ decreases, the system converges to a savanna periodic equilibrium (see panel $\mathbf{a,b,c}$). We note also that $\gamma_{TG}$ has an impact on the amplitude of the periodic savanna equilibrium and the maximal amount of grass biomass.  \par

Increasing the impact of fire on trees via an increase of the "$\lambda_{fT}$" coefficient, a bistability between forest and savanna occurs in R2 (see figure \ref{figR2_b} below).
 \begin{table}[H]
 \caption{Parameters values related to figure \ref{figR2_b}}
 \begin{center}
 \renewcommand{\arraystretch}{1.5}
 \begin{tabular}{cccccccccc}
 \hline
 $K_{G}$ & $\gamma_{G}$ & $\delta_{G0}$ &  $K_{T}$&  $\gamma_{T}$ & $\delta_{T}$& $\alpha$ & $\tau$ & $\lambda_{fT}$ & $\lambda_{fG}$\\
 \hline\hline
 $8$ & $1.5$ & $0$ & $30$ & $0.9$ & $0$ & $2$ & $2.2$ & $0.8$ & $0.5$\\
 \hline\hline
 \end{tabular}
   \label{params_R2b}
 \end{center}	
 \end{table}
Using the same values for $\gamma_{TG}$, we derive Table \ref{R2b} and Figure  \ref{figR2_b}.
 \begin{table}[H]
 \caption{Thresholds Table related to Table \ref{params_R2b} and Figure \ref{figR2_b}}	
 \begin{center}
 \renewcommand{\arraystretch}{1.5}
 \begin{tabular}{cccccccccc}
 \hline
 \bf{Panel} & $\mathcal{R}_{01}$ &
 $\mathcal{R}_{0,pulse}^{\tilde{G}_{e}}$& $\tilde{\mathcal{R}}_{0,\mathcal{R}_{01}}$ &    $\mathcal{R}_{0,pulse}^{*}$ & $\mathcal{R}_{0,stable}^{*}$ & $\mathcal{R}_{0,stable}^{**}$& \bf{Case}\\
 \hline\hline
 \bf{a, b} & $>1$ & $ >1$ & $ >1$ & $ >1$ & $ >1$&$ >1$ & \bf{XIII} \\
 \hline
 \bf{c} & $ >1$ & $ >1$ & $ <1$ & $ >1$ & $ >1$ & $ >1$ & \bf{VIII} \\
 \hline
 \bf{d} & $<1$ & $>1$ & $ - $ & $ >1$ & $ -$& $-$ & \bf{V} \\
 \hline\hline
 \end{tabular}
 \label{R2b}
 \end{center}
 \end{table}
 \begin{figure}[H]
 \begin{center}
\includegraphics[scale=0.5]{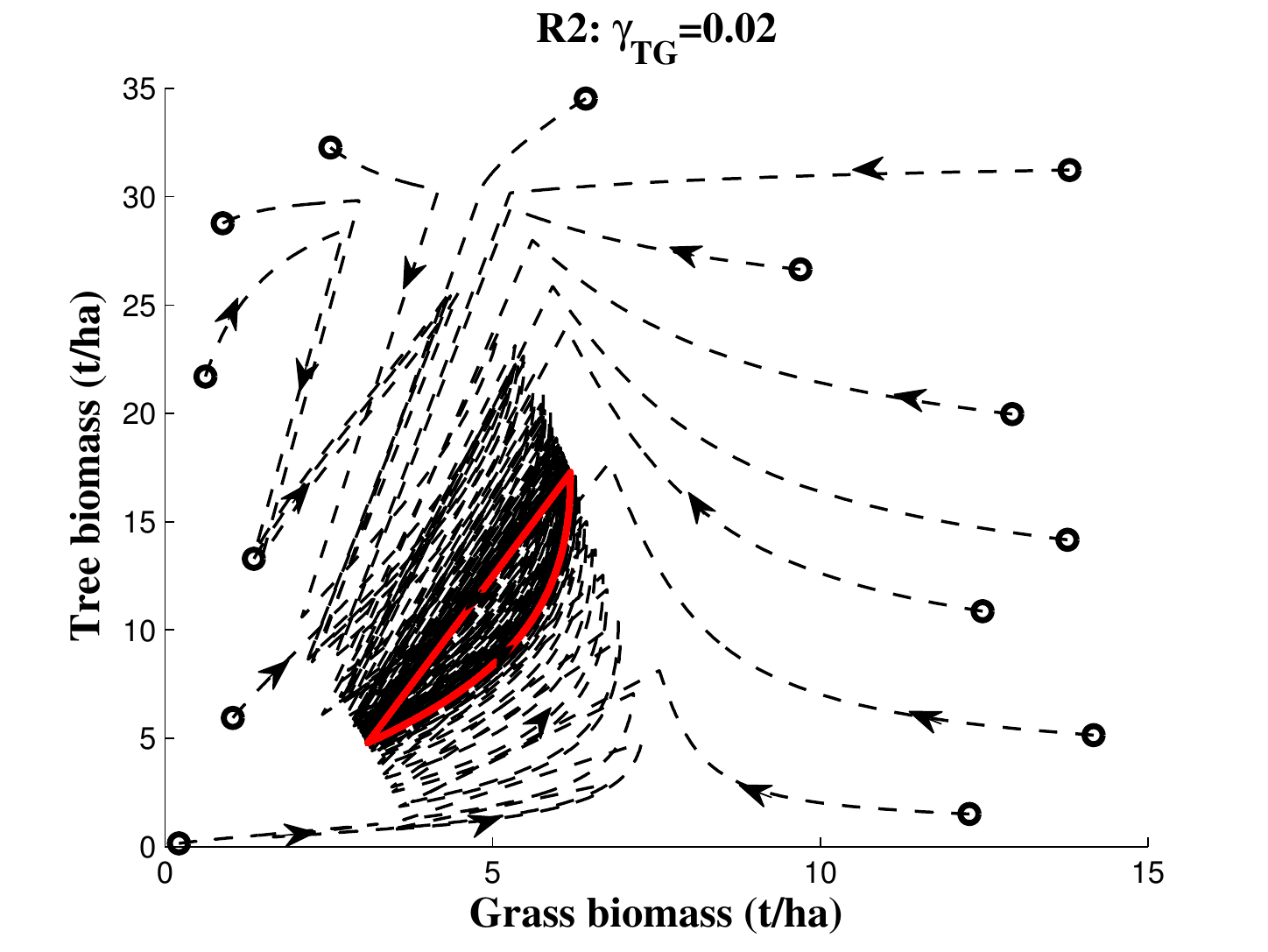}
\includegraphics[scale=0.5]{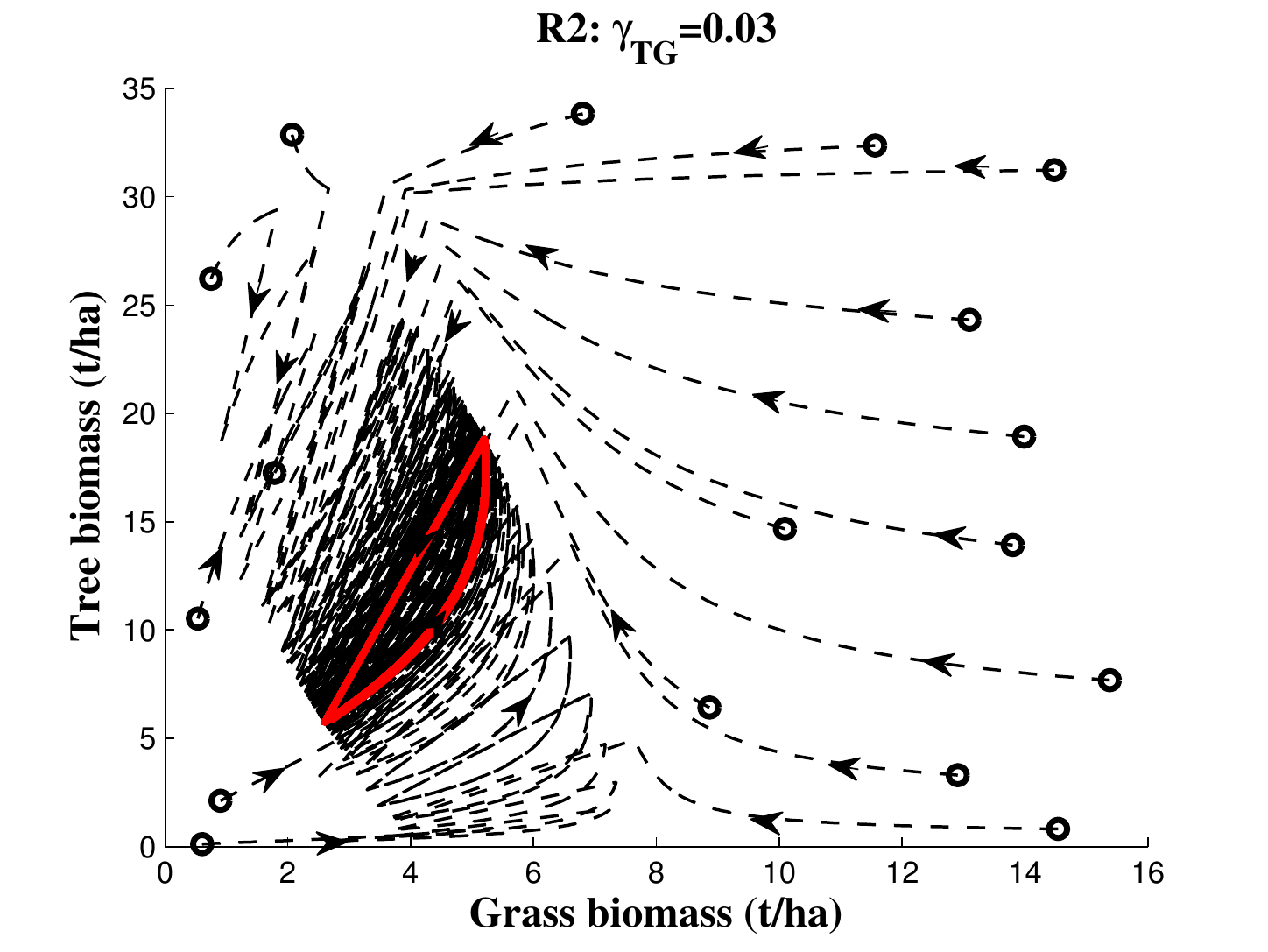}
\vspace{1cm}
\includegraphics[scale=0.5]{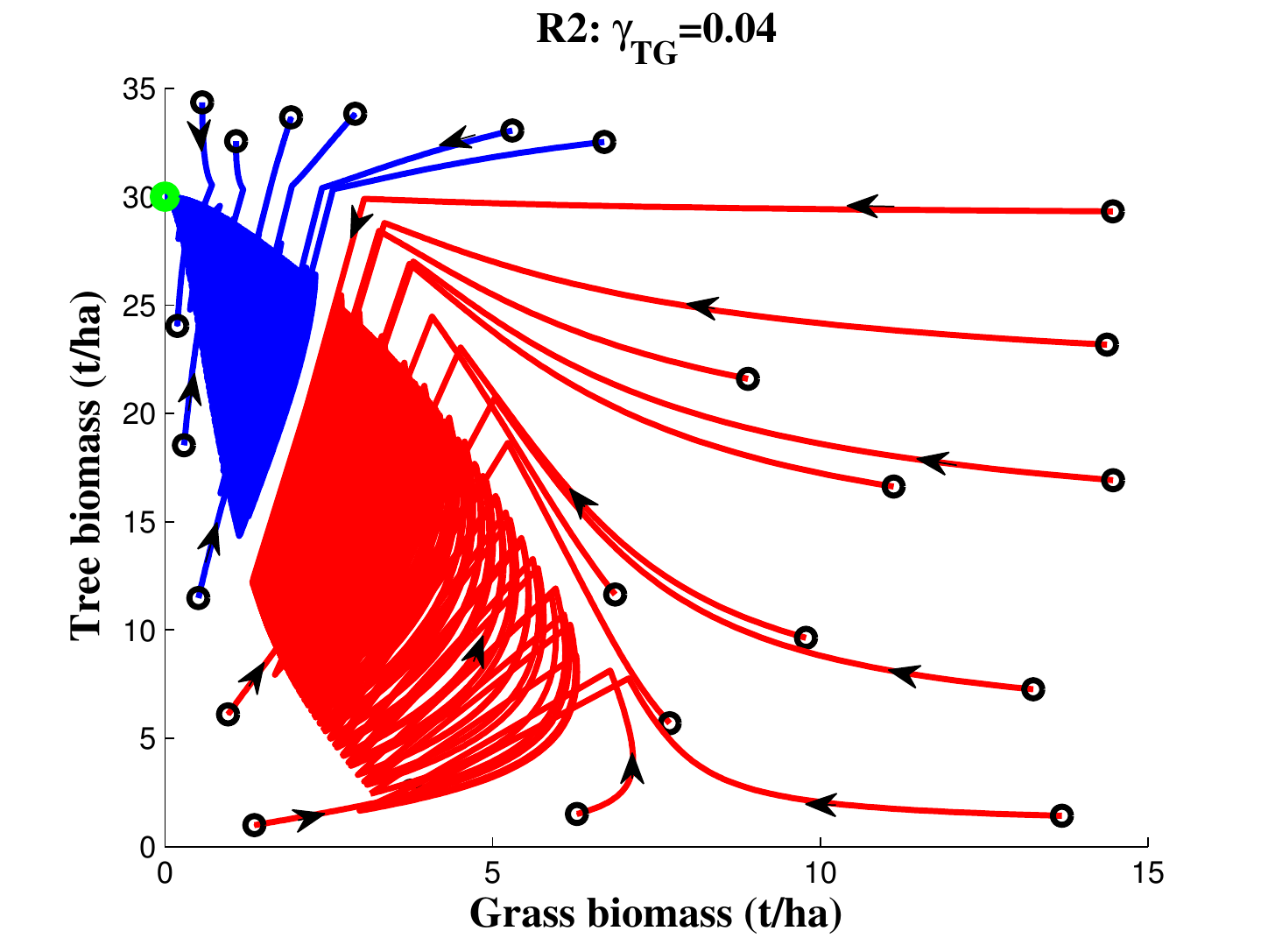}
\includegraphics[scale=0.5]{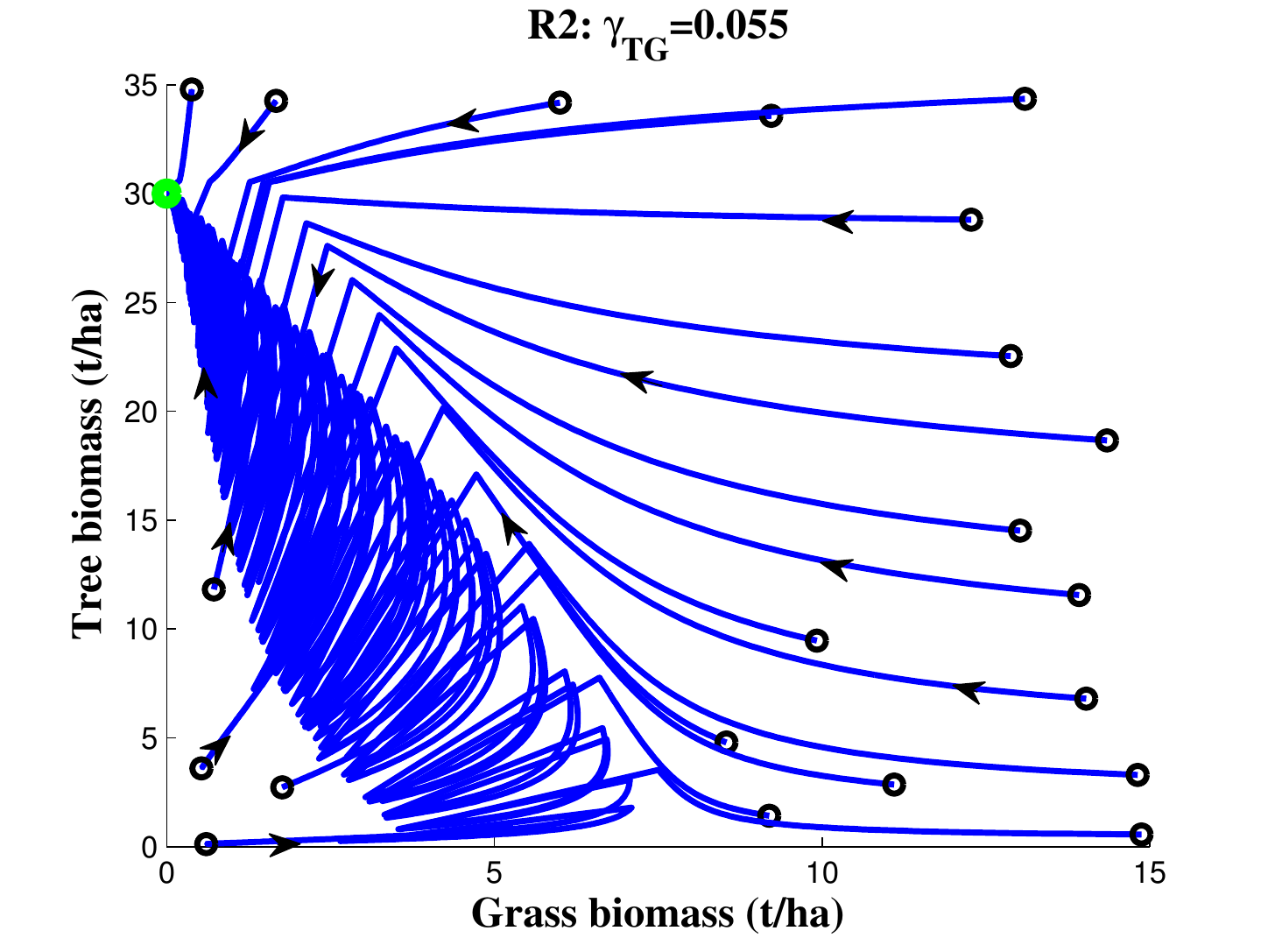}
 \caption{Phase diagrams in R2 with a period of fire of $\tau=5$ years (Table \ref{params_R2b}) and a stronger impact of fire than in Figure 3 via $\lambda_{fT}=0.8$. Panels $\mathbf{a}$, $\mathbf{b}$, $\mathbf{c}$ and $\mathbf{d}$ show the impact of the tree/grass competition parameter.}
  \label{figR2_b}
 \end{center}
 \end{figure}

Figure \ref{figR2_b} shows interesting behaviours. In particular, in panel $\mathbf{c}$, we obtain a bistable situation, where the system can either converge to the forest equilibrium or the periodic Tree-Grass equilibrium, depending on the initial values. In that case, we don't have analytic results that allow us to know what are the basins of attraction of each equilibrium. That is why the use of a well fitted numerical scheme is of utmost importance, in order to capture this essential information. We will show other examples of bistability in the next section.

\subsubsection{Simulations in region R3}

According to table  \ref{3_regions}, we first consider the following values for simulations in region R3:

\begin{table}[H]	
 	\caption{Parameters values related to figure \ref{figR3_a}}
\begin{center}
\renewcommand{\arraystretch}{1.5}
\begin{tabular}{cccccccccccc}
\hline
$K_{G}$ & $\gamma_{G}$ & $\delta_{G0}$ &  $K_{T}$&  $\gamma_{T}$ & $\delta_{T}$& $\alpha$ & $\tau$ & $\lambda_{fT}$ & $\lambda_{fG}$\\
\hline\hline
$17$ & $4.5$ & $0$ & $45$& $6$ & $0$ & $2$ & $0.6$ &$0.4$ & $0.4$\\
\hline\hline
\end{tabular}
\label{params_R3}
\end{center}
\end{table}

\begin{figure}[H]
\begin{center}
\includegraphics[scale=0.5]{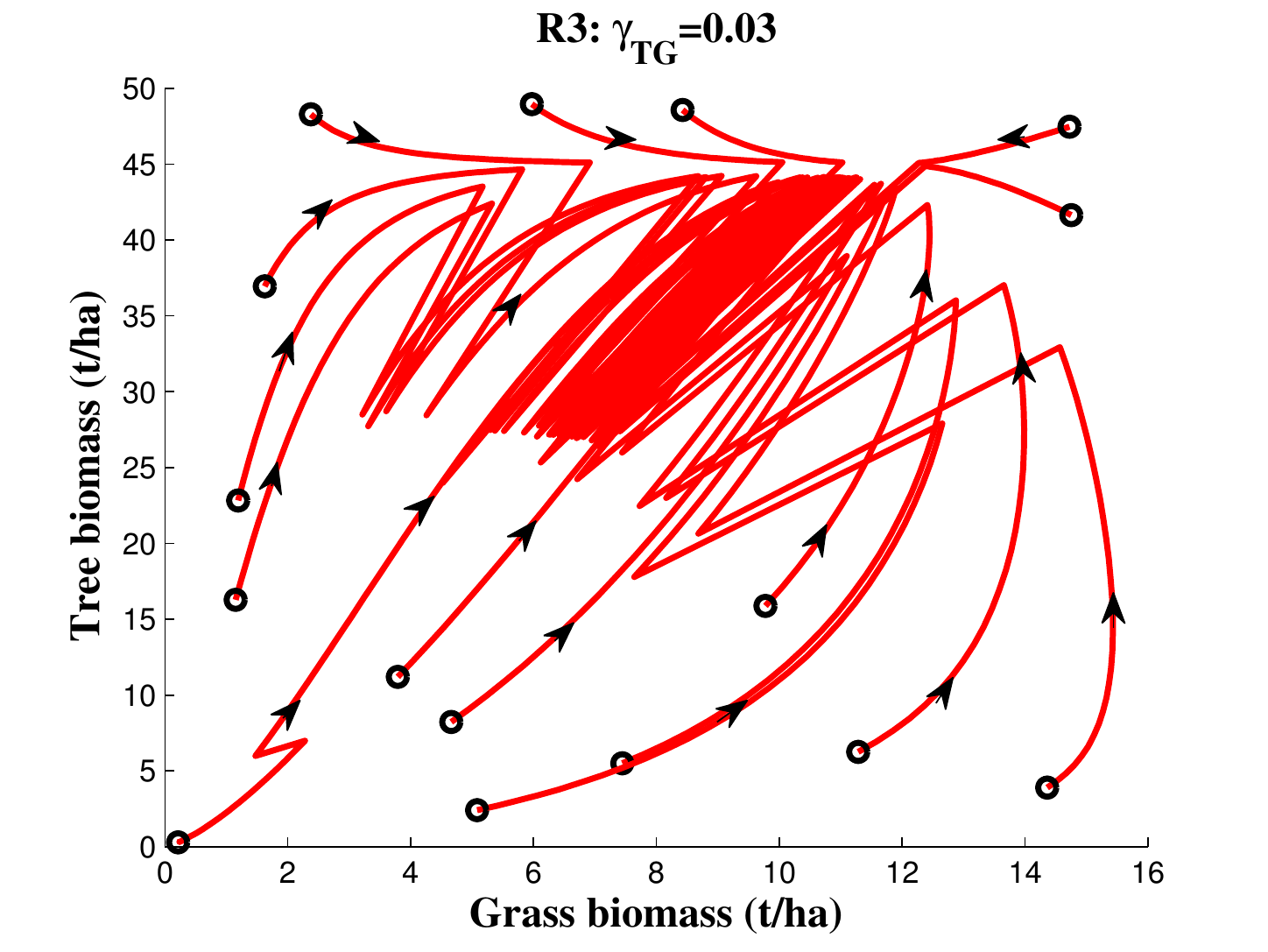}
\includegraphics[scale=0.5]{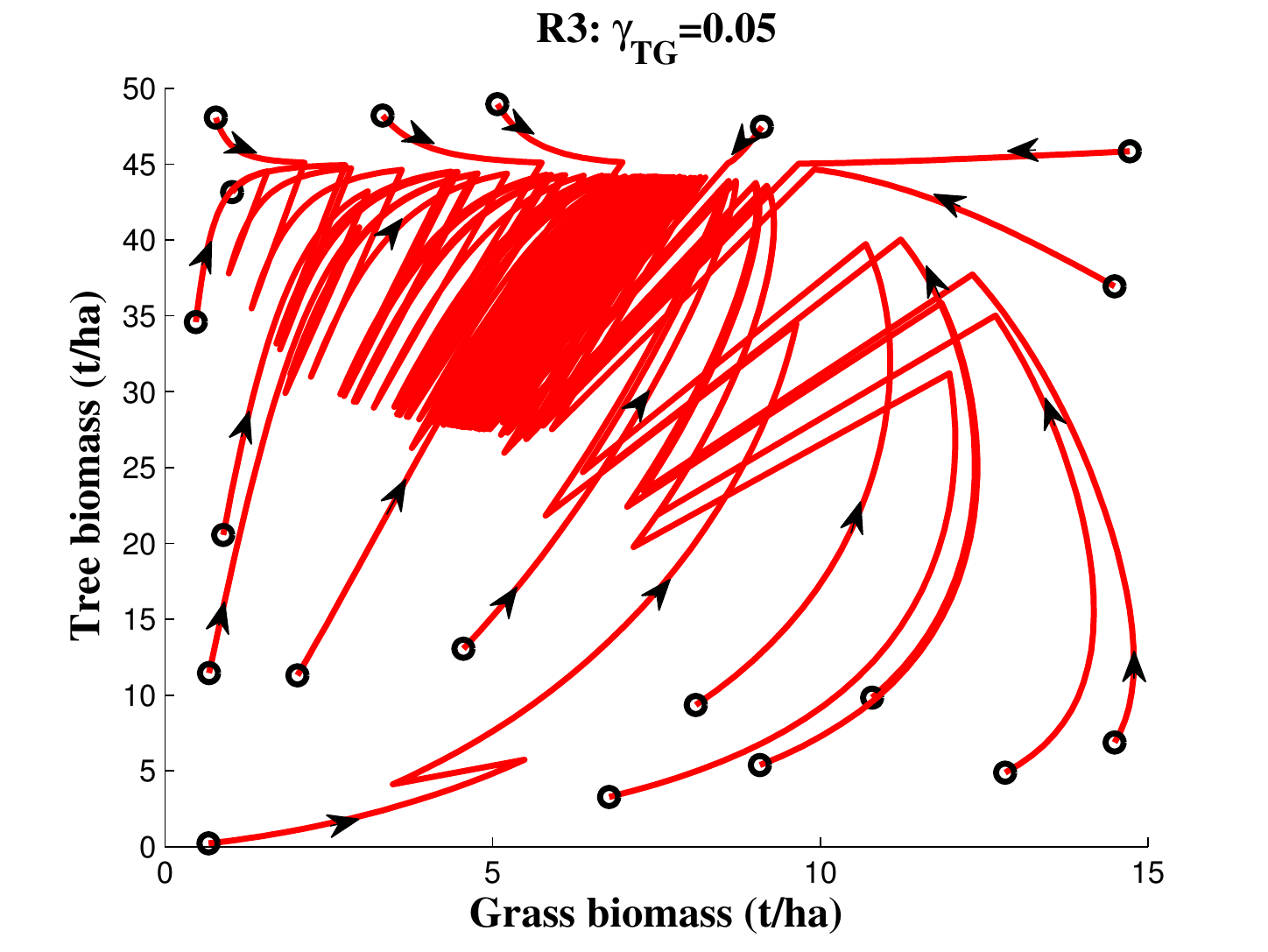}
\vspace{1cm}
\includegraphics[scale=0.5]{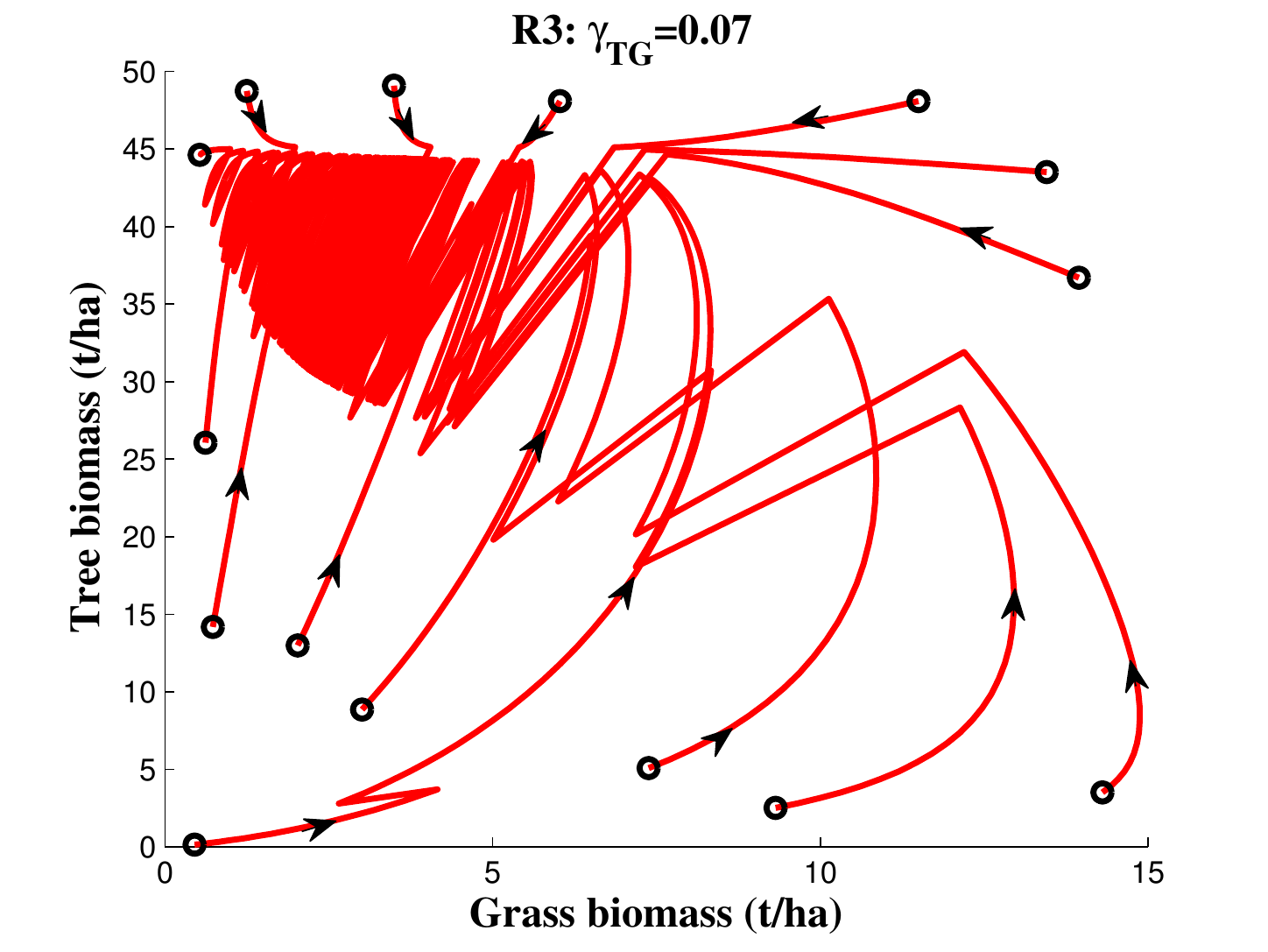}
\includegraphics[scale=0.5]{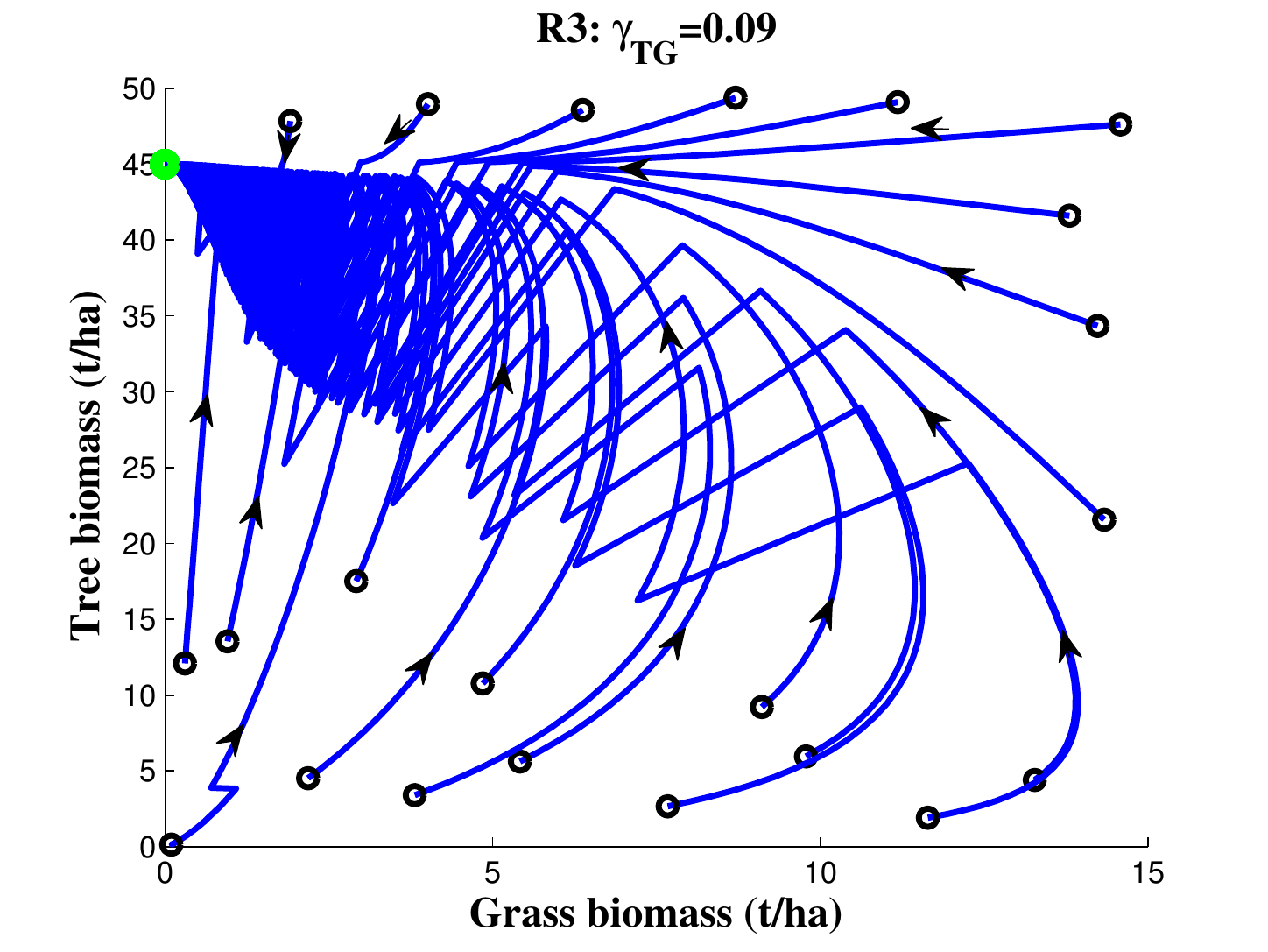}
\caption{\footnotesize {Phase diagrams in R3. This figure shows two equilibria: forest and periodic savanna.}} 
\label{figR3_a}
\end{center}
\end{figure}

Using Table \ref{params_R3_b} with $\gamma_{TG}=[0.03,0.05,0.07,0.09]$ leads to Table \ref{R3} and figure )\ref{figR3_a}.

\begin{table}[H]
\caption{Thresholds Table related to Table \ref{params_R3_b} and Figure \ref{figR3_a}}	
\begin{center}
\renewcommand{\arraystretch}{1.5}
\begin{tabular}{ccccccccc}
\hline
\bf{Panel} & $\mathcal{R}_{01}$ &
$\mathcal{R}_{0,pulse}^{\tilde{G}_{e}}$& $\tilde{\mathcal{R}}_{0,\mathcal{R}_{01}}$ &  
$\mathcal{R}_{0,pulse}^{*}$ & $\mathcal{R}_{0,stable}^{*}$&  $\mathcal{R}_{0,stable}^{**}$ &\bf{Case}\\
\hline\hline
\bf{a,b,c} & $>1$ & $ >1$ & $>1$ & $>1$ & $>1$ & $>1$& \bf{XIII}\\
\hline 
\bf{d} &$>1$ & $>1$ & $ <1$ & $>1$ & $-$ & $-$ & \bf{X} \\
\hline\hline
\end{tabular}
\label{R3}
\end{center}
\end{table}

In the humid zone, the vegetation is intrinsically  dominated by trees which exert competitive pressure on grasses, such that grasses are suppressed or even out-competed (Scholes and Walker 1993 \cite{Scholes1993african}, see panel \textbf{d}). Figure \ref{figR3_a} illustrates the effect of the competition parameter in R3.\par

\begin{note} Using realistic ranges for parameters, we show that in R3, a stable periodic savanna equilibrium may appear but also a stable forest equilibrium, for sufficiently high values of $\gamma_{TG}$. However we
cannot have a periodic grassland equilibrium. In R2, it is possible
to have forest equilibrium and periodic savanna. 
\end{note}

Like in Region R2, bistable situations can occur in region R3. Let us first consider the following parameter values

\begin{table}[H]	
 	\caption{Parameters values related to figure \ref{fig2}}
\begin{center}
\renewcommand{\arraystretch}{1.5}
\begin{tabular}{cccccccccccc}
\hline
$K_{G}$ & $\gamma_{G}$ & $\delta_{G0}$ & $\lambda_{fG}$ &  $K_{T}$&  $\gamma_{T}$ & $\delta_{T}$& $\lambda_{fT}$&$\tau$ & $\alpha$ &$\gamma_{TG}$\\
\hline\hline
$19$ & $3.1$ & $0.1$ & $0.5$& $50$ & $1.5$ & $0.015$ & $0.6$ &$0.5$ & $2$ & $0.09$\\
\hline\hline
\end{tabular}
\label{params_R3_a}
\end{center}
\end{table}

\begin{table}[H]	
\caption{Threshold values related to figure \ref{fig2}}
\begin{center}
\renewcommand{\arraystretch}{1.5}
\begin{tabular}{ccccccccc}
\hline
$\mathcal{R}_{01}$ & 
$\mathcal{R}_{0,pulse}^{\tilde{G}_{e}}$&
$\tilde{\mathcal{R}}_{0,\mathcal{R}_{01}}$ &  
$\mathcal{R}_{0,pulse}^{*}$ & $\mathcal{R}_{0,stable}^{*}$& $\mathcal{R}_{0,stable}^{**}$ &\bf{Case}\\
\hline\hline
$ <1$ & $>1$ & $ - $ & $<1$ &$-$ & $-$ & \bf{II} \\
\hline\hline
\end{tabular}
\label{diag_R3_1}
\end{center}
\end{table}

\begin{figure}[H]
\begin{center}
\includegraphics[scale=0.5]{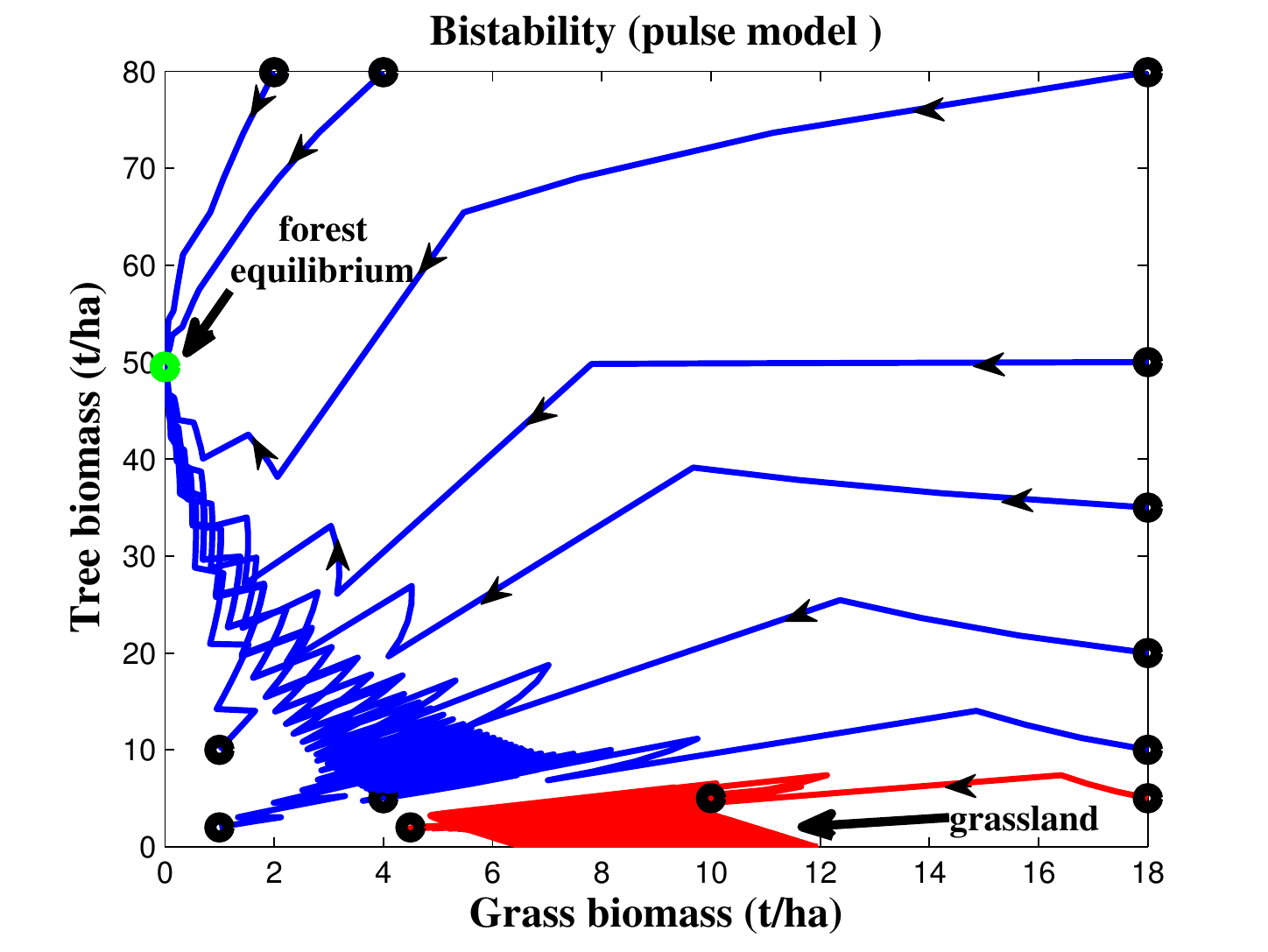}
\includegraphics[scale=0.5]{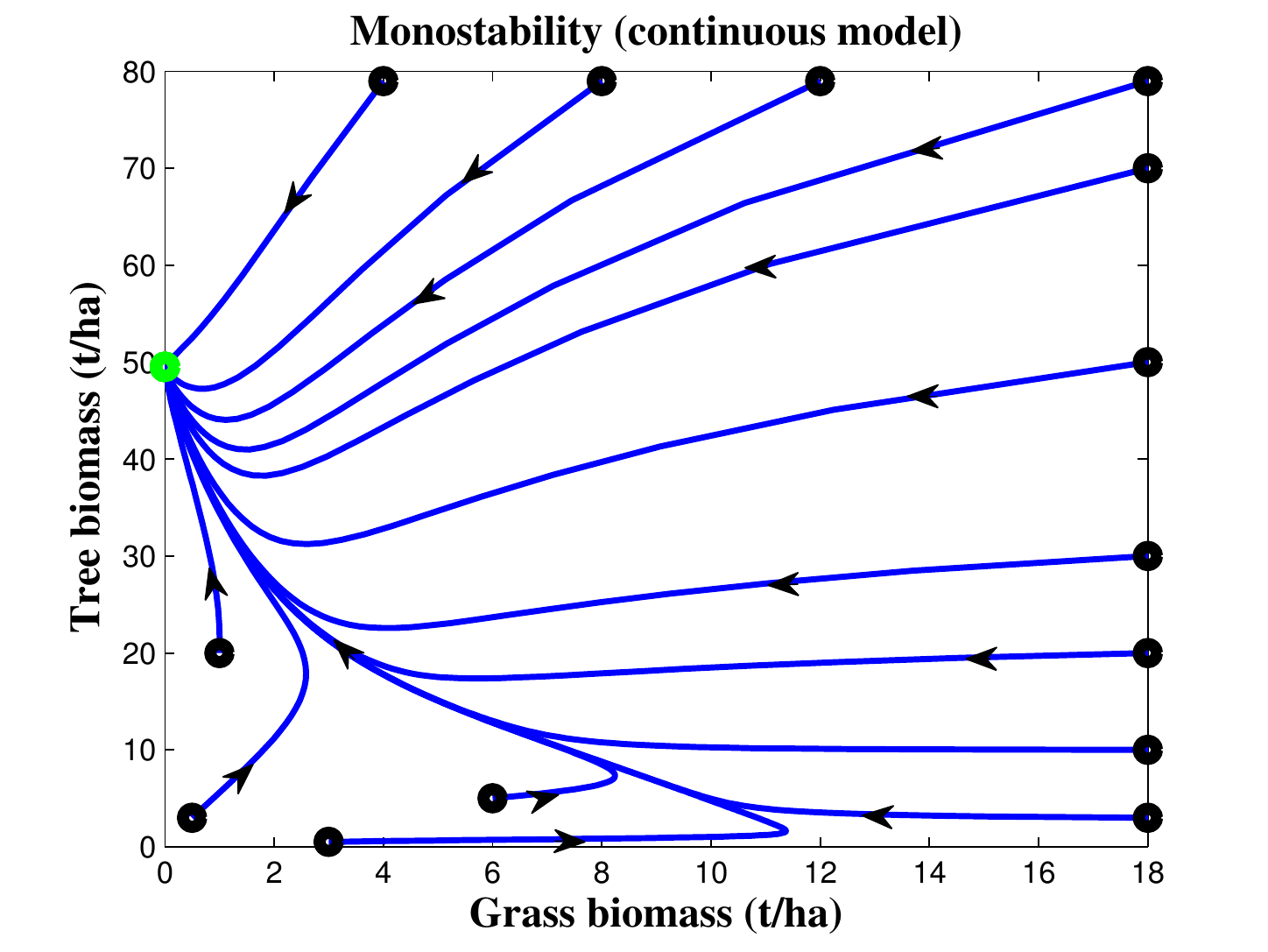}
\caption{\footnotesize {Comparison of the pulse vs. continuous models in reference to R3. The continuous
model presents a forest equilibrium which is GAS (see panel
$\mathbf{b}$). In contrast, the discrete model shows two 
equilibria: the forest (as for the continuous model) and the
periodic grassland (see panel $\mathbf{a}$). Depending of the
initial conditions the system converges to a stable periodic
grassland or to a stable forest equilibrium.}} \label{fig2}
\end{center}
\end{figure}

The mathematical analysis shows that there is a bistability between the forest equilibrium and a periodic grassland equilibrium (see line $\mathbf{II}$ in table \ref{tabfinal}). Panel $\mathbf{a}$ (pulse model) in figure \ref{fig2} illustrates two basins of attraction: one in favour of forest equilibrium;  another in favour of the periodic grassland equilibrium (bistability as in line $\mathbf{II}$ in table \ref{tabfinal}). For the same values of parameters,  the continuous model does not yield bistability (see panel $\mathbf{b}$ in \ref{fig2}): the only equilibrium is the forest equilibrium \cite{Tchuinte2014}.

We have shown that the interspecific parameter, $\gamma_{TG}$ plays a great role in the dynamics. However, fire period can also have an impact on the dynamics of the Tree-Grass system in R3 as in the two previous zones. Let us consider the following parameter values:

\begin{table}[H]
 	\caption{Parameters values related to figure \ref{BifR3}}	
\begin{center}
\renewcommand{\arraystretch}{1.5}
\begin{tabular}{cccccccccccc}
\hline
$K_{G}$ & $\gamma_{G}$ & $\delta_{G0}$ & $\lambda_{fG}$ &  $K_{T}$&  $\gamma_{T}$ & $\delta_{T}$& $\lambda_{fT}$ & $\alpha$ &$\gamma_{TG}$\\
\hline\hline
$19$ & $3.1$ & $0.1$ & $0.5$& $65$ & $1.5$ & $0.015$ & $0.6$  & $2$ & $0.04$\\
\hline\hline
\end{tabular}
\label{params_R3_b}
\end{center}
\end{table}
According to Tables \ref{tabfinal} and \ref{params_R3_b}, we derive, in Table \ref{diag_R3_2}, the different possible dynamics ot the Tree-Grass system:
\begin{table}[H]	
\caption{Threshold values related to figure \ref{BifR3}	}
\begin{center}
\renewcommand{\arraystretch}{1.5}
\begin{tabular}{cccccccccc}
\hline
\bf{Panel} & $\mathcal{R}_{01}$ &
$\mathcal{R}_{0,pulse}^{\tilde{G}_{e}}$& $\tilde{\mathcal{R}}_{0,\mathcal{R}_{01}}$ &  
$\mathcal{R}_{0,pulse}^{*}$ & $\mathcal{R}_{0,stable}^{*}$ & $\mathcal{R}_{0,stable}^{**}$ & \bf{Case}\\
\hline\hline
\bf{a} & $>1$ &  $>1$  & $<1$ & $<1$ & $-$ & $-$&  \bf{VII} \\
\hline 
\bf{b} & $>1$ & $ >1$ & $<1$ & $>1$ & $>1$ & $>1$ & \bf{VIII} \\
\hline\hline
\end{tabular}
\label{diag_R3_2}
\end{center}
\end{table}

\begin{figure}[h]
\begin{center}
\includegraphics[scale=0.4]{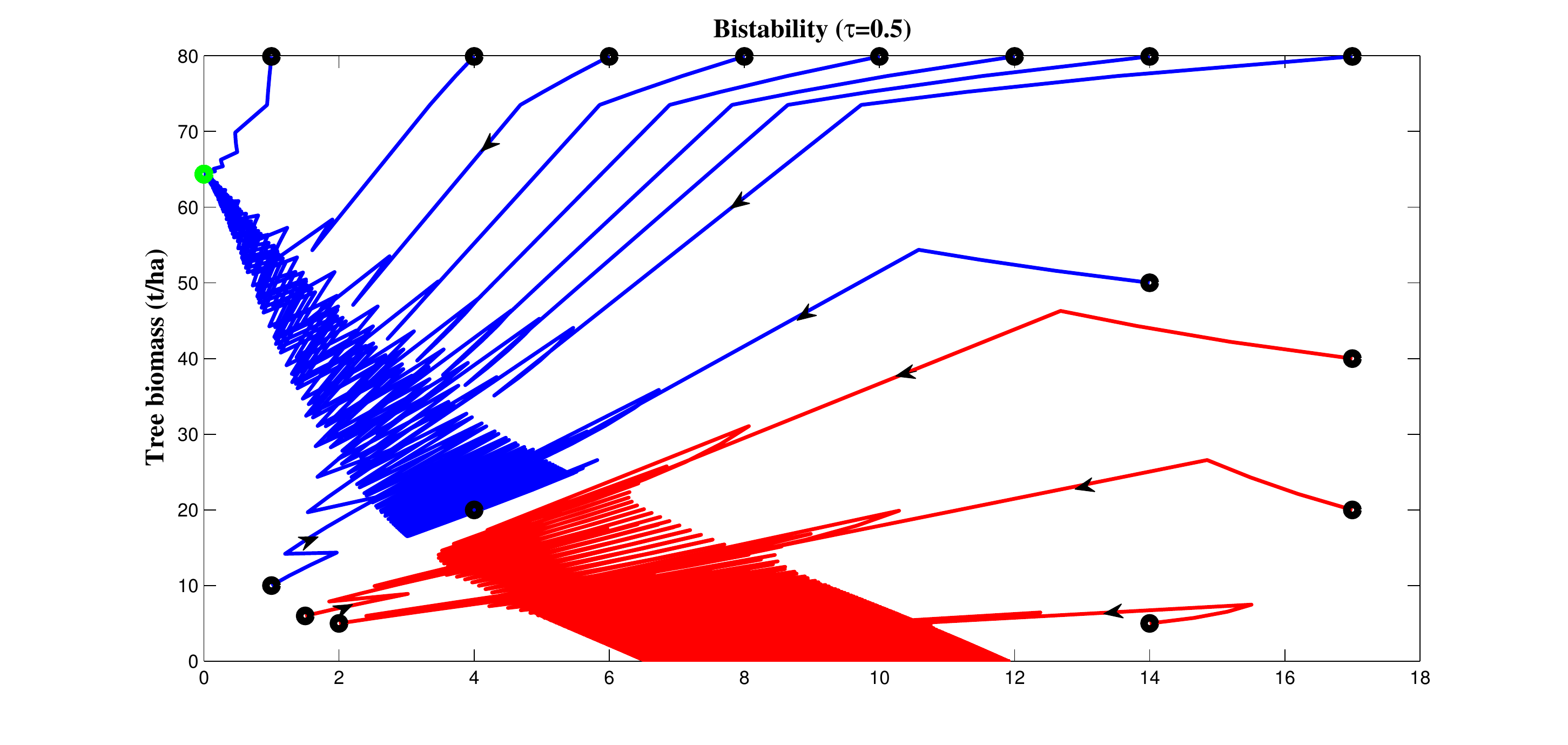}
\includegraphics[scale=0.4]{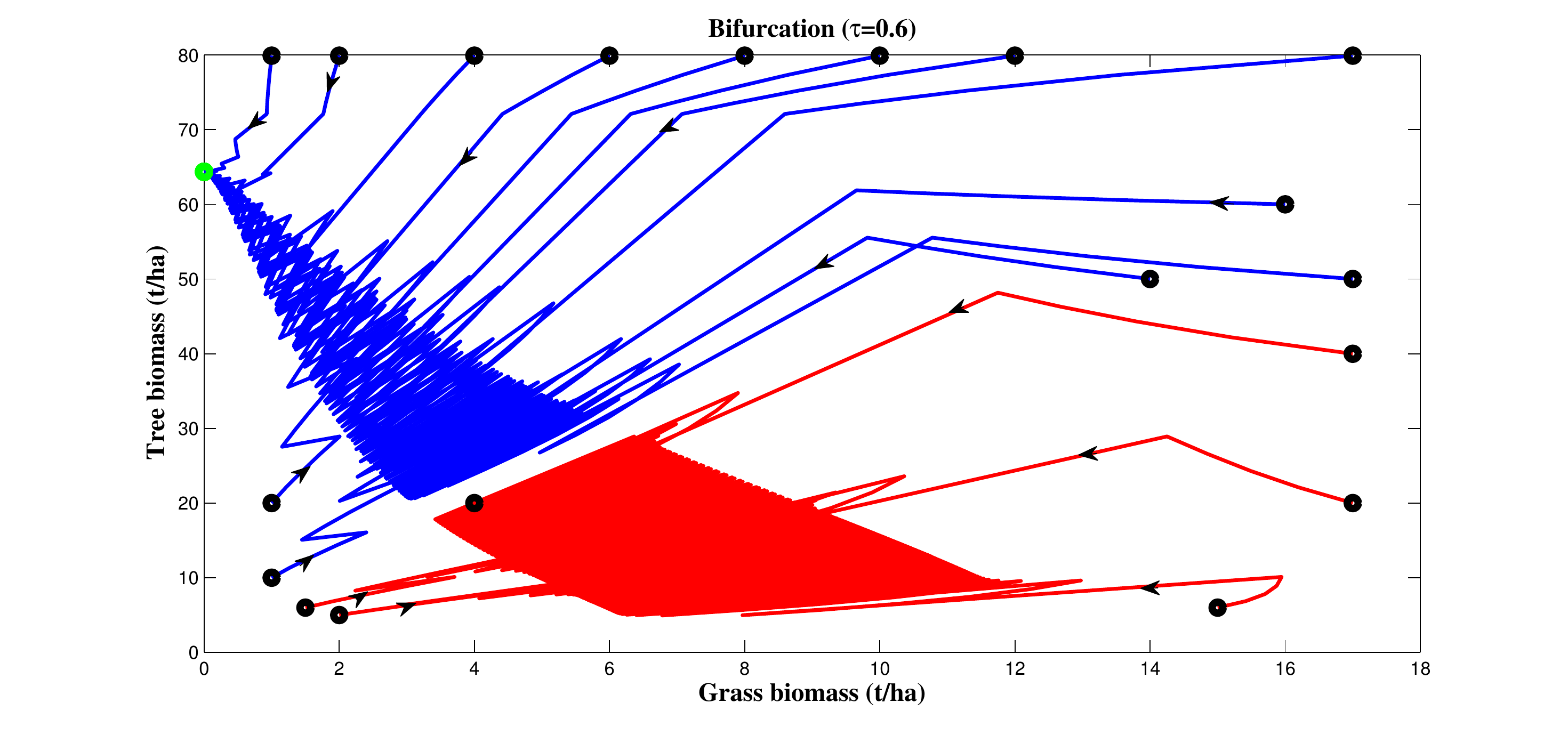}
\caption{Fire mediated Tree-Grass mixtures and changes
in their physiognomies according to the fire period. Panel $\mathbf{a}$ shows two bistable equilibria (a forest equilibrium and a periodic grassland
equilibrium (see $\mathbf{VII}$ in table \ref{tabfinal})). From $\tau=0.5$
to $0.6$, the ecosystem changes. Panel $\mathbf{b}$ illustrates  two
bistable equilibria (the periodic savanna and the forest equilibrium
(see $\mathbf{VIII}$ in table \ref{tabfinal})).}
 \label{BifR3}
\end{center}
\end{figure}

It has been evidenced that in humid savannas, fire is necessary to establish the Tree-Grass coexistence equilibrium (Sankaran et al. 2005 \cite{Sankaran2005determinants}). However, various bistabilities occur: between forest and grassland (see panel $\mathbf{a}$ in figure \ref{BifR3}); between forest and savanna (see panel $\mathbf{b}$ in figure \ref{BifR3}). The system can shift from $\mathbf{a}$ to $\mathbf{b}$ (bifurcation)  when the fire period increases. When $\tau=0.5$ (two fires a year in a sub-equatorial context), we have a bistability case: with a forest equilibrium and a periodic grassland equilibrium (see panel $\mathbf{a}$ in figure \ref{BifR3} which illustrates case $\mathbf{VII}$ in table \ref{tabfinal}). When the fire period increases by $10\%$ i.e. $\tau=0.6$, there is still bistability case, but the periodic grassland equilibrium is replaced by a periodic savanna equilibrium. This corresponds to line $\mathbf{VIII}$ in table \ref{tabfinal} and to panel $\mathbf{b}$ in figure \ref{BifR3}.

In R3, both forest-grassland, and forest-savanna equilibria are predicted (see  panels $\mathbf{a}$ and $\mathbf{b}$ in Figure
	\ref{BifR3}).  The fire
	period is a bifurcation parameter that shapes the
	Tree-Grass dynamics in R3.

\section{Conclusion}
Savannas are complex systems due to the interaction of trees and grasses which are frequently mediated by disturbances and notably by fires.  The
broad objective of this work was to develop a predictive
understanding of Tree-Grass dynamics across rainfall gradients in
Africa on the basis of a minimalistic model. This is done using specific features of three ecological contexts:
semi-arid, mesic and humid. They represent different ecological
conditions in terms of  rainfall amount and deriving variations of most of the  parameters used in the model. We formalize a
new model of Tree-Grass interactions. The novelty of the
model, with respect to other models (Staver et al., 2011
\cite{Staver2011tree}, Accatino et al., 2010
\cite{Accatino2010tree}, \cite{Tchuinte2014}) is that fire is considered as discrete
events with high or low return times. In addition, fire frequencies
and fire intensities are decoupled. Discrete events are typically
modelled by impulsive differential equations. We show here that this framework yields richer qualitative behaviours than continuous modelling (see figure \ref{fig2}). Several
authors, in order to deal with the stochastic occurrence of fire
have  modelled fire through purely stochastic differential
equations (D'Odoricco et al. 2006 \cite{DOdorico2006probabilistic},
Beckage et al. 2011 \cite{Beckage2011grass}, De Michele and Accatino 2014 \cite{DeMichele2014tree}). However results of those stochastic  models are often obtained numerically by iterating the equations given parameter values and initial conditions and the authors have difficulties to verify mathematically that the qualitative properties (e.g existence of equilibria and their stability) of the model are preserved in their simulations. By contrast,
at least qualitatively, our impulsive model is lend itself to a comprehensive qualitative analysis of the possible dynamic outcomes and properties of the system.
For this raison, it is perhaps interesting to used impulsive framework to look  at least qualitatively at the problem of predictability of discrete fire impacts in tropical savannas. The impulsive modelling of fire suggests ways for deriving from a
minimal continuous fire model (e.g. Tchuinte et al., 2014
\cite{Tchuinte2014}) more realistic discrete fire model. The
theory of IDE allows for more detailed analyses  of the system than
stochastic differential equations. Further, using IDE technique, we
can highlight mathematically thresholds that summarize the
Tree-Grass interaction (see table \ref{tabfinal}) and point out
bifurcation parameters.\par

The impulsive model
illustrates  different kinds of dynamics which can be observed across
of  rainfall gradients in Africa. The model generates savanna
equilibria in all the three regions. It means that trees are able to persist while not reaching $80$-$100\%$ cover all over the wide range of rainfall considered. However,
at  high rainfall sites such like the boundaries of  the tropical
rainforests of central Africa, the vegetation is due to be dominated by trees and may even reach a closed tree cover in the absence of recurrent fires of low return times which are fostered by high grass production. On the other hand, trees can reach high individual biomass and exert intense light
competition on grasses, such that grasses are  suppressed or
even out-competed (Scholes and Walker 1993
\cite{Scholes1993african}). This explains why forest may be stable in R3 and even in R2. However, tree cover can also facilitate the growth of grasses, i.e. $\gamma_{TG}<0$ in arid conditions (i.e. R1) by limiting soil water losses from transpiration in the topsoil. In R1, we show that the existence of savanna is mainly due to the competition parameter. Thus, only two equilibria are possible: a forest equilibrium and a periodic savanna equilibrium.  The bifurcation from the forest to the periodic
savanna is related to the resource competition parameter, $\gamma_{TG}$. This result joins those of  Sankaran et al., (2005)
\cite{Sankaran2005determinants} that argued that in drier sites,
savannas are stable in the sense that tree cover is intrinsically limited by resource (and this contributes to the low value of KT) and fire is not necessary for Tree-Grass coexistence. Recent modelling studies
by Baudena et al. (2014) \cite{Baudena2014forests} confirm that in
semi-arid savannas, while trees are water-limited,  the water
competition with grasses is also a key factor determining savanna
existence. 

We show that in all regions, the competition parameter plays a great role: it is a bifurcation parameter whatever the context.\par

Grass-fire feedback principally occurs in mesic and humid areas. No savanna or grassland would emerge without this positive feedback in Regions R2 and R3. Using realistic parameter values, we show that in R2 (mesic area) forest and savanna are more present than grassland. In R3, in humid area, two
bistability  situations may occur: bistability between
forest and grassland and bistability between forest and 
(periodic) savanna. In the equatorial climate in Southern Cameroon (Central
Africa), there is forest-savanna contact in three sites inside R3:
Akonolinga, Bertoua and Mbam-Kim (Youta 1998
\cite{Happi1998arbres}).
Bistability of forest and grassland
can be found in the equatorial and tropical climate of transitions
($1100$-$1500$ mm/yr) where fires are usually occurring every
$1$-$5$ years (Frost and Robertson, 1987 \cite{Frost1987fire}, Favier et al. 2012 \cite{Favier2012abrupt}). In
this area and due to the large grass biomass, flame height is usually $2$-$3$ m high (Frost and
Robertson, 1987 \cite{Frost1987fire}). Therefore, a severe fire
could have a great impact on young trees/shrubs. In R3, water availability enables high fuel production. As a result,
fire is severe and may occur as frequently as $0.5$-$1$ yr. Thus, fire has
a stabilizing role of grassland and savannas by preventing tree invasion on long time scales,
freezing the forest-savanna boundary in a historical position
(Gillon, 1983 \cite{Gillon1983fire}). We show from  our model that in R3 the Tree-Grass system can
shift from bistability between forest and grassland to bistability between forest and savanna due to fire period. 

To conclude: (i) in all regions
Tree-Grass competition is the most important parameter for
Tree-Grass co-existence; (ii) in R2 and R3, fire period can also be a
bifurcation parameter. Thus, the fire period (and fire intensities) and the competition parameters are the main determinant in the Tree-Grass dynamics; (iii) Modelling fire as pulse events provides more realistic situations than modelling fire as continuous events.

%



\bibliography{Biblio}

\newpage

\section*{APPENDICES}

 Some detail of the proofs of results associated with system
 (\ref{Impuleq1}) are provided.
\section*{Appendix A: Proof of Lemma \ref{lem_compact}}
\label{AppendixA}
It is obvious that $G=0$, and $T=0$ are vertical and horizontal
null-clines respectively. Then, no trajectory can cut these axes.
Thus, the positive cone $\mathbf{R}_{+}^{2}$ is positively invariant
for $(\ref{Impuleq1})$ because, all trajectories that start in
$\mathbf{R}^{2}_{+}$ remain in $\mathbf{R}^{2}_{+}$ for all positive
time. From system $(\ref{Impuleq1})$, with the initial conditions
$T(t_{0})=T_{0}>0$ and $G(t_{0})=G_{0}>0$, we obtain the following
system

\begin{equation}
\label{Impuleq13} \left\{
\begin{array}{lcl}
\displaystyle\frac{dG}{dt}\leq(\gamma_{G}-\delta_{G0})G-\mu_{G}G^{2},\\
\\
\displaystyle\frac{dT}{dt}\leq(\gamma_{T}-\delta_{T})T-\mu_{T}T^{2},\\
\\
T(t_{0})=T_{0},\\
\\
G(t_{0})=G_{0}.
\end{array}
\right.
\end{equation}

Using the maximum principle, we deduce that

\begin{equation}
\label{Impuleq14} \left\{
\begin{array}{lcl}
G\leq\displaystyle\frac{G_{0}}{\displaystyle\frac{G_{0}}{X_{G}}+\left(1-\displaystyle\frac{G_{0}}{X_{G}}\right)\exp\left\{-X_{G}\mu_{G}t\right\}},\\
\\
T\leq\displaystyle\frac{T_{0}}{\displaystyle\frac{T_{0}}{Y_{T}}+\left(1-\displaystyle\frac{T_{0}}{Y_{T}}\right)\exp\left\{-Y_{T}\mu_{T}t\right\}},\\
\end{array}
\right.
\end{equation}

 When
$t\rightarrow\infty$, we obtain

\begin{equation}
\label{Impuleq15} \left\{
\begin{array}{lcl}
\lim\limits_{t\rightarrow\infty}G(t)\leq X_{G}=\displaystyle\frac{\gamma_{G}-\delta_{G0}}{\mu_{G}}=K_{G}\left(1-\displaystyle\frac{\delta_{G0}}{\gamma_{G}}\right),\\
\\
\lim\limits_{t\rightarrow\infty}T(t)\leq Y_{T}=\displaystyle\frac{\gamma_{T}-\delta_{T}}{\mu_{T}}=K_{T}\left(1-\displaystyle\frac{\delta_{T}}{\gamma_{T}}\right).\\
\end{array}
\right.
\end{equation}

Hence, when  $\gamma_{G}>\delta_{G0}$ and $\gamma_{T}>\delta_{T}$,
all trajectories of system $(\ref{Impuleq1})$ that reach the
neighbourhood of $\mathcal{B}$ converge inside  as $t$ tends to
infinity. Since $\mathcal{B}\subseteq \mathbf{R}_{+}^{2}$, then
$\mathcal{B}$  is positively invariant and attracting for system
$(\ref{Impuleq1})$.

\section*{Appendix B: Proof of theorem \ref{thm_grassland} (existence of the semi-trivial periodic equilibrium)} 
\label{AppendixB}

Let $T(t)\equiv 0$, from system (\ref{Impuleq1}), we have the
following simple logistic impulsive differential system:

 \begin{equation}
  \label{Impuleq19}
 \left\{
 \begin{array}{lcl}
\displaystyle \frac{dG}{dt}&=& r_{G}G-\mu_{G}G^{2},\hspace{0.5cm}t\neq t_{n},\\
\\
G(t_{n}^{+})&=& G(t_{n})-\lambda_{fG}G(t_{n})\hspace{0.5cm}t=t_{n}.
 \end{array}
 \right.
 \end{equation}

The solution of system (\ref{Impuleq19}) is given in \cite{Nie2010}.
Here, is shown in detail the proof.

Setting $X=\displaystyle\frac{1}{G}$ in the first equation of system
$(\ref{Impuleq19})$, we have the following differential equation
\begin{equation*}
\displaystyle\frac{dX}{dt}=-r_{G}X+\mu_{G}.
\end{equation*}
 Integrating
 \begin{equation*}
\displaystyle\frac{dX}{dt}=-r_{G}X
 \end{equation*}
 from $n\tau$ to $t$, we obtain
 \begin{equation}
 \label{Impuleq20}
 X(t)=ae^{-r_{G}(t-n\tau)},\hspace{0.5cm}a\in\mathbf{R}.
 \end{equation}
 Using the variation of the constant $a$, we have the following differential equation
 \begin{equation*}
 \displaystyle\frac{da}{dt}=\mu_{G}e^{r_{G}(t-n\tau)}.
 \end{equation*}
 Thus, we have
 \begin{equation}
 \label{Impuleq21}
 a(t)=\displaystyle\frac{1}{X_{G}}e^{r_{G}(t-n\tau)}+b,\hspace{0.5cm}b\in\mathbf{R}.
 \end{equation}
 Substituting $(\ref{Impuleq21})$ in $(\ref{Impuleq20})$, we obtain

  \begin{equation}
  \label{Impuleq22}
 X(t)=\displaystyle\frac{1}{X_{G}}+be^{-r_{G}(t-n\tau)},\hspace{0.5cm}b\in\mathbf{R}.
  \end{equation}
 Considering $X(t)=\displaystyle\frac{1}{G(t)}$ in $(\ref{Impuleq22})$, we have
 \begin{equation*}
 \displaystyle\frac{1}{G(t)}=\displaystyle\frac{1}{X_{G}}+be^{-r_{G}(t-n\tau)},\hspace{0.5cm}b\in\mathbf{R},
 \end{equation*}
 this implies that,

  \begin{equation}
  \label{Impuleq23}
 G(t)=\displaystyle\frac{1}{\displaystyle\frac{1}{X_{G}}+be^{-r_{G}(t-n\tau)}},\hspace{0.5cm}b\in\mathbf{R}.
  \end{equation}

At $t=n\tau$, we have
\begin{equation*}
G(n\tau)=\displaystyle\frac{1}{\displaystyle\frac{1}{X_{G}}+b},\hspace{0.5cm}b\in\mathbf{R},
\end{equation*}
and therefore
\begin{equation*}
b=\displaystyle\frac{1}{G(n\tau)}-\displaystyle\frac{1}{X_{G}}.
\end{equation*}
 Substituting the expression of $b$ in $(\ref{Impuleq23})$, we obtain
   \begin{equation}
   \label{Impuleq24}
  G(t)=\displaystyle\frac{G(n\tau)e^{r_{G}(t-n\tau)}}{1+\displaystyle\frac{G(n\tau)}{X_{G}}\left(e^{r_{G}(t-n\tau)}-1\right)}.
   \end{equation}
 We have
    \begin{equation}
    \label{Impuleq25}
   G(t_{n+1}^{+})=G((n+1)\tau)=\displaystyle\frac{G(n\tau)e^{r_{G}\tau}}{1+\displaystyle\frac{G(n\tau)}{X_{G}}\left(e^{r_{G}\tau}-1\right)}.
    \end{equation}
 By substituting $(\ref{Impuleq25})$ into the second equation of $(\ref{Impuleq19})$, we obtain

     \begin{equation*}
    G(n\tau)=(1-\lambda_{fG})G((n+1)\tau)=\displaystyle\frac{(1-\lambda_{fG})G(n\tau)e^{r_{G}\tau}}{1+\displaystyle\frac{G(n\tau)}{X_{G}}\left(e^{r_{G}\tau}-1\right)},
     \end{equation*}

which implies that
 \begin{equation}
     \label{Impuleq26}
    G(n\tau)=\displaystyle\frac{X_{G}\{(1-\lambda_{fG})e^{r_{G}\tau}-1\}}{e^{r_{G}\tau}-1}.
     \end{equation}

 Substituting the expression of $G(n\tau)$ into $(\ref{Impuleq24})$, we have

 \begin{align*}
 G(t)&=\displaystyle\frac{X_{G}\{(1-\lambda_{fG})e^{r_{G}\tau}-1\}e^{r_{G}(t-n\tau)}}{(e^{r_{G}\tau}-1)
\left\{1+\displaystyle\frac{\{(1-\lambda_{fG})e^{r_{G}\tau}-1\}(e^{r_{G}(t-n\tau)}-1)}{(e^{r_{G}\tau}-1)}\right\}}\\
&=\displaystyle\frac{X_{G}\{(1-\lambda_{fG})e^{r_{G}\tau}-1\}e^{r_{G}(t-n\tau)}}{(e^{r_{G}\tau}-1)+
\{(1-\lambda_{fG})e^{r_{G}\tau}-1\}(e^{r_{G}(t-n\tau)}-1)}\\
&=\displaystyle\frac{X_{G}\{(1-\lambda_{fG})e^{r_{G}\tau}-1\}e^{r_{G}(t-n\tau)}}{\{(1-\lambda_{fG})e^{r_{G}\tau}-1\}e^{r_{G}(t-n\tau)}+\lambda_{fG}e^{r_{G}\tau}}.
 \end{align*}

Thus, when $\mathcal{R}_{0,pulse}^{\tilde{G}_{e}}>1$, there exists

\begin{equation}
\label{Impuleq27}
\tilde{G}_{e}(t)=X_{G}\displaystyle\frac{\{(1-\lambda_{fG})e^{r_{G}\tau}-1\}e^{r_{G}(t-n\tau)}}{\{(1-\lambda_{fG})e^{r_{G}\tau}-1\}e^{r_{G}(t-n\tau)}+\lambda_{fG}e^{r_{G}\tau}}>0,\hspace{0.5cm}t\in
[n\tau, (n+1)\tau[, n=0,1,2,....
\end{equation}
This completes the proof.

\section*{Appendix C: Proof of theorem \ref{thm11}.}
\label{AppendixC}
We calculate the unique periodic equilibrium. \\
Set $X=\displaystyle\frac{1}{T}$ in  the second equation of
$(\ref{Impuleq1})$. We obtain the following differential equation

    \begin{equation}
     \label{Impuleq45}
  \begin{array}{lcl}
 \displaystyle \frac{dX}{dt}=-r_{T}X+\mu_{T}.
  \end{array}
    \end{equation}

Integrate system (\ref{Impuleq45}) in $n\tau\leq t<(n+1)\tau$, we
have

     \begin{equation}
   \label{Impuleq50}
\begin{array}{lcl}
T(t)&=&\displaystyle\frac{T(n\tau)e^{r_{T}(t-n\tau)}}{1+\displaystyle\frac{T(n\tau)}{Y_{T}}\left[e^{r_{T}(t-n\tau)}-1\right]}.
\end{array}
   \end{equation}

Now, we solve the first equation of $(\ref{Impuleq1})$

  \begin{equation}
  \label{Impuleq51}
  \begin{array}{lcl}
  \displaystyle \frac{dG}{dt}=(\gamma_{G}-\delta_{G0})G-\mu_{G}G^{2}-\gamma_{TG}TG,
  \end{array}
  \end{equation}
where  the expression of $T$ is given by $(\ref{Impuleq50})$.\par
Set $Y=\displaystyle\frac{1}{G}$ in  $(\ref{Impuleq51})$. We have
the following differential equation

    \begin{equation}
    \label{Impuleq52}
 \begin{array}{lcl}
\displaystyle\frac{dY}{dt}&=&-(r_{G}-\gamma_{TG}T(t))Y+\mu_{G}.
 \end{array}
    \end{equation}
Integrating

      \begin{equation}
      \label{Impuleq53}
   \begin{array}{lcl}
  \displaystyle\frac{dY}{dt}&=&-(r_{G}-\gamma_{TG}T(t))Y,
   \end{array}
      \end{equation}

 from $n\tau$ to $t$, we have

  \begin{align*}
  \ln(Y)&=-r_{G}(t-n\tau)+\gamma_{TG}\int_{n\tau}^{t}T(u)du+a,~~a\in\mathbf{R}\\
        &=-r_{G}(t-n\tau)+\displaystyle\frac{\gamma_{TG}}{\mu_{T}}\ln\left[1+\displaystyle\frac{T(n\tau)}{Y_{T}}\left[e^{r_{T}(t-n\tau)}-1\right]\right]+a,~~a\in\mathbf{R},
  \end{align*}

  which implies that,

         \begin{equation}
         \label{Impuleq54}
      \begin{array}{lcl}
    Y(t)&=& Pe^{-m(t,n\tau,T(n\tau))},~~P\in\mathbf{R},
      \end{array}
         \end{equation}

 where,
\begin{equation}
 \label{Impuleq55}
 \begin{array}{lcl}
 m(t,n\tau,T(n\tau))&=&r_{G}(t-n\tau)+\displaystyle\frac{\gamma_{TG}}{\mu_{T}}\ln\left[\displaystyle\frac{1}{1+\displaystyle\frac{T(n\tau)}{Y_{T}}\left[e^{r_{T}(t-n\tau)}-1\right]}\right].
 \end{array}
 \end{equation}
Variation of $P$ gives

           \begin{equation}
            \label{Impuleq56}
         \begin{array}{lcl}
P(t)=\mu_{G}\int_{n\tau}^{t}\chi(u,n\tau,T(n\tau))du+b,~~b\in\mathbf{R}
         \end{array}
            \end{equation}

 where,

             \begin{equation}
             \label{Impuleq57}
          \begin{array}{lcl}
 \chi(t,n\tau,T(n\tau))=e^{m(t,n\tau,T(n\tau))}.
          \end{array}
             \end{equation}

  Then, we have

  \begin{displaymath}
  Y(t)=\displaystyle\frac{P(t)}{\chi(t,n\tau,T(n\tau))},
  \end{displaymath}

 which implies that

  \begin{equation}
  \label{Impuleq58}
  \begin{array}{lcl}
  G(t)&=&\displaystyle\frac{G(n\tau)\chi(t,n\tau,T(n\tau))}{1+\mu_{G}G(n\tau)\int_{n\tau}^{t}\chi(u,n\tau,T(n\tau))du}.
  \end{array}
  \end{equation}

From $(\ref{Impuleq50})$ and $(\ref{Impuleq58})$, we have

\begin{equation}
  \label{Impuleq59}
  \left\{
  \begin{array}{lcl}
  G(t)&=&\displaystyle\frac{G(n\tau)\chi(t,n\tau,T(n\tau))}{1+\mu_{G}G(n\tau)\int_{n\tau}^{t}\chi(u,n\tau,T(n\tau))du},\\
  \\
  T(t)&=&\displaystyle\frac{T(n\tau)e^{r_{T}(t-n\tau)}}{1+\displaystyle\frac{T(n\tau)}{Y_{T}}(e^{r_{T}(t-n\tau)}-1)},
  \end{array}
  \right.
\end{equation}

where $G(n\tau)$ and $T(n\tau)$ are values of grasses and trees
biomasses respectively, immediately after the $n^{th}$ pulse  of
fire at the time $n\tau$.  $G(n\tau)$ and $T(n\tau)$ can be viewed
as the initial values  of $(\ref{Impuleq1})$ in the interval
$[n\tau, (n+1)\tau[$. The initial values may change in different
intervals. For all $t\in [n\tau, (n+1)\tau[$,   using the fact that,

\begin{equation}
     \label{Impuleq60}
     \left\{
     \begin{array}{lcl}
      G(t_{n}^{+})=G(t_{n})-\lambda_{fG}G(t_{n}),\\
           \\
      T(t_{n}^{+})= T(t_{n})-\lambda_{fT}\omega(\lambda_{fG}G(t_{n}))T(t_{n}),\\
     \end{array}
     \right.
\end{equation}

we have the following discrete system for $G(n\tau)$ and $T(n\tau)$

    \begin{equation}
    \label{Impuleq61}
    \left\{
    \begin{array}{lcl}
     G((n+1)\tau)&=& \displaystyle\frac{(1-\lambda_{fG})G(n\tau)\chi((n+1)\tau,n\tau,T(n\tau))}{1+\mu_{G}G(n\tau)\int_{n\tau}^{(n+1)\tau}\chi(u,n\tau,T(n\tau))du},\\
     \\
    T((n+1)\tau)&=& \displaystyle\frac{(1-\lambda_{fT}\omega(\lambda_{fG}G(n\tau)))T(n\tau)e^{r_{T}\tau}}{1+\displaystyle\frac{T(n\tau)}{Y_{T}}(e^{r_{T}\tau}-1)}.\\
    \end{array}
    \right.
    \end{equation}

Setting

\begin{equation}
\label{Impuleq62} \left\{
\begin{array}{lcl}
 U(G(n\tau),T(n\tau))&=&\displaystyle\frac{(1-\lambda_{fG})G(n\tau)\chi((n+1)\tau,n\tau,T(n\tau))}{1+\mu_{G}G(n\tau)\int_{n\tau}^{(n+1)\tau}\chi(u,n\tau,T(n\tau))du},\\
\\
V(G(n\tau),T(n\tau))&=& \displaystyle\frac{(1-\lambda_{fT}\omega(G(n\tau)))T(n\tau)e^{r_{T}\tau}}{1+\displaystyle\frac{T(n\tau)}{Y_{T}}(e^{r_{T}\tau}-1)},\\
\end{array}
\right.
\end{equation}
implies that system $(\ref{Impuleq61})$ is equivalent to
\begin{equation}
\label{Impuleq63} \left\{
\begin{array}{lcl}
G((n+1)\tau)&=& U(G(n\tau),T(n\tau)),\\
\\
T((n+1)\tau)&=& V(G(n\tau),T(n\tau)).\\
\end{array}
\right.
\end{equation}

The existence of a periodic solution of $(\ref{Impuleq1})$ (with
period $\tau$) is equivalent to an existence of the equilibrium of
the discrete system  $(\ref{Impuleq63})$. This leads to solve the
following system

  \begin{equation}
  \label{Impuleq64}
  \left\{
  \begin{array}{lcl}
  U(x,y)&=& x,\\
  \\
  V(x,y)&=& y.\\
  \end{array}
  \right.
  \end{equation}
We have
\begin{displaymath}
U(x,y)=x\Leftrightarrow
1+\mu_{G}x\int_{n\tau}^{(n+1)\tau}\chi(u,n\tau,y)du=(1-\lambda_{fG})\chi((n+1)\tau,n\tau,y),
\end{displaymath}
which implies that

\begin{equation}
\label{Impuleq65}
x=\displaystyle\dfrac{(1-\lambda_{fG})\chi((n+1)\tau,n\tau,y)-1}{\mu_{G}\int_{n\tau}^{(n+1)\tau}\chi(u,n\tau,y)du}:=\phi(y).
\end{equation}

On the order hand, we also have

\begin{displaymath}
V(x,y)=y\Leftrightarrow
y=\displaystyle\frac{Y_{T}\{(1-\lambda_{fT}\omega(\lambda_{fG}x))e^{r_{T}\tau}-1\}}{(e^{r_{T}\tau}-1)},
\end{displaymath}
which implies that

\begin{equation}
\label{Impuleq66} \displaystyle\frac{y}{Y_{T}}(e^{r_{T}\tau}-1)
+\lambda_{fT}\omega(\lambda_{fG}x)e^{r_{T}\tau}=(e^{r_{T}\tau}-1).
\end{equation}

From (\ref{Impuleq66}) and the nonnegativity of the variable $y$ and
the function $\omega(\lambda_{fG}x)$ with $x>0$, it follows that $y$
must belong to the interval $D=[0, Y_{T}]$. Substituting $x=\phi(y)$
into the left side of (\ref{Impuleq66}) yields an equation for $y$,

$$\displaystyle\frac{y}{Y_{T}}(e^{r_{T}\tau}-1)
+\lambda_{fT}\omega(\lambda_{fG}\phi(y))e^{r_{T}\tau}-(e^{r_{T}\tau}-1)=0.$$

Set
\begin{equation}
\label{Impuleq67} h(y)=\displaystyle\frac{y}{Y_{T}}(e^{r_{T}\tau}-1)
+\lambda_{fT}\omega(\lambda_{fG}\phi(y))e^{r_{T}\tau}-(e^{r_{T}\tau}-1).
\end{equation}

It is obvious that  $h(y)$ is nonnegative and continuously
differentiable with respect to $y$. The algebraic calculation shows
that
\begin{equation}
\label{Impuleq68} \left\{
\begin{array}{lcl}
h(0)&=& \lambda_{fT}\omega(\lambda_{fG}\phi(0))e^{r_{T}\tau}-(e^{r_{T}\tau}-1),\\
\\
h(Y_{T})&=&
\lambda_{fT}\omega(\lambda_{fG}\phi(Y_{T}))e^{r_{T}\tau}.
\end{array}
\right.
\end{equation}

We have
\begin{displaymath}
\phi(0)=\displaystyle\dfrac{(1-\lambda_{fG})\chi((n+1)\tau,n\tau,0)-1}{\mu_{G}\int_{n\tau}^{(n+1)\tau}\chi(u,n\tau,0)du},
\end{displaymath}
where
\begin{displaymath}
\chi((n+1)\tau,n\tau,0)=e^{r_{G}\tau},\hspace{0.25cm}\mbox{and}
\end{displaymath}
\begin{displaymath}
\int_{n\tau}^{(n+1)\tau}\chi(u,n\tau,0)du=\int_{n\tau}^{(n+1)\tau}e^{r_{G}(u-n\tau)}du=\displaystyle\frac{1}{r_{G}}(e^{r_{G}\tau}-1).
\end{displaymath}
Then
\begin{displaymath}
\phi(0)=X_{G}\displaystyle\frac{(1-\lambda_{fG})e^{r_{G}\tau}-1}{e^{r_{G}\tau}-1}=\tilde{G}_{e}(\tau).
\end{displaymath}

Note that
\begin{equation}
\label{Impuleq69} \tilde{G}_{e}(\tau)>0\Leftrightarrow
\mathcal{R}_{0,pulse}^{\tilde{G}_{e}}>1.
\end{equation}

From (\ref{Impuleq68}), it  is easy to see that
\begin{equation}
\label{Impuleq70} h(Y_{T})>0.
\end{equation}
We also have
\begin{align*}
h(0)<0 &\Leftrightarrow \lambda_{fT}\omega(\lambda_{fG}\phi(0))e^{r_{T}\tau}- e^{r_{T}\tau}-1<0\\
&\Leftrightarrow 1<e^{r_{T}\tau}(1-\lambda_{fT}\omega(\lambda_{fG}\tilde{G}_{e}(\tau)))\\
&\Leftrightarrow\displaystyle\frac{1}{1-\lambda_{fT}\omega(\lambda_{fG}\tilde{G}_{e}(\tau))}<e^{r_{T}\tau}\\
&\Leftrightarrow \ln\left(\displaystyle\frac{1}{1-\lambda_{fT}\omega(\lambda_{fG}\tilde{G}_{e}(\tau))}\right)<r_{T}\tau\\
&\Leftrightarrow 1<\displaystyle\frac{r_{T}}{\displaystyle\frac{1}{\tau}\ln\left(\displaystyle\frac{1}{1-\lambda_{fT}\omega(\lambda_{fG}\tilde{G}_{e}(\tau))}\right)}\\
&\Leftrightarrow 1< \mathcal{R}_{0,pulse}^{*}.
\end{align*}

Then

\begin{equation}
\label{Impuleq71} h(0)<0\Leftrightarrow \mathcal{R}_{0,pulse}^{*}>1.
\end{equation}

Thus, when (\ref{Impuleq70}), and
$(\ref{Impuleq71})$ are verified, there exists at least one positive zero of $h(y)=0$ in the
interval $D$. To have the uniqueness, we show that $h$ is a highly monotone function.\par The derivative of $h(y)$ with respect to $y$ is
\begin{displaymath}
\displaystyle\frac{dh(y)}{dy}=\displaystyle\frac{1}{Y_{T}}(e^{r_{T}\tau}-1)+\phi^{'}(y)\lambda_{fT}\lambda_{fG}\omega^{'}(\lambda_{fG}\phi(y))e^{r_{T}\tau},
\end{displaymath}

where

\begin{displaymath}
\phi^{'}(y)=\phi_{1}^{'}(y)-\phi_{2}^{'}(y),
\end{displaymath}
with
\begin{displaymath}
\phi_{1}^{'}(y)=\displaystyle\frac{(1-\lambda_{fG})\chi^{'}((n+1)\tau,n\tau,y)\mu_{G}\int_{n\tau}^{(n+1)\tau}\chi(u,n\tau,y)du}{(\mu_{G}\int_{n\tau}^{(n+1)\tau}\chi(u,n\tau,y)du)^{2}},
\end{displaymath}
and
\begin{displaymath}
\phi_{2}^{'}(y)=\displaystyle\frac{((1-\lambda_{fG})\chi((n+1)\tau,n\tau,y)-1)\mu_{G}\int_{n\tau}^{(n+1)\tau}\chi^{'}(u,n\tau,y)du}{(\mu_{G}\int_{n\tau}^{(n+1)\tau}\chi(u,n\tau,y)du)^{2}}.
\end{displaymath}
We have
\begin{align*}
\chi^{'}(u,n\tau,y)&=-\displaystyle\frac{\gamma_{TG}}{\mu_{T}}\displaystyle\frac{\displaystyle\frac{(e^{r_{T}(u-n\tau)}-1)}{Y_{T}}}{1+\displaystyle\frac{y}{Y_{T}}(e^{r_{T}(u-n\tau)}-1)}e^{-m(u,n\tau,y)}\\
&=-\displaystyle\frac{\gamma_{TG}}{r_{T}}\displaystyle\frac{(e^{r_{T}(u-n\tau)}-1)}{1+\displaystyle\frac{y}{Y_{T}}(e^{r_{T}(u-n\tau)}-1)}\chi(u,n\tau,y).
\end{align*}

Substituting $\chi^{'}(u,n\tau,y)$ in $\phi_{1}^{'}(y)$ and
$\phi_{2}^{'}(y)$, we obtain

\begin{align*}
\phi_{1}^{'}(y)&=\displaystyle\frac{(1-\lambda_{fG})\left\{-\displaystyle\frac{\gamma_{TG}}{r_{T}}\displaystyle\frac{(e^{r_{T}\tau}-1)}{1+\displaystyle\frac{y}{Y_{T}}(e^{r_{T}\tau}-1)}\chi((n+1)\tau,n\tau,y)\right\}}{\mu_{G}\int_{n\tau}^{(n+1)\tau}\chi(u,n\tau,y)du}\\
&=-(1-\lambda_{fG})\chi((n+1)\tau,n\tau,y)\varTheta,
\end{align*}
where
\begin{displaymath}
\varTheta=\left(\displaystyle\frac{\displaystyle\frac{\gamma_{TG}}{r_{T}}\displaystyle\frac{(e^{r_{T}\tau}-1)}{1+\displaystyle\frac{y}{Y_{T}}(e^{r_{T}\tau}-1)}}{\mu_{G}\int_{n\tau}^{(n+1)\tau}\chi(u,n\tau,y)du}\right),
\end{displaymath}
and we also have
\begin{align*}
\phi_{2}^{'}(y)&=\displaystyle\frac{-\{(1-\lambda_{fG})\chi((n+1)\tau,n\tau,y)-1\}\displaystyle\frac{\gamma_{TG}}{r_{T}}\int_{n\tau}^{(n+1)\tau}\displaystyle\frac{(e^{r_{T}(u-n\tau)}-1)}{1+\displaystyle\frac{y}{Y_{T}}(e^{r_{T}(u-n\tau)}-1)}\chi(u,n\tau,y)du}{(\mu_{G}\int_{n\tau}^{(n+1)\tau}\chi(u,n\tau,y)du)(\int_{n\tau}^{(n+1)\tau}\chi(u,n\tau,y)du)}\\
&=\varTheta\displaystyle\frac{-\{(1-\lambda_{fG})\chi((n+1)\tau,n\tau,y)-1\}\int_{n\tau}^{(n+1)\tau}\displaystyle\frac{(e^{r_{T}(u-n\tau)}-1)\{1+\displaystyle\frac{y}{Y_{T}}(e^{r_{T}\tau}-1)\}}{(e^{r_{T}\tau}-1)\{1+\displaystyle\frac{y}{Y_{T}}(e^{r_{T}(u-n\tau)}-1)\}}\chi(u,n\tau,y)du}{\int_{n\tau}^{(n+1)\tau}\chi(u,n\tau,y)du}.
\end{align*}

Setting
\begin{displaymath}
\pi(u,n\tau,y)=\displaystyle\frac{(e^{r_{T}(u-n\tau)}-1)\{1+\displaystyle\frac{y}{Y_{T}}(e^{r_{T}\tau}-1)\}}{(e^{r_{T}\tau}-1)\{1+\displaystyle\frac{y}{Y_{T}}(e^{r_{T}(u-n\tau)}-1)\}},
\end{displaymath}

we obtain

\begin{align*}
\phi^{'}(y)&=\phi_{1}^{'}(y)-\phi_{2}^{'}(y)\\
&=\varTheta Q(y),
\end{align*}
where
\begin{displaymath}
 Q(y)=\displaystyle\frac{\{(1-\lambda_{fG})\chi((n+1)\tau,n\tau,y)-1\}\int_{n\tau}^{(n+1)\tau}\pi(u,n\tau,y)\chi(u,n\tau,y)du}{\int_{n\tau}^{(n+1)\tau}\chi(u,n\tau,y)du}-(1-\lambda_{fG})\chi((n+1)\tau,n\tau,y).
\end{displaymath}

From the expression of $\displaystyle\frac{dh(y)}{dy}$, we see that
that $h(y)$ is monotonous for all $\tau$. Hence, there is only one
positive root of $h(y)=0$ for  $\tau>0$. Therefore, when
$\mathcal{R}_{0,pulse}^{\tilde{G}_{e}}>1$ and
$\mathcal{R}_{0,pulse}^{*}>1$,  equation $h(y)=0$ has a unique
positive zero for $\tau>0$. Thus system (\ref{Impuleq1}) admits a
unique non trivial periodic solution.

\section*{Appendix D: Proof of theorem \ref{thm12} (Local stability of constant
equilibria).}
\label{AppendixD}
The proof of the stability is on the basis of the linearization to
(\ref{Impuleq63}). Letting $(G^{e},T^{e})$ be the equilibrium of
(\ref{Impuleq63}), $(0; 0)$ or $\left(0; Y_{T}\right)$. Set
$X(t)=G(t)-G^{e}$, and $Y(t)=T(t)-T^{e}$, then the linearized system
of (\ref{Impuleq63}) is

\begin{equation}
\label{Impuleq72} \left\{
\begin{array}{lcl}
X((n+1)\tau)&=& a_{11}X(n\tau)+a_{12}Y(n\tau),\\
\\
Y((n+1)\tau)&=& a_{21}X(n\tau)+a_{22}Y(n\tau),
\end{array}
\right.
\end{equation}

where,

\begin{displaymath}
a_{11}=\displaystyle\frac{\partial U}{\partial
X}(G^{e},T^{e}),\hspace{0.5cm}a_{12}=\displaystyle\frac{\partial
U}{\partial Y}(G^{e},T^{e}),
\end{displaymath}
\begin{displaymath}
a_{21}=\displaystyle\frac{\partial V}{\partial
X}(G^{e},T^{e}),\hspace{0.5cm}a_{22}=\displaystyle\frac{\partial
V}{\partial Y}(G^{e},T^{e}),
\end{displaymath}

with

\begin{equation}
\label{Impuleq73} \left\{
\begin{array}{lcl}
 U(x,y)&=&\displaystyle\frac{(1-\lambda_{fG})x\chi((n+1)\tau,n\tau,y)}{1+\mu_{G}x\int_{n\tau}^{(n+1)\tau}\chi(u,n\tau,y)du},\\
\\
V(x,y)&=& \displaystyle\frac{(1-\lambda_{fT}\omega(\lambda_{fG}x))ye^{r_{T}\tau}}{1+\displaystyle\frac{y}{Y_{T}}(e^{r_{T}\tau}-1)}.\\
\end{array}
\right.
\end{equation}
We have
\begin{displaymath}
 a_{11}=\displaystyle\frac{(1-\lambda_{fG})\chi((n+1)\tau,n\tau,T^{e})}{(1+\mu_{G}G^{e}\int_{n\tau}^{(n+1)\tau}\chi(u,n\tau,T^{e})du)^{2}}>0,
\end{displaymath}
\begin{displaymath}
a_{21}=\displaystyle\frac{-\lambda_{fT}\lambda_{fG}T^{e}\omega^{'}(\lambda_{fG}G^{e})e^{r_{T}\tau}}{\left(1+\displaystyle\frac{T^{e}}{Y_{T}}(e^{r_{T}\tau}-1)\right)^{2}}<0.
\end{displaymath}

Setting
$\xi=\left(\displaystyle\frac{-\gamma_{TG}(e^{r_{T}\tau}-1)}{r_{T}\left(1+\displaystyle\frac{T^{e}}{Y_{T}}(e^{r_{T}\tau}-1)\right)}\right)<0$,
we obtain
\begin{displaymath}
a_{12}=\displaystyle\frac{(1-\lambda_{fG})G^{e}\chi((n+1)\tau,n\tau,T^{e})\xi\left\{1+\mu_{G}G^{e}\int_{n\tau}^{(n+1)\tau}\left(1-\pi(u,n\tau,T^{e})\right)\chi(u,n\tau,T^{e})du\right\}}{(1+\mu_{G}G^{e}\int_{n\tau}^{(n+1)\tau}\chi(u,n\tau,T^{e})du)^{2}},
\end{displaymath}
where,

\begin{displaymath}
\pi(u,n\tau,T^{e})=\displaystyle\frac{(e^{r_{T}(u-n\tau)}-1)\{1+\displaystyle\frac{T^{e}}{Y_{T}}(e^{r_{T}\tau}-1)\}}{(e^{r_{T}\tau}-1)\{1+\displaystyle\frac{T^{e}}{Y_{T}}(e^{r_{T}(u-n\tau)}-1)\}}.
\end{displaymath}
It is easy to show that for $n\tau\leq u<(n+1)\tau$, we have
$\pi(u,n\tau,T^{e})<1$. Then, $a_{12}<0$.\par

We have,
\begin{displaymath}
a_{22}=\displaystyle\frac{(1-\lambda_{fT}\omega(\lambda_{fG}G^{e}))e^{r_{T}\tau}}{\left(1+\displaystyle\frac{T^{e}}{Y_{T}}(e^{r_{T}\tau}-1)\right)^{2}}>0.
\end{displaymath}
The stability of the equilibrium of (\ref{Impuleq63}) can be
determined by eigenvalues of the linearized matrix
$A=(a_{ij})_{2\times 2}$.

\begin{enumerate}
\item The two eigenvalues for $A$ at $(0; 0)$ are
\begin{align*}
\rho_{1}&=(1-\lambda_{fG})\chi((n+1)\tau,n\tau,0)\\
&=(1-\lambda_{fG})e^{r_{G}\tau},
\end{align*}
and
\begin{displaymath}
\rho_{2}=e^{r_{T}\tau}>1.
\end{displaymath}
Then, the trivial equilibrium $(0, 0)$ is always unstable.
\item Concerning the local stability of forest equilibrium, we have the following two eigenvalues
\begin{align*}
\nu_{1}&=(1-\lambda_{fG})\chi((n+1)\tau,n\tau,Y_{T})\\
&=(1-\lambda_{fG})e^{\left(r_{G}\tau-\displaystyle\frac{\gamma_{TG}}{\mu_{T}}r_{T}\tau\right)}\\
&=(1-\lambda_{fG})e^{r_{G}\tau\left(1-\displaystyle\frac{1}{\mathcal{R}_{01}}\right)},
\end{align*}
where
\begin{displaymath}
\mathcal{R}_{01}=\displaystyle\frac{r_{G}}{r_{T}}\times\displaystyle\frac{\mu_{T}}{\gamma_{TG}}.
\end{displaymath}

The second eigenvalue is
\begin{displaymath}
\nu_{2}=\displaystyle\frac{1}{e^{r_{T}\tau}}<1.
\end{displaymath}

\end{enumerate}
Thus, we have the following results:
\begin{itemize}
\item If $\mathcal{R}_{01}\leq 1$, then $\nu_{1}<1$, which implies that the forest equilibrium $E_{01}=(0; Y_{T})$ is locally asymptotically stable (LAS) (similarly in the continuous model (see
Tchuinte et al. 2014 \cite{Tchuinte2014})).
\item If $\mathcal{R}_{01}>1$ and  $\tilde{\mathcal{R}}_{0,\mathcal{R}_{01}}=\mathcal{R}_{0,pulse}^{\tilde{G}_{e}}\left(1-\displaystyle\frac{1}{\mathcal{R}_{01}}\right)<1$, then $(0; Y_{T})$ is LAS since
$\nu_{1}<1$. This situation is specific for the impulse model. The
continuous model does not imply the stability of the forest
equilibrium when $\mathcal{R}_{01}>1$.
\item If  $\mathcal{R}_{01}>1$ and $\tilde{\mathcal{R}}_{0,\mathcal{R}_{01}}>1$, then $(0; Y_{T})$ is unstable.
\end{itemize}

\section*{Appendix E: Proof of theorem \ref{thm13} (Local stability of the periodic grassland equilibrium)}
\label{AppendixE}
To show that $(\tilde{G}_{e}(t), 0)$ is LAS, we consider Floquet's
theory.  Set $G(t)=\tilde{G}_{e}(t)+x(t)$, and $T(t)=0+y(t)$, where
$x(t)$, and $y(t)$ are small perturbations, every solution of the
linearized equations can be written as

\begin{displaymath}
 \left(
 \begin{array}{ccccc}
 x(t)\\
 \\
 y(t)\\
 \end{array}
 \right)
 =\Phi(t)\left(
 \begin{array}{ccccc}
 x(0)\\
 \\
 y(0)\\
 \end{array}
 \right),
 \end{displaymath}

 where, $\Phi(t)=\left(\varphi_{ij}(t)\right), i,j=1,2$ is a fundamental matrix and satisfies

     \begin{equation}
        \label{Impuleq74}
        \begin{array}{lcl}
         \displaystyle\frac{d\Phi(t)}{dt}&=& \mathcal{A}(t)\Phi(t),\\
        \end{array}
        \end{equation}
 with,

  \begin{displaymath}
 \mathcal{ A}(t)=\left(
  \begin{array}{ccccc}
  r_{G}-2\mu_{G}\tilde{G}_{e}(t)&-\gamma_{TG}\tilde{G}_{e}(t)\\
  & \\
  0& r_{T}\\
  \end{array}
  \right).
  \end{displaymath}

Since $\Phi(t)$ is the principal fundamental matrix, then
$\Phi(0)=I_{2},$ where $I_{2}$  is the identity matrix of
$\mathcal{M}_{2}(\mathbf{R}).$  Then,

 \begin{displaymath}
 \left(
 \begin{array}{ccccc}
 x(n\tau^{+})\\
 \\
 y(n\tau^{+})\\
 \end{array}
 \right)
 =\left(
 \begin{array}{ccccc}
 1-\lambda_{fG} & 0\\
 & \\
 0 & 1-\lambda_{fT}\omega(\lambda_{fG}\tilde{G}_{e}(n\tau))\\
 \end{array}
 \right)
 \left(
 \begin{array}{ccccc}
 x(n\tau)\\
 \\
 y(n\tau)\\
 \end{array}
 \right),
 \end{displaymath}

and hence, if the absolute value of all eigenvalues (Floquet
multipliers) of the monodromy matrix

\begin{displaymath}
M=\left(
\begin{array}{ccccc}
1-\lambda_{fG} & 0\\
&\\
0 & 1-\lambda_{fT}\omega(\lambda_{fG}\tilde{G}_{e}(\tau))\\
\end{array}
\right)\Phi(\tau)
\end{displaymath}

are less than one, the periodic solution $(\tilde{G}_{e}(t);0)$ is
locally asymptotically stable.\par By calculation, we obtain

\begin{displaymath}
\Phi(t)=\left(
\begin{array}{ccccc}
\varphi_{11}(t)&\varphi_{12}(t)\\
&\\
0&\varphi_{22}(t)\\
\end{array}
\right),
\end{displaymath}

where,

\begin{displaymath}
\varphi_{11}(t)=\exp\left\{r_{G}t-2\mu_{G}\int_{0}^{t}\tilde{G}_{e}(u)du\right\},\\
\end{displaymath}
and
 \begin{displaymath}
 \varphi_{22}(t)=\exp\left\{\int_{0}^{t}r_{T}du\right\}.
 \end{displaymath}

 We deduce that, eigenvalues $\lambda_{i}, i=1, 2$ of $M$ are

 \begin{equation}
  \label{Impuleq75}
  \left\{
  \begin{array}{lcl}
\lambda_{1}=(1-\lambda_{fG})\varphi_{11}(\tau),\\
 \\
 \lambda_{2}=(1-\lambda_{fT}\omega(\lambda_{fG}\tilde{G}_{e}(\tau)))\varphi_{22}(\tau).\\
   \end{array}
 \right.
 \end{equation}

From the first equation of (\ref{Impuleq1}), we have

\begin{equation}
 \label{Impuleq76}
 \displaystyle\frac{d\tilde{G}}{\tilde{G}}=(r_{G}-\mu_{G}\tilde{G}-\gamma_{TG}\tilde{T})dt.
\end{equation}

Integrating (\ref{Impuleq76}) in $[0, \tau]$, we obtain

\begin{equation}
 \label{Impuleq77}
 \tilde{G}(\tau)=\tilde{G}(0)\exp\left\{r_{G}\tau-\mu_{G}\int_{0}^{\tau}\tilde{G}(u)du-\gamma_{TG}\int_{0}^{\tau}\tilde{T}(u)du\right\}.
\end{equation}
At $(\tilde{G}_{e}(t); 0)$, equation (\ref{Impuleq77}) becomes

\begin{equation}
 \label{Impuleq78}
 \tilde{G}_{e}(\tau)=\tilde{G}_{e}(0)\exp\left\{r_{G}\tau-\mu_{G}\int_{0}^{\tau}\tilde{G}(u)du\right\}.
\end{equation}

From the third equation of (\ref{Impuleq1}), we have

\begin{equation}
\label{Impuleq79} \tilde{G}_{e}(0)=(1-\lambda_{fG})
\tilde{G}_{e}(\tau).
\end{equation}

Substituting (\ref{Impuleq79}) in (\ref{Impuleq78}), we obtain
\begin{displaymath}
 \tilde{G}_{e}(\tau)=(1-\lambda_{fG}) \tilde{G}_{e}(\tau)\exp\left\{r_{G}\tau-\mu_{G}\int_{0}^{\tau}\tilde{G}(u)du\right\},
\end{displaymath}
 which implies that
 \begin{displaymath}
 (1-\lambda_{fG}) \exp\left\{r_{G}\tau-\mu_{G}\int_{0}^{\tau}\tilde{G}(u)du\right\}=1.
 \end{displaymath}

 Thus,
\begin{align*}
\lambda_{1}&=(1-\lambda_{fG})\varphi_{11}(\tau)\\
           &=(1-\lambda_{fG})\exp\left\{r_{G}\tau-2\mu_{G}\int_{0}^{\tau}\tilde{G}(u)du\right\}\\
           &=(1-\lambda_{fG})\exp\left\{r_{G}\tau-\mu_{G}\int_{0}^{\tau}\tilde{G}(u)du\right\}\exp\left\{-\mu_{G}\int_{0}^{\tau}\tilde{G}(u)du\right\}\\
           &=\exp\left\{-\mu_{G}\int_{0}^{\tau}\tilde{G}(u)du\right\}<1.
\end{align*}
On the other hand, we have
 \begin{displaymath}
 \lambda_{2}<1\Leftrightarrow \mathcal{R}^{*}_{0,pulse}=\displaystyle\frac{r_{T}}{\displaystyle\frac{1}{\tau}\ln\left(\displaystyle\frac{1}{1-\lambda_{fT}\omega(\lambda_{fG}\tilde{G}_{e}(\tau))}\right)}<1.
 \end{displaymath}
Then, the periodic grassland equilibrium $\tilde{E}_{10}$ is LAS if
$\mathcal{R}^{*}_{0,pulse}<1$. This completes the proof.

\section*{Appendix F: Proof of theorem \ref{thm14} (local stability of the periodic savanna equilibrium)}
\label{AppendixF}
Now, we investigate local properties of the periodic savanna
equilibrium. Similarly to the proof of the local stability of the
periodic grassland equilibrium, we set $G(t)=\tilde{G}^{*}(t)+x(t)$,
and $T(t)=\tilde{T}^{*}(t)+y(t)$, where $x(t)$, and $y(t)$ are small
perturbations and they are solutions of the linearized equations

\begin{displaymath}
 \left(
 \begin{array}{ccccc}
 x(t)\\
 \\
 y(t)\\
 \end{array}
 \right)
 =\Phi^{*}(t)\left(
 \begin{array}{ccccc}
 x(0)\\
 \\
 y(0)\\
 \end{array}
 \right).
 \end{displaymath}

 $\Phi^{*}(t)$ is the fundamental principal matrix which satisfies

 \begin{displaymath}
\displaystyle\frac{d\Phi^{*}(t)}{dt}=\left(
 \begin{array}{ccccc}
 r_{G}-2\mu_{G}\tilde{G}^{*}(t)-\gamma_{TG}\tilde{T}^{*}(t)&-\gamma_{TG}\tilde{G}^{*}(t)\\
 & \\
 0& r_{T}-2\mu_{T}\tilde{T}^{*}(t)\\
 \end{array}
 \right)\Phi^{*}(t).
 \end{displaymath}

 By calculation, we have

 \begin{displaymath}
 \Phi^{*}(t)=\left(
 \begin{array}{ccccc}
\varphi^{*}_{11}(t)&\varphi^{*}_{12}(t)\\
 &\\
 0&\varphi^{*}_{22}(t)\\
 \end{array}
 \right),
 \end{displaymath}

 where,

 \begin{displaymath}
 \varphi^{*}_{11}(t)=\exp\left\{\int_{0}^{t}(r_{G}-2\mu_{G}\tilde{G}^{*}(u)-\gamma_{TG}\tilde{T}^{*}(u))du\right\},\\
 \end{displaymath}
 and
 \begin{displaymath}
 \varphi^{*}_{22}(t)=\exp\left\{\int_{0}^{t}(r_{T}-2\mu_{T}\tilde{T}^{*}(u))du\right\}.
 \end{displaymath}

From system (\ref{Impuleq1}), at $t=t_{n}$, we have the following
system

\begin{displaymath}
 \left(
 \begin{array}{ccccc}
 x(n\tau^{+})\\
 \\
 y(n\tau^{+})\\
 \end{array}
 \right)
 =\left(
 \begin{array}{ccccc}
 1-\lambda_{fG} & 0\\
 & \\
 0 & 1-\lambda_{fT}\omega(\lambda_{fG}\tilde{G}^{*}(n\tau))\\
 \end{array}
 \right)
 \left(
 \begin{array}{ccccc}
 x(n\tau)\\
 \\
 y(n\tau)\\
 \end{array}
 \right).
 \end{displaymath}

Hence, according to the Floquet theory, if all eigenvalues (Floquet
multipliers) $\lambda_{1}^{*}$ and $\lambda_{2}^{*}$ of

\begin{displaymath}
M^{*}=\left(
\begin{array}{ccccc}
1-\lambda_{fG} & 0\\
&\\
0 & 1-\lambda_{fT}\omega(\lambda_{fG}\tilde{G}^{^{*}}(\tau))\\
\end{array}
\right)\Phi^{*}(\tau)
\end{displaymath}
are less than one, then the coexistence Tree-Grass periodic
equilibrium is locally asymptotically stable. We have
\begin{equation}
 \label{Impuleq80}
 \left\{
 \begin{array}{lcl}
\lambda^{*}_{1}&=(1-\lambda_{fG})\varphi^{*}_{11}(\tau)= (1-\lambda_{fG})\exp\left\{r_{G}\tau-2\mu_{G}\int_{0}^{\tau}\tilde{G}^{*}(u)du-\gamma_{TG}\int_{0}^{\tau}\tilde{T}^{*}(u)du\right\},\\
\\
\lambda^{*}_{2}&=\{1-\lambda_{fT}\omega(\lambda_{fG}\tilde{G}^{*}(\tau))\}\varphi^{*}_{22}(\tau)=(1-\lambda_{fT}\omega(\lambda_{fG}\tilde{G}^{*}(\tau)))\exp\left\{r_{T}\tau-2\mu_{T}\int_{0}^{\tau}\tilde{T}^{*}(s)ds\right\}.
 \end{array}
\right.
\end{equation}
Starting with $\lambda^{*}_{1}$, we have
\begin{align*}
\lambda^{*}_{1}<1&\Leftrightarrow(1-\lambda_{fG})\exp\left\{r_{G}\tau-2\mu_{G}\int_{0}^{\tau}\tilde{G}^{*}(u)du-\gamma_{TG}\int_{0}^{\tau}\tilde{T}^{*}(u)du\right\}<1\\
 &\Leftrightarrow 1-\displaystyle\frac{2\mu_{G}}{r_{G}}\left(\displaystyle\frac{1}{\tau}\int_{0}^{\tau}\tilde{G}^{*}(u)du\right)-\displaystyle\frac{\gamma_{TG}}{r_{G}}\left(\displaystyle\frac{1}{\tau}\int_{0}^{\tau}\tilde{T}^{*}(u)du\right)<\displaystyle\frac{\displaystyle\frac{1}{\tau}\ln\left(\displaystyle\frac{1}{1-\lambda_{fG}}\right)}{r_{G}}\\
  &\Leftrightarrow 1-\displaystyle\frac{2}{X_{G}}\left(\displaystyle\frac{1}{\tau}\int_{0}^{\tau}\tilde{G}^{*}(u)du\right)-\displaystyle\frac{\gamma_{TG}}{r_{G}}\left(\displaystyle\frac{1}{\tau}\int_{0}^{\tau}\tilde{T}^{*}(u)du\right)<\displaystyle\frac{1}{\mathcal{R}_{0,pulse}^{\tilde{G}_{e}}}.\\
\end{align*}

Integrating $\tilde{G}^{*}(t)$ and $\tilde{T}^{*}(t)$ from $0$ to
$\tau$, we obtain

\begin{equation}
\label{Impuleq81} \left\{
\begin{array}{lcl}
\int_{0}^{\tau}\tilde{G}^{*}(u)du&=& \displaystyle\frac{1}{\mu_{G}}\ln\left(1+\mu_{G}G^{*}\int_{0}^{\tau}\chi(u,0,T^{*})du\right),\\
\\
\int_{0}^{\tau}\tilde{T}^{*}(u)du&=& \displaystyle\frac{1}{\mu_{T}}\ln\left(1+\displaystyle\frac{T^{*}}{Y_{T}}(e^{r_{T}\tau}-1)\right).\\
\end{array}
\right.
\end{equation}
Substituting (\ref{Impuleq44c}) which is expression of $T^{*}$ in
the second equation of (\ref{Impuleq81}), we have

\begin{equation}
\label{Impuleq82}
\int_{0}^{\tau}\tilde{T}(u)du=\displaystyle\frac{1}{\mu_{T}}\left[\ln(1-\lambda_{fT}\omega(\lambda_{fG}G^{*}))+r_{T}\tau\right].
\end{equation}

Then,
\begin{align*}
\displaystyle\frac{\gamma_{TG}}{r_{G}}\left(\displaystyle\frac{1}{\tau}\int_{0}^{\tau}\tilde{T}^{*}(u)du\right)&=\displaystyle\frac{\gamma_{TG}}{r_{G}}\displaystyle\frac{1}{\mu_{T}}\left\{r_{T}-\displaystyle\frac{1}{\tau}\ln\displaystyle\frac{1}{(1-\lambda_{fT}\omega(\lambda_{fG}G^{*}))}\right\}\\
&=\displaystyle\frac{\gamma_{TG}}{r_{G}}\displaystyle\frac{r_{T}}{\mu_{T}}\left\{1-\displaystyle\frac{\displaystyle\frac{1}{\tau}\ln\displaystyle\frac{1}{(1-\lambda_{fT}\omega(\lambda_{fG}G^{*}))}}{r_{T}}\right\}\\
&=\displaystyle\frac{1}{\mathcal{R}_{01}}\left(1-\displaystyle\frac{1}{\mathcal{R}_{0,stable}^{*}}\right):=\mathcal{\tilde{R}}_{0,\mathcal{R}_{01}}^{G^{*}}.
\end{align*}

Thus we have
\begin{align*}
\lambda^{*}_{1}<1&\Leftrightarrow 1-\displaystyle\frac{2}{X_{G}}\left(\displaystyle\frac{1}{\tau}\int_{0}^{\tau}\tilde{G}^{*}(u)du\right)<\mathcal{\tilde{R}}_{0,\mathcal{R}_{01}}^{G^{*}}+\displaystyle\frac{1}{\mathcal{R}_{0,pulse}^{\tilde{G}_{e}}}\\
&\Leftrightarrow 1<\mathcal{\tilde{R}}_{0,\mathcal{R}_{01}}^{G^{*}}+\displaystyle\frac{1}{\mathcal{R}_{0,pulse}^{\tilde{G}_{e}}}+\displaystyle\frac{2}{X_{G}}\left(\displaystyle\frac{1}{\tau}\int_{0}^{\tau}\tilde{G}^{*}(u)du\right):=\mathcal{\tilde{R}}_{0,stable}^{**}.\\
\end{align*}

Therefore,
\begin{equation}
\label{Impuleq83} \lambda^{*}_{1}<1\Leftrightarrow
\mathcal{\tilde{R}}_{0,stable}^{**}>1.
\end{equation}

We have $\chi(\tau,\tau,T^{*})=1$, which implies that
\begin{align*}
\tilde{G}^{*}(\tau)&=\displaystyle\frac{\chi(\tau,\tau,T^{*})G^{*}}{1+\mu_{G}G^{*}\int_{\tau}^{\tau}\chi(u,\tau,T^{*})du}\\
&=\chi(\tau,\tau,T^{*})G^{*}\\
&=G^{*}.
\end{align*}

Then,
\begin{align*}
\lambda_{2}&=(1-\lambda_{fT}\omega(G^{*}))\exp\left\{-r_{T}\tau-2\ln(1-\lambda_{fT}\omega(\lambda_{fG}G^{*}))\right\}\\
&=\exp\left\{-r_{T}\tau-\ln(1-\lambda_{fT}\omega(\lambda_{fG}G^{*}))\right\}\\
&=\exp\left\{-r_{T}\tau+\ln\displaystyle\frac{1}{(1-\lambda_{fT}\omega(\lambda_{fG}G^{*}))}\right\}.\\
&=\exp\left\{r_{T}\tau\left(\displaystyle\frac{\displaystyle\frac{1}{\tau}\ln\displaystyle\frac{1}{(1-\lambda_{fT}\omega(\lambda_{fG}G^{*}))}}{r_{T}}-1\right)\right\}\\
&=\exp\left\{r_{T}\tau\left(\displaystyle\frac{1}{\mathcal{R}_{0,stable}^{*}}-1\right)\right\}.\\
\end{align*}

Thus,
\begin{equation}
\label{Impuleq84}
\lambda_{2}<1\Longleftrightarrow\mathcal{R}_{0,stable}^{*}>1.
\end{equation}

From (\ref{Impuleq83}) and (\ref{Impuleq84}), it follows that the
coexistence Tree-Grass periodic equilibrium $\tilde{E}^{*}_{11}$ is
locally asymptotically stable if $\mathcal{R}_{0,stable}^{*}>1$ and
$\mathcal{\tilde{R}}_{0,stable}^{**}>1$. This proof completes this
section.

\section*{Appendix G: Proof of theorem \ref{thm15} (Global stability of the forest
equilibrium)}
\label{AppendixG}
From system (\ref{Impuleq1}), we have

\begin{equation}
\label{Impuleq85} \left\{
\begin{array}{lcl}
G^{'}(t)&\leq & r_{G}G(t),\hspace{0.5cm}t\neq t_{k},\\
\\
G(t_{k}^{+})&=&(1-\lambda_{fG})G(t_{k})\hspace{0.5cm}t=t_{k}.
\end{array}
\right.
\end{equation}
Using  impulsive differential inequations (Lakshmikantham et al 1989 \cite{Lakshmikantham1989}), we have

\begin{align*}
G(t)&\leq G(t_{0}^{+})\left(\prod_{t_{0}<t_{k}<t}(1-\lambda_{fG})\right)\exp\left(\int_{t_{0}}^{t}r_{G}ds\right)\\
&=G_{0}(1-\lambda_{fG})^{\left(\left[\displaystyle\frac{t}{\tau}\right]-\left[\displaystyle\frac{t_{0}}{\tau}\right]\right)}\exp\left\{r_{G}(t-t_{0})\right\}\\
&\leq G_{0}(1-\lambda_{fG})^{\left(\left[\displaystyle\frac{t}{\tau}\right]-\left[\displaystyle\frac{t_{0}}{\tau}\right]\right)}\exp\left\{r_{G}([t-t_{0}]+1)\right\}\\
&= G_{0}e^{r_{G}\tau}(1-\lambda_{fG})^{\left(\left[\displaystyle\frac{t}{\tau}\right]-\left[\displaystyle\frac{t_{0}}{\tau}\right]\right)}\exp\left\{r_{G}([t-t_{0}])\right\}\\
&=G_{0}e^{r_{G}\tau}(1-\lambda_{fG})^{\left(\left[\displaystyle\frac{t}{\tau}\right]-\left[\displaystyle\frac{t_{0}}{\tau}\right]\right)}\exp\left\{r_{G}\tau\left(\left[\displaystyle\frac{t}{\tau}\right]-\left[\displaystyle\frac{t_{0}}{\tau}\right]\right)\right\}\\
&=G_{0}e^{r_{G}\tau}\exp\left\{r_{G}\tau\left(1-\displaystyle\frac{\displaystyle\frac{1}{\tau}\ln\left(\displaystyle\frac{1}{1-\lambda_{fG}}\right)}{r_{G}}\right)\right\}^{\left(\left[\displaystyle\frac{t}{\tau}\right]-\left[\displaystyle\frac{t_{0}}{\tau}\right]\right)}\\
&=G_{0}e^{r_{G}\tau}\exp\left\{r_{G}\tau\left(1-\displaystyle\frac{1}{\mathcal{R}_{0,pulse}^{\tilde{G}_{e}}}\right)\right\}^{\left(\left[\displaystyle\frac{t}{\tau}\right]-\left[\displaystyle\frac{t_{0}}{\tau}\right]\right)}.\\
\end{align*}

Then,  $\mathcal{R}_{0,pulse}^{\tilde{G}_{e}}<1$ implies that
$G(t)\rightarrow 0$, as $t\rightarrow\infty$.\par

Now, we will prove that
\begin{equation}
\label{Impuleq86} \lim\limits_{t\longrightarrow \infty}T(t)=Y_{T}.
\end{equation}

When $G(t)\equiv 0$ and $n\tau<t\leq(n+1)\tau,$
$n=1,2,...,N_{\tau}$, system $(\ref{Impuleq1})$ becomes

\begin{equation}
\label{Impuleq87} \left\{
\begin{array}{lcl}
\displaystyle\frac{dT}{dt}=r_{T}T-\mu_{T}T^{2},~~~~~t\neq t_{n},\\
\\
T(t_{n}^{+})= T(t_{n})-\lambda_{fT}\omega(0)T(t_{n}),~~~~~t= t_{n}.\\
\end{array}
\right.
\end{equation}
since $\omega(0)=0$, system $(\ref{Impuleq87})$ is equivalent to

\begin{equation}
\label{Impuleq88}
 \displaystyle\frac{dT}{dt}=r_{T}T-\mu_{T}T^{2},
\end{equation}

which has two equilibria $0$, and $Y_{T}$. It is obvious that
$Y_{T}$ is GAS. Thus,

\begin{displaymath}
\lim\limits_{t\rightarrow+\infty }T(t)=Y_{T}.
\end{displaymath}
This complete the proof. Thus, when
$\mathcal{R}_{0,pulse}^{\tilde{G}_{e}}<1$ the forest equilibrium
$E_{01}$ is GAS.

\section*{Appendix H: Proof of theorem \ref{thm16}}
\label{AppendixH}
We prove  the global stability of $\tilde{E}_{10}$ in the following
two steps:\par
 \begin{itemize}
 \item  Step 1.
First, we show that $\lim\limits_{t\rightarrow\infty}T(t)=0$ if
$\mathcal{R}_{0,stable}^{\tilde{G}_{e}}<1$. In fact, from system
(\ref{Impuleq1}), we obtain

\begin{equation}
\label{Impuleq91}
\left\{
\begin{array}{lcl}
T^{'}(t)&\leq &r_{T}T(t), \hspace{0.5cm}t\neq t_{n},\\
\\
T(t_{n}^{+})&=&(1-\lambda_{fT}\omega(\lambda_{fG}\tilde{G}_{e}(t_{n})))T(t_{n}),
\hspace{0.5cm}t= t_{n}.
\end{array}
\right.
\end{equation}

Using impulsive differential inequations (Lakshmikantham et al 1989 \cite{Lakshmikantham1989}) we show that

\begin{displaymath}
T(t)\leq
T(t_{0}^{+})\left(\prod_{t_{0}<t_{n}<t}(1-\lambda_{fT}\omega(\lambda_{fG}\tilde{G}_{e}(t_{n})))\right)\exp\left(\int_{t_{0}}^{t}r_{T}ds\right).
\end{displaymath}

It is obvious that  $\tilde{G}_{e}$ increases monotonically in
$[n\tau,(n+1)\tau[$, $n=1,2,...$
On the other hand $\omega$ is an
increasing function, then  for all $n=1,2,...$ we have
\begin{displaymath}
(1-\lambda_{fT}\omega(\lambda_{fG}\tilde{G}_{e}(n\tau)))\leq
(1-\lambda_{fT}\omega(\lambda_{fG}\tilde{G}_{e}(\tau))).
\end{displaymath}

Then, we have
\begin{align*}
T(t)&\leq T(t_{0}^{+})\left(\prod_{t_{0}<t_{k}<t}(1-\lambda_{fT}\omega(\lambda_{fG}\tilde{G}_{e}(\tau)))\right)\exp\left(\int_{t_{0}}^{t}r_{T}ds\right)\\
   &=T_{0}(1-\lambda_{fT}\omega(\lambda_{fG}\tilde{G}_{e}(\tau)))^{\left(\left[\displaystyle\frac{t}{\tau}\right]-\left[\displaystyle\frac{t_{0}}{\tau}\right]\right)}\exp\{r_{T}(t-t_{0})\}\\
   &\leq T_{0}(1-\lambda_{fT}\omega(\lambda_{fG}\tilde{G}_{e}(\tau)))^{\left(\left[\displaystyle\frac{t}{\tau}\right]-\left[\displaystyle\frac{t_{0}}{\tau}\right]\right)}\exp\{r_{T}([t-t_{0}]+1)\}\\
   &=T_{0}e^{r_{T}}(1-\lambda_{fT}\omega(\lambda_{fG}\tilde{G}_{e}(\tau)))^{\left(\left[\displaystyle\frac{t}{\tau}\right]-\left[\displaystyle\frac{t_{0}}{\tau}\right]\right)}\exp\left\{r_{T}\tau\left(\left[\displaystyle\frac{t}{\tau}\right]-\left[\displaystyle\frac{t_{0}}{\tau}\right]\right)\right\}\\
   &=T_{0}e^{r_{T}}\exp\left\{r_{T}\tau-\ln\left(\displaystyle\frac{1}{1-\lambda_{fT}\omega(\lambda_{fG}\tilde{G}_{e}(\tau))}\right)\right\}^{\left(\left[\displaystyle\frac{t}{\tau}\right]-\left[\displaystyle\frac{t_{0}}{\tau}\right]\right)}\\
   &=T_{0}e^{r_{T}}\exp\left\{r_{T}\tau \left(1-\displaystyle\frac{\displaystyle\frac{1}{\tau}\ln\left(\displaystyle\frac{1}{1-\lambda_{fT}\omega(\lambda_{fG}\tilde{G}_{e}(\tau))}\right)}{r_{T}}\right)\right\}^{\left(\left[\displaystyle\frac{t}{\tau}\right]-\left[\displaystyle\frac{t_{0}}{\tau}\right]\right)}\\
   &=T_{0}e^{r_{T}}\exp\left\{r_{T}\tau \left(1-\displaystyle\frac{1}{\mathcal{R}_{0,pulse}^{*}}\right)\right\}^{\left(\left[\displaystyle\frac{t}{\tau}\right]-\left[\displaystyle\frac{t_{0}}{\tau}\right]\right)}\\
\end{align*}

Thus, when $\mathcal{R}_{0,pulse}^{*}<1$, $T(t)\rightarrow 0$ as
$t\rightarrow \infty$.\par

\item Step 2.
We prove  that
$\lim\limits_{t\rightarrow\infty}|G(t)-\tilde{G}_{e}(t)|=0$.\par
Since $\lim\limits_{t\rightarrow\infty}T(t)=0$, then for any given
$\epsilon_{1}>0$, there exists $t_{1}>0$, such that

\begin{equation}
   \label{Impuleq92}
   -\epsilon_{1}\leq T(t)\leq \epsilon_{1},
\end{equation}
For all $t>t_{1}$. Using (\ref{Impuleq92}) into the first equation
of (\ref{Impuleq1}), we obtain

\begin{equation}
\label{Impuleq93} \left\{
\begin{array}{lcl}
\displaystyle\frac{dG}{dt}\geq r_{G}G-\mu_{G}G^{2}-\gamma_{TG}\epsilon_{1}G,~~~~~t\neq t_{n},\\
\\
G(t_{n}^{+})= G(t_{n})-\lambda_{fG}G(t_{n}),~~~~~t= t_{n}.\\
\end{array}
\right.
\end{equation}
Let $z=G^{-1}$, we have
$\displaystyle\frac{dz}{dt}=-\displaystyle\frac{1}{G^{2}}\displaystyle\frac{dG}{dt}$.
Then system (\ref{Impuleq93}) changes into following
\begin{equation}
\label{Impuleq94} \left\{
\begin{array}{lcl}
\displaystyle\frac{dz}{dt}\leq -(r_{G}-\gamma_{TG}\epsilon_{1})z+\mu_{G},~~~~~t\neq t_{n},\\
\\
z(t_{n}^{+})=\displaystyle\frac{1}{1-\lambda_{fG}} z(t_{n}),~~~~~t= t_{n}.\\
\end{array}
\right.
\end{equation}
Set $A_{\epsilon_{1}}=r_{G}-\gamma_{TG}\epsilon_{1}$. Using impulsive differential inequations (Lakshmikantham et al 1989 \cite{Lakshmikantham1989}), we have

\begin{align*}
z(t)&\leq z(t_{1}^{+})\prod_{t_{1}<t_{k}<t}\left(\displaystyle\frac{1}{1-\lambda_{fG}}\right)\exp\left(-\int_{t_{1}}^{t}A_{\epsilon_{1}}ds\right)\\
&+\int_{t_{1}}^{t}\prod_{s<t_{k}<t}\left(\displaystyle\frac{1}{1-\lambda_{fG}}\right)\mu_{G}\exp\left(-\int_{s}^{t}A_{\epsilon_{1}}d\sigma\right)ds\\
&=z(t_{1}^{+})\left(\displaystyle\frac{1}{1-\lambda_{fG}}\right)^{\left(\left[\displaystyle\frac{t}{\tau}\right]-\left[\displaystyle\frac{t_{1}}{\tau}\right]\right)}\exp\left(-A_{\epsilon_{1}}(t-t_{1})\right)\\
&+\mu_{G}\exp\left(-A_{\epsilon_{1}}t\right) \int_{t_{1}}^{t}\prod_{s<t_{k}<t}\left(\displaystyle\frac{1}{1-\lambda_{fG}}\right)\exp\left(A_{\epsilon_{1}}s\right)ds,\\
\end{align*}
which implies that

\begin{align*}
z(t)&\leq
z(t_{1}^{+})\left(\displaystyle\frac{1}{1-\lambda_{fG}}\right)^{\left(\left[\displaystyle\frac{t}{\tau}\right]-\left[\displaystyle\frac{t_{1}}{\tau}\right]\right)}\exp\left(-A_{\epsilon_{1}}(t-t_{1})\right)\\
&+\mu_{G}\exp\left(-A_{\epsilon_{1}}t\right)[
\int_{t_{1}}^{[t_{1}]+\tau}\prod_{s<t_{k}<t}\left(\displaystyle\frac{1}{1-\lambda_{fG}}\right)\exp\left(A_{\epsilon_{1}}s\right)ds \\
& +\int_{[t_{1}]+\tau}^{[t_{1}]+2\tau}\prod_{s<t_{k}<t}\left(\displaystyle\frac{1}{1-\lambda_{fG}}\right)\exp\left(A_{\epsilon_{1}}s\right)ds\\
&+ ...+\int_{[t]-\tau}^{[t]}\prod_{s<t_{k}<t}\left(\displaystyle\frac{1}{1-\lambda_{fG}}\right)\exp\left(A_{\epsilon_{1}}s\right)ds\\
&  +\int_{[t]}^{t}\prod_{s<t_{k}<t}\left(\displaystyle\frac{1}{1-\lambda_{fG}}\right)\exp\left(A_{\epsilon_{1}}s\right)ds]\\
&= z(t_{1}^{+})\left(\displaystyle\frac{1}{1-\lambda_{fG}}\right)^{\left(\left[\displaystyle\frac{t}{\tau}\right]-\left[\displaystyle\frac{t_{1}}{\tau}\right]\right)}\exp\left(-A_{\epsilon_{1}}(t-t_{1})\right)\\
&+\displaystyle\frac{\mu_{G}}{A_{\epsilon_{1}}}\exp\left(-A_{\epsilon_{1}}t\right)[
\int_{t_{1}}^{[t_{1}]+\tau}\prod_{s<t_{k}<t}\left(\displaystyle\frac{1}{1-\lambda_{fG}}\right)\exp\left(A_{\epsilon_{1}}s\right)d(A_{\epsilon_{1}}s) \\
& +\int_{[t_{1}]+\tau}^{[t_{1}]+2\tau}\prod_{s<t_{k}<t}\left(\displaystyle\frac{1}{1-\lambda_{fG}}\right)\exp\left(A_{\epsilon_{1}}s\right)d(A_{\epsilon_{1}}s)\\
&+ ...+\int_{[t]-\tau}^{[t]}\prod_{s<t_{k}<t}\left(\displaystyle\frac{1}{1-\lambda_{fG}}\right)\exp\left(A_{\epsilon_{1}}s\right)d(A_{\epsilon_{1}}s)\\
&  +\int_{[t]}^{t}\prod_{s<t_{k}<t}\left(\displaystyle\frac{1}{1-\lambda_{fG}}\right)\exp\left(A_{\epsilon_{1}}s\right)d(A_{\epsilon_{1}}s)]\\
&= z(t_{1}^{+})\left(\displaystyle\frac{1}{1-\lambda_{fG}}\right)^{\left(\left[\displaystyle\frac{t}{\tau}\right]-\left[\displaystyle\frac{t_{1}}{\tau}\right]\right)}e^{(-A_{\epsilon_{1}}(t-t_{1}))}\\
&+\displaystyle\frac{\mu_{G}}{A_{\epsilon_{1}}}e^{(-A_{\epsilon_{1}}t)}[
\left(\displaystyle\frac{1}{1-\lambda_{fG}}\right)^{\left[\displaystyle\frac{t-t_{1}}{\tau}\right]}(e^{A_{\epsilon_{1}}([t_{1}]+\tau)}-e^{A_{\epsilon_{1}}t_{1}})\\
& +\left(\displaystyle\frac{1}{1-\lambda_{fG}}\right)^{\left(\left[\displaystyle\frac{t}{\tau}\right]-\left[\displaystyle\frac{t_{1}}{\tau}\right]-1\right)}(e^{A_{\epsilon_{1}}([t_{1}]+2\tau)}-e^{A_{\epsilon_{1}}([t_{1}]+\tau)})\\
&+ ...+\left(\displaystyle\frac{1}{1-\lambda_{fG}}\right)(e^{A_{\epsilon_{1}}}-1)e^{A_{\epsilon_{1}}([t]-\tau)}\\
&  +e^{A_{\epsilon_{1}}t}-e^{A_{\epsilon_{1}}[t]}].\\
\end{align*}
Using the fact that $[t_{1}]\leq t_{1}$, it follows that
\begin{align*}
z(t)&\leq z(t_{1}^{+})\left(\displaystyle\frac{1}{1-\lambda_{fG}}\right)^{\left(\left[\displaystyle\frac{t}{\tau}\right]-\left[\displaystyle\frac{t_{1}}{\tau}\right]\right)}e^{(-A_{\epsilon_{1}}(t-t_{1}))}\\
&+\displaystyle\frac{\mu_{G}}{A_{\epsilon_{1}}}e^{(-A_{\epsilon_{1}}t)}[
\left(\displaystyle\frac{1}{1-\lambda_{fG}}\right)^{\left[\displaystyle\frac{t-t_{1}}{\tau}\right]}(e^{A_{\epsilon_{1}}([t_{1}]+\tau)}-e^{A_{\epsilon_{1}}[t_{1}]})\\
& +\left(\displaystyle\frac{1}{1-\lambda_{fG}}\right)^{\left(\left[\displaystyle\frac{t}{\tau}\right]-\left[\displaystyle\frac{t_{1}}{\tau}\right]-1\right)}(e^{A_{\epsilon_{1}}([t_{1}]+2\tau)}-e^{A_{\epsilon_{1}}([t_{1}]+\tau)})\\
&+ ...+\left(\displaystyle\frac{1}{1-\lambda_{fG}}\right)(e^{A_{\epsilon_{1}}}-1)e^{A_{\epsilon_{1}}([t]-\tau)}\\
&  +e^{A_{\epsilon_{1}}t}-e^{A_{\epsilon_{1}}[t]}]\\
&=z(t_{1}^{+})\left(\displaystyle\frac{1}{1-\lambda_{fG}}\right)^{\left(\left[\displaystyle\frac{t}{\tau}\right]-\left[\displaystyle\frac{t_{1}}{\tau}\right]\right)}e^{(-A_{\epsilon_{1}}(t-t_{1}))}\\
&+\displaystyle\frac{\mu_{G}}{A_{\epsilon_{1}}}(e^{A_{\epsilon_{1}}\tau}-1)e^{(-A_{\epsilon_{1}}t)}\left(\displaystyle\frac{1}{1-\lambda_{fG}}\right)^{\left[\displaystyle\frac{t}{\tau}\right]}[
\left((1-\lambda_{fG})e^{A_{\epsilon_{1}}\tau}\right)^{\left[\displaystyle\frac{t_{1}}{\tau}\right]}\\
& +\left((1-\lambda_{fG})e^{A_{\epsilon_{1}}\tau}\right)^{\left(\left[\displaystyle\frac{t_{1}}{\tau}\right]+1\right)}\\
&+ ...+\left((1-\lambda_{fG})e^{A_{\epsilon_{1}}\tau}\right)^{\left(\left[\displaystyle\frac{t}{\tau}\right]-1\right)}]\\
&  +\displaystyle\frac{\mu_{G}}{A_{\epsilon_{1}}}e^{(-A_{\epsilon_{1}}t)}[e^{A_{\epsilon_{1}}t}-e^{A_{\epsilon_{1}}[t]}]\\
\end{align*}

Let
\begin{equation}
\label{Impuleq95} a=(1-\lambda_{fG})e^{A_{\epsilon_{1}}\tau}.
\end{equation}
We have
\begin{align*}
S&=a^{\left[\displaystyle\frac{t_{1}}{\tau}\right]}+a^{\left(\left[\displaystyle\frac{t_{1}}{\tau}\right]+1\right)}+...+a^{\left(\left[\displaystyle\frac{t_{1}}{\tau}\right]+n\right)}\\
&=a^{\left[\displaystyle\frac{t_{1}}{\tau}\right]}\left\{1+a+a^{2}+...+a^{n}\right\}\\
&=a^{\left[\displaystyle\frac{t_{1}}{\tau}\right]}\times\frac{1-a^{n+1}}{1-a}.
\end{align*}
In the expression of $S$, we have
$\left[\displaystyle\frac{t_{1}}{\tau}\right]+n=\left[\displaystyle\frac{t}{\tau}\right]-1$,
which implies that

\begin{equation}
\label{Impuleq96}
n=\left[\displaystyle\frac{t}{\tau}\right]-\left[\displaystyle\frac{t_{1}}{\tau}\right]-1.
\end{equation}

Substituting (\ref{Impuleq96}) in the expression of $S$ leads to

\begin{equation}
\label{Impuleq97}
S=a^{\left[\displaystyle\frac{t_{1}}{\tau}\right]}\times\frac{1-a^{\left[\displaystyle\frac{t}{\tau}\right]-\left[\displaystyle\frac{t_{1}}{\tau}\right]}}{1-a}=\frac{a^{\left[\displaystyle\frac{t_{1}}{\tau}\right]}-a^{\left[\displaystyle\frac{t}{\tau}\right]}}{1-a}.
\end{equation}

Using (\ref{Impuleq95}) and (\ref{Impuleq97}), we have

\begin{align*}
z(t)&\leq
L(t)+\displaystyle\frac{\mu_{G}}{A_{\epsilon_{1}}}e^{(-A_{\epsilon_{1}}t)}
\left\{-\displaystyle\frac{(e^{A_{\epsilon_{1}}\tau}-1)\left(e^{A_{\epsilon_{1}}\tau}\right)^{\left[\displaystyle\frac{t}{\tau}\right]}}{1-(1-\lambda_{fG})e^{A_{\epsilon_{1}}\tau}}+e^{A_{\epsilon_{1}}t}-e^{A_{\epsilon_{1}}[t]}
\right\}\\
\end{align*}

where

\begin{align*}
\label{Impuleq98}
L(t)&=z(t_{1}^{+})\left(\displaystyle\frac{1}{1-\lambda_{fG}}\right)^{\left(\left[\displaystyle\frac{t}{\tau}\right]-\left[\displaystyle\frac{t_{1}}{\tau}\right]\right)}e^{(-A_{\epsilon_{1}}(t-t_{1}))}\\
&+\displaystyle\frac{\mu_{G}}{A_{\epsilon_{1}}}(e^{A_{\epsilon_{1}}\tau}-1)e^{(-A_{\epsilon_{1}}t)}\left(\displaystyle\frac{1}{1-\lambda_{fG}}\right)^{\left(\left[\displaystyle\frac{t}{\tau}\right]-\left[\displaystyle\frac{t_{1}}{\tau}\right]\right)}\left\{\displaystyle\frac{\left(e^{A_{\epsilon_{1}}\tau}\right)^{\left[\displaystyle\frac{t_{1}}{\tau}\right]}}{1-(1-\lambda_{fG})e^{A_{\epsilon_{1}}\tau}}
\right\}.\\
\end{align*}

Then we have for all $t>t_{1}$,
\begin{align*}
z(t) &\leq L(t)+\displaystyle\frac{\mu_{G}}{A_{\epsilon_{1}}}
\left\{1-\displaystyle\frac{(e^{A_{\epsilon_{1}}\tau}-1)\left(e^{A_{\epsilon_{1}}\tau}\right)^{\left(\left[\displaystyle\frac{t}{\tau}\right]-\displaystyle\frac{t}{\tau}\right)}}{1-(1-\lambda_{fG})e^{A_{\epsilon_{1}}\tau}}-\left(e^{A_{\epsilon_{1}}\tau}\right)^{\left(\left[\displaystyle\frac{t}{\tau}\right]-\displaystyle\frac{t}{\tau}\right)}  \right\}\\
\end{align*}

Since $\epsilon_{1}>0$ is arbitrary, it is obvious that
\begin{align*}
\lim\limits_{t\rightarrow\infty}L(t)&=\lim\limits_{t\rightarrow\infty}z(t_{1}^{+})\left(\displaystyle\frac{1}{1-\lambda_{fG}}\right)^{\left(\left[\displaystyle\frac{t}{\tau}\right]-\left[\displaystyle\frac{t_{1}}{\tau}\right]\right)}e^{(-A_{\epsilon_{1}}(t-t_{1}))}\\
&+\lim\limits_{t\rightarrow\infty}\displaystyle\frac{\mu_{G}}{A_{\epsilon_{1}}}(e^{A_{\epsilon_{1}}\tau}-1)e^{(-A_{\epsilon_{1}}t)}\left(\displaystyle\frac{1}{1-\lambda_{fG}}\right)^{\left(\left[\displaystyle\frac{t}{\tau}\right]-\left[\displaystyle\frac{t_{1}}{\tau}\right]\right)}\left\{\displaystyle\frac{\left(e^{A_{\epsilon_{1}}\tau}\right)^{\left[\displaystyle\frac{t_{1}}{\tau}\right]}}{1-(1-\lambda_{fG})e^{A_{\epsilon_{1}}\tau}}
                                \right\}\\
                                &=0.
\end{align*}

Thus, it follows that
\begin{equation}
\label{Impuleq99}
\lim\limits_{t\rightarrow\infty}z(t)\leq\displaystyle\frac{\mu_{G}}{A_{\epsilon_{1}}}.
\end{equation}

On the other hand  $\lim\limits_{\epsilon_{1}\rightarrow
0}A_{\epsilon_{1}}=\lim\limits_{\epsilon_{1}\rightarrow
0}(r_{G}-\gamma_{TG}\epsilon_{1})=r_{G}$. Thus  for
$\epsilon_{1}\rightarrow 0$, (\ref{Impuleq99}) changes to

\begin{equation}
\label{Impuleq101}
\lim\limits_{t\rightarrow\infty}z(t)\leq\displaystyle\frac{1}{X_{G}}.
\end{equation}

Because $z=\displaystyle\frac{1}{G}$, from (\ref{Impuleq101}) it
follows that
\begin{align*}
\lim\limits_{t\rightarrow\infty}\displaystyle\frac{1}{G(t)}&\leq\displaystyle\frac{1}{X_{G}},\\
\end{align*}

which implies that
\begin{equation}
\label{Impuleq102} \lim\limits_{t\rightarrow\infty}G(t)\geq X_{G}.
\end{equation}

Thus for any $\epsilon_{2}>0$, there exists $t_{2}>0$  such that
\begin{equation}
\label{Impuleq103} G(t)\geq X_{G}-\epsilon_{2}
\end{equation}
for all $t>t_{2}$.\par
On the other hand into $\mathcal{B}$, we have
\begin{equation}
\label{Impuleq104} G(t)\leq X_{G},
\end{equation}
for all $t>0$.
Coupling (\ref{Impuleq103}) and (\ref{Impuleq104}),
for all $t>t_{2}$ we obtain

\begin{equation}
\label{Impuleq105} X_{G}\geq G(t)\geq X_{G}-\epsilon_{2}
\end{equation}

Let $\epsilon=\min\{\epsilon_{1},\epsilon_{2}\}$, and
$t_{*}=\max\{t_{1},t_{2}\}$. Then from (\ref{Impuleq105}), we obtain
\begin{displaymath}
X_{G}\geq G(t)\geq X_{G}-\epsilon,
\end{displaymath}
for all $t>t_{*}$ that is

\begin{equation}
\label{Impuleq106} G(t)\rightarrow
X_{G},\hspace{0.5cm}\mbox{as}\hspace{0.5cm}t\rightarrow\infty,
\epsilon\rightarrow 0^{+}.
\end{equation}

It is obvious that
\begin{equation}
\label{Impuleq107} \lim\limits_{t\rightarrow\infty}\tilde{G}_{e}(t)=
X_{G}.
\end{equation}

Using (\ref{Impuleq106}) and (\ref{Impuleq107}), we have

\begin{displaymath}
\lim\limits_{t\rightarrow\infty}|G(t)-\tilde{G}_{e}(t)|=0.
\end{displaymath}
Therefore, the grassland savanna periodic equilibrium is globally
asymptotically stable when $\mathcal{R}_{0,pulse}^{*}<1$. This
completes  the proof of the theorem.

 \end{itemize}

\end{document}